\documentclass[final,3p]{elsarticle}
\usepackage{amssymb}
 \usepackage{amsthm}
\usepackage{amsmath,bbm,dsfont, color, xfrac, accents}
\usepackage{hyperref}

\usepackage[OT2, T1]{fontenc}

\usepackage[cp1250]{inputenc}

\begin{document}
\newcommand{\supp}       {{\rm supp\,} }
\newcommand{\csubset}         {\hookrightarrow}
\newcommand*{\id}{\mathds{1}}
\newcommand*{\zahlen}{\mathbb{Z}}
\newcommand*{\en}{\mathbb{N}}
\newcommand*{\er}{\mathbb{R}}
\newcommand{\divv}{\mathop{\text{div}\,}}
\newcommand{\esssup}{\mathop{\text{ess\:sup}\,}}
\newcommand{\htimes}{\mathop{\text{\large$ďż˝$}}}
\newcommand{\hexists}{\mathop{\text{\LARGE$\exists$}}}
\newcommand{\hforall}{\mathop{\text{\LARGE$\forall$}}}
\newtheorem{lem}{Lemma}
\newtheorem{theo}{Theorem}
\newtheorem{cor}{Corollary}
\newtheorem{prop}{Proposition}
\newtheorem{mydef}{Definition}
\newtheorem{rem}{Remark}
\newtheorem{ass}{Assumption}

\newcommand{\comentario}[1]{\textcolor{red}{\sc[*** #1 ***]}}

\begin{frontmatter}

 \title{Evolutionary, symmetric p-Laplacian. Interior regularity of time derivatives and its consequences}

\author{Jan Burczak} \address{Institute of Mathematics, Polish Academy of Sciences, \'Sniadeckich 8, 00-656 Warsaw}

\author{Petr Kaplick\'y }  \address{Department of Mathematical Analysis, Faculty of Mathematics and Physics, \\
Charles University in Prague,
186 75 Praha 8,
Czech Republic
}



\begin{abstract}
We consider the evolutionary symmetric $p$-Laplacian with safety $1$. By \emph{symmetric} we mean that the full gradient of $p$-Laplacian is replaced by its symmetric part, which causes breakdown of the Uhlenbeck structure. We derive the interior regularity of time derivatives of its local weak solution. To circumvent the space-time growth mismatch, we devise a new local regularity technique of iterations in Nikolskii-Bochner spaces. It is interesting by itself, as it may be modified to provide new regularity results for the full-gradient $p$-Laplacian case with lower-order dependencies. Finally, having the regularity result for  time derivatives,  we obtain respective regularity of the main part. The Appendix on Nikolskii-Bochner spaces, that includes theorems on their embeddings and interpolations, may be of independent interest. 
\end{abstract}




\end{frontmatter}






\newenvironment{rcases}
  {\left.\begin{aligned}}
  {\end{aligned}\right\rbrace}


%

%


\newcommand{\ep}{\varepsilon}
\newcommand{\eps} {\varepsilon}


\newcommand{\Div}       {{\rm div}_x}
\newcommand{\toc}       {{\stackrel{b}{\longrightarrow\,}}}

\newcommand{\A}     {{\mathcal A}}

\newcommand{\V}     {{\mathcal V}}

\newcommand{\D}     {{\rm D}}

\def\tens#1{\pmb{\mathsf{#1}}}
\def\vec#1{\boldsymbol{#1}}

\newcommand{\eN}	{{\mathcal{N}}} 
\newcommand{\M}	{{\mathcal{M}}} 
\newcommand{\eO}	{{\mathcal{O}}} 
\newcommand{\T}	{{\mathcal{T}}}
\newcommand{\W}	{\mathbb W} 
\newcommand{\ti}	{\tilde}
\newcommand{\wti}	{\widetilde} 
\newcommand{\vep}{\varepsilon}

\newcommand{\vdd}{\tens{\bar{v}}} 
\newcommand{\Ref} {\eqref}





%

\pagestyle{headings} 

%


%

%

\newcommand{\Sym}{{\rm Sym}}

\newcommand{\vfi}{{\varphi}}
\newcommand{\Vfi}       {\bar \varphi}

 \def\Xint#1{\mathchoice
 {\XXint\displaystyle\textstyle{#1}}%
 {\XXint\textstyle\scriptstyle{#1}}%
 {\XXint\scriptstyle\scriptscriptstyle{#1}}%
 {\XXint\scriptscriptstyle\scriptscriptstyle{#1}}%
 \!\int}
 \def\XXint#1#2#3{{\setbox0=\hbox{$#1{#2#3}{\int}$}
 \vcenter{\hbox{$#2#3$}}\kern-.5\wd0}}
 \def\ddashint{\Xint=}
 \def\dashint{\Xint-}

\def\B{\mathbb{B}}
\def\test{\mathcal{D}}
\def\F{\mathbb{F}}
\def\N{\mathbb{N}}
\def\O{\Omega}
\def\R{\mathbb{R}}
\def\Kor{{\rm Kor}}
\def\pod#1{\mathop{#1}\limits}
\def\diagin{-\hskip-11.0truept\intop}
\def\diagint{{\raise-.1pt\hbox{--}\hskip-7.9pt\intop}}
\def\diagintop{\mathop{\mathchoice
{{\diagin}}%
{{\diagint}}%
{{\diagint}}%
{{\diagint}}%
}\limits}
 \def\ddashint{\Xint=}
 \def\dashint{\Xint-}

%
%
%
%
\section{Introduction}
\noindent
Let $\D u$ denote the symmetric part of the gradient, \emph{\emph{i.e.}} $\D u = \frac{\nabla u + \nabla^T u }{2}$. We consider the following symmetric $p$-Laplace system
\begin{equation}\label{eq:jb_structure:weak}
u,_{t}-\,\divv \A(\D u)=0,
\end{equation}
that generalizes the generic examples
\begin{equation}\label{eq:jb_structure:generic}
\A^1= (1 + | \D u|)^{p-2} \D u, \qquad   \A^2 = (1 + | \D u|^2)^\frac{p-2}{2} \D u
\end{equation}
with $p \ge 2$. The precise assumptions on $\A$ are presented in Section \ref{seq:structure}.

 In this paper we obtain local regularity results connected with the second-order energy estimates in time (\emph{\emph{i.e.}} testing \eqref{eq:jb_structure:weak} with the localized $u,_{tt}$, roughly speaking). Contrary to the result of \cite{BurKap14a}, where we derive local regularity related to the second-order energy estimates in space, now we face serious difficulties, even on the level of formal estimates. These difficulties arise from the space-time mismatch in localization terms (for more on this, see subsection \ref{ssec:apriori}). We cope with these difficulties by means of a new interior iteration technique in Nikolskii-Bochner spaces,  inspired by Bul\'\i\v{c}ek, Ettwein, Kaplick\'y \& Pra\v{z}ak \cite{BEKP11}. We believe that our approach is interesting in itself and may be fruitful for obtaining new regularity results for other $p$-nonlinear systems of PDEs, for instance for the full $p$-Laplacian with lower-order terms ($\A(x,\nabla u)$ or  $\A(x, u, \D u)$). In addition to obtaining the quantitative result, we also provide (qualitative) inequalities.

\section{Motivation and known results}
From the perspective of  non-Newtonian hydrodynamics and nonlinear elasticity, it would be essential to repeat for the symmetric  $p$-Laplacian the $C^{1, \alpha}_{loc}$-regularity result, available for the full gradient $p$-Laplacian since the works of Uhlenbeck \cite{Uhl77}, Tolksdorff \cite{Tol83} (stationary case) and DiBenedetto \& Friedman  \cite{DiBFri85} (evolutionary case). Unfortunately, the pointwise structure of the symmetric  $p$-Laplacian seems to be resistant to  the methods used in the full $p$-Laplacian case\footnote{For more on this, see \cite{Burphd} and \cite{Bur14}} to get boundedness of gradients. 

Nevertheless, one can obtain certain regularity results for the symmetric  $p$-Laplacian. For the stationary case, Beir\~ao  da Veiga   Beir\~ao  da Veiga \& Crispo \cite{BdV} and \cite{daVCri13dcds} provide smoothness of a periodic-boundary value problem, provided $p$ is close to $2$. Results for generic boundary-value problems, developed for the full $p$-Navier-Stokes system, are of course available for the symmetric $p$-Laplacian. In particular, one has smoothness of solutions to basic initial-boundary value problems in $2$d case, see Kaplick\'y, M\'alek \& Star\'a \cite{KapMalSta02} and  Kaplick\'y \cite{Kap05}, \cite{Kap08} as well as existence of strong solutions for $3$d case, compare \cite{MalNecRuz01} by M\'alek, Ne\v cas \& R\r u\v zi\v cka and also \cite{Vei09a}, \cite{Vei09b} by Beir\~ao  da Veiga and \cite{VeiKapRuz11} by Beir\~ao  da Veiga,  Kaplick\'y \& R\r u\v zi\v cka. Let us mention also here a recent local regularity study for the $p$-Laplacian, symmetric $p$-Laplacian and $p$-Stokes type problems by Frehse \& Schwarzacher \cite{FreSchARXIV}. Since their results corresponds strongly to ours, we compare them in a more detailed manner in Subsection \ref{ssec:comp}. For small data regularity results, one may refer to Crispo \& Grisanti \cite{CriGri08}. 

As remarked, the regularity results cited above concern certain basic boundary-value problems. Local (interior) regularity results are much more scarce. In \cite{Bur14} the partial $C^{1, \alpha}_{loc}$ regularity theory has been developed. One should mention also \cite{FucSer} by Fuchs and Seregin. 

On the other hand, the regularity results for the full-gradient case are abundant. For the evolutionary  $p$-case, let us restrict ourselves to referring to the classical monograph by DiBenedetto \cite{DiB93} and a simple proof of $C^{1, \alpha}_{loc}$ regularity by Gianazza, Surnachev \& Vespri \cite{GiaSurVes10}.

\section{Plan of the paper}
In Section \ref{sec:prelim} we provide needed preliminaries, including notation, assumptions on growth of stress tensor $\A$,  definition of local weak solution to \eqref{eq:jb_structure:weak} and elements of Nikolskii-Bochner spaces. For traceability, details concerning Nikolskii-Bochner, including the equivalence between their definitions, their embeddings and interpolations, are gathered in Section \ref{sec:app} -- Appendix. Another reason for creating Appendix is the fact that its results may be of independent interest. Section \ref{sec:res} presents our main results concerning \eqref{eq:jb_structure:weak} and  Section \ref{sec:pfs} is devoted to their proofs.

\section{Preliminaries}\label{sec:prelim}

\subsection{Notation}
Constants denoted by $C, K$ may change from line to line of estimates and are larger than $1$. If a more careful control over a constant is needed, we denote its dependence on certain parameters writing $C(parameter)$ and generally suppress marking its dependence on irrelevant parameters. Such constant may also vary.

A space-time point $z=(x,t)$ is taken from $\Omega \times I=: \O_I$, where $\Omega$ is a spatial domain (an open, bounded connected set) in  $\er^d$ and $I = (I_L, I_R)$ is an open, bounded time interval in $\er$. As we develop the interior regularity theory, further assumptions on $\O_I $ are immaterial. $B_r (x)$ denotes the ball with the radius $r$ centered at a point $x$,  $I_\rho (t) $ --- the interval centered at $t$ and with radius $\rho$ and  $Q^r_{I_\rho} (z)$ denotes the space-time cylinder $B_r (x) \times I_\rho (t)$. We use, in particular, cylinders with parabolic scaling, \emph{i.e.} cylinders of the type $Q^r_{I_{r^2}} (z)$, which we denote briefly by $Q_r (z)$ and refer to as parabolic cylinders. Where possible, we drop dependence on $x, z, t$.

Symbol '$\Subset $' denotes embedding of the closure, \emph{i.e.} $A \Subset B$ iff $\overline{A} \subset B$. For a set  $S$ in an Euclidean space we denote its $\eps$-interior by '$S^\eps$'$ := \{ s \in S | \; dist(s, \partial S) \ge \eps \}$. 

We will write $ \Delta_h u := u \circ T_h - u,$ where $T_h f (x,t) :=  f (x,t+h)$.

We use the standard notation for function spaces, often dropping the underlying domain, when there is no danger of confusion. For example, we write
$L^q (W^{k,p}), W^{r,s} (W^{1,q}) $ for  Sobolev-Bochner spaces and
$N^{k,p} (I;X) $ for Nikolskii-Bochner space, see Section \ref{jb:strong:nbs}.

As we are interested in  systems of PDEs, the underlying spaces are vector valued. For instance $L^q (W^{k,p}) = L^q \left(I; \, (W^{k,p} (\O))^N \right) $. In fact we will have the equality of the dimension of the spatial domain $\O$ and of the target space, \emph{i.e.} $N=d$.

We will introduce additional notation where necessary.

\subsection{Structure of nonlinearity}\label{seq:structure}
We will work within

\begin{ass}\label{ass:struc} Fix arbitrary positive numerical constants $c,  C, K$. 
The tensor $\A$ satisfies for any symmetric ${d \times d}$ tensor  $P, Q$
\begin{equation*}
\begin{aligned}
(\A(P) - \A(Q) ) \! : \! (P- Q) & \ge  c\, \vfi'' (|P| + |Q|) |P-Q|^2, \\
| \A(P) - \A(Q)|& \le C\, \vfi'' (|P| + |Q|) |P-Q|.
\end{aligned}
\end{equation*}
where the function $\vfi \in C^{2} \left([ 0, \infty) \right)$ satisfies $\vfi'' (0) > 0$,  $\vfi' (0) = 0$,  $\vfi (0) = 0$ and
\[ K^{-1} \vfi'' (0) (1 + t^{p-2}) \le  \vfi'' (t) \le K (1 +  t^{p-2}) \] with $p \ge 2$.
\end{ass}
The prototype tensors $\A^1, \A^2$ satisfy Assumption \ref{ass:struc}. They are given through 
\[\A(Q) : = \partial_Q \vfi ( |Q|)\]
by the following $p$-potentials 
\begin{equation*}
\vfi^1 (t) = \int_0^t \! (\mu + s^{p-2}) \, s \,ds,  \qquad  \vfi^2 (t) = \int_0^t \!(\mu + s^2)^\frac{p-2}{2} \, s \,ds.
\end{equation*}
The formulation of Assumption \ref{ass:struc} is clearly inspired by the Orlicz-structure-type assumptions. In this context let us mention that in this paper we could have used Boyd indices $q_1, q_2$ of $\vfi$ instead of the pure $p$-growth. This would result in less tangible results and further technical complications, hence we do not follow this direction.

In order to gain a quadratic structure in our estimates, we introduce the following \emph{square root} of an $\eN$-function $\vfi$ and the associated tensor $\V$ 

\begin{mydef}[$\Vfi$ and $\V$] \label{def:jb_weak:V}
Let $\vfi$ be an $\eN$-function. We define 
\[\Vfi' := \sqrt{ t \vfi' (t)}, \qquad \V := \partial_Q \Vfi ( |Q|).\]
\end{mydef}
We have
\begin{equation}\label{eq:lem:jb_weak:mon_equiv2}
(\A(P) - \A(Q) ) \!:\! (P- Q)  \sim  \vfi'' (|P| + |Q|) |P-Q|^2  \sim | \V(P) - \V(Q) |^2,
\end{equation}
see for instance Proposition 5 of \cite{BurKap14a}.

Let us also make a technical assumption
\begin{ass}\label{rem:uberunter}
We fix now an \"ubercylinder $\ddot{Q}$ and a nonempty untercylinder $\underaccent{\ddot}{Q}$ such that $\underaccent{\ddot}{Q} \Subset \ddot{Q}$ and $ \ddot{Q} \Subset \O_I$, $ \ddot{Q} \subset Q_1$. All the subsequent analysis will happen in-between them. 
\end{ass}
Assumption \ref{rem:uberunter} serves only the purpose of not writing explicitly all the dependences of constants on the size of a domain.  Indeed, their scaling is now controlled by the choice of $\underaccent{\ddot}{Q}$, $\ddot{Q}$ and hence not written explicitly. As we derive in this paper local (interior) estimates, the important part of  Assumption \ref{rem:uberunter} lies in the untercylinder  $\underaccent{\ddot}{Q}$. We need this lower bound in the Poincar\'e-type estimates, for instance in Proposition \ref{lem:jb_strong:nikolski_ared}.

\subsection{Weak solutions}
We will use the following standard definition of the weak solution
\begin{mydef}\label{def:jb_weak:st_ws} Take an open interval $I$, a domain $\O \subset \er^d$ and a stress tensor $\A$ compatible with Assumption \ref{ass:struc}. A function $u \in C (I; L^2 (\O)) \cap  L^p (I; W^{1, p} (\O))$ is the {\emph {local weak solution of the system \eqref{eq:jb_structure:weak}} on $\O_I$\!} iff  for any subinterval $(t_1, t_2) \subset I$
\begin{equation}\label{eq:jb_weak:st_ws}
\int_\O u (t ) \cdot w (t){ \Big|}^{t_2}_{t_1} + \int^{t_2}_{t_1} \int_\O -\, u \cdot w,_t \, + \, \A (\D u): \D w = 0
\end{equation} 
for an arbitrary test function $w \in W^{1,1}_{0} (I; L^2 (\O)) \cap L^p (I; W^{1,p}_{0}(\O) )$.
\end{mydef}
The above definition implies that $u$ has the generalized time derivative in $ (L^p (I; W_0^{1, p} (\O)))^*$, compare Theorem 2.8 in \cite{Burphd}. Let us mention here, that instead of assuming $C (I; L^2 (\O))$ in Definition \ref{def:jb_weak:st_ws} we can start with $L^\infty (I; L^2 (\O))$ there and regain $C (I; L^2 (\O))$ via duality and just-mentioned regularity of the time derivative. For a time difference
\[
 \Delta_h u := u \circ T_h - u,
\]
where $T_h f (x,t) :=  f (x,t+h)$, we have
\begin{lem}\label{cor:jb_strong:intro} Take a local weak solution $u$  of \eqref{eq:jb_structure:weak} on $\O_I$.
Fix small $\delta > 0$ such that $I^\delta , \O^\delta $ are nonempty. For almost any real $\tau, \, | \tau| \le  \frac{\delta }{2}$ the generalized time derivative of the difference $\Delta_\tau u$ belongs to
\[
 \big(L^p (I^\delta ; W_0^{1, p} (\O^\delta )) \big)^*
\]
with the estimate
\begin{equation}\label{eq:jb_weak:lemdiff:2t}
 |  (\Delta_\tau u)_t   |_{ (L^p (I^\delta ; \,W_0^{1, p} (\O^\delta )))^*} \le C (p) \! \int_{\O_I}  \left( 1 +|\D u|^p  \right).
\end{equation}
For every $(t_1, t_2) \subset I^\delta$ and for an arbitrary test function ${w \in L^p (I^\delta ; W^{1,p}_{0} (\Omega^\delta  ))} $
\begin{equation}\label{eq:jb_weak:lemdiff:pointwise:t}
\int^{t_2}_{t_1} \left( \langle (\Delta_\tau u)_t , w \rangle_{(W^{1,p}_{0} (\Omega^\delta ))^*\!, \, W^{1,p}_{0} (\Omega^\delta  )} + \int_{\Omega^\delta } \Delta_\tau \!\left( \A ( \D u ) \right)  \!:\!  \D w  \right)=0.
\end{equation}
For almost every  $t \in I^\delta$ and for an arbitrary test function $v \in  W^{1,p}_{0} (\Omega^\delta ) $ we have
\begin{equation}\label{eq:jb_weak:theo:pointwise_t:t}
\langle  (\Delta_\tau  u)_t (t) , v \rangle_{(W^{1,p}_{0} (\Omega^\delta  ))^*\!, W^{1,p}_{0} (\Omega^\delta )}   + \int_{\Omega^\delta }  \Delta_\tau \! \left(  \A ( \D u (t) ) \right)  \!:\! \D v  =0
\end{equation}
for almost every $\tau$.
\end{lem}
We will need this result to derive rigorous estimates in Section \ref{sec:res}. It is shown as Lemma 2.9 in \cite{Burphd}.

\subsection{Nikolskii-Bochner spaces}\label{jb:strong:nbs}
As already mentioned, due to technical reasons, explained in Subsection \ref{ssec:apriori} we resort in proofs to Nikolskii-Bochner spaces. We provide here only their definition, as some of our results use them. The broader presentation of Nikolskii-Bochner spaces is shifted to Section \ref{sec:app} - Appendix for the sake of clarity of our exposition.

Take a function in  $f \in L^{1} ( I; X)$ and an $r \in 2, 3, \dots$. We define for any \[t \in  I_{rh} = \{ t \in I |\; t+rh \in I \}\]  the $r$-th order difference recursively by
\[
 \Delta^r_h f (t) :=  \Delta_h (\Delta^{r-1}_h f (t)).
 \]
 Let us introduce for any $\alpha \in \R_{+}$, natural $r > \alpha$ and $\delta \in (0, \infty)$
\begin{equation}\label{eq:jb_strong:nik_semi:gen}
[f]_{{r, \delta,} N^{\alpha, p} ( I; X)} := \sup_{h \in ( 0, \delta) } h^{- \alpha} | \Delta_h^r f |_{L^{p} ( I_{rh}; X)}.
\end{equation}
For $h$ exceeding $\frac{|I|}{r}$ the underlying time interval $I_{rh} $ becomes empty, so for any  $\delta \in (0, \infty)$
\begin{equation*}
[f]_{{r, \delta,} N^{\alpha, p} ( I; X)} \le [f]_{{r, \frac{|I|}{r},} N^{\alpha, p} ( I; X)}.
\end{equation*}
Let $[\alpha]$ be the largest natural number not larger than $\alpha$, \emph{i.e.}
\[
[\alpha] := \sup \{n \in \N | \; n \le \alpha \}.
\]
We introduce
\begin{mydef}\label{def:nik}
Let $X$ be a Banach space. Choose $\alpha \in \er_+$, $p \in [1, \infty]$ and the smallest integer $r_0> [\alpha]$, \emph{i.e.}
\[
r_0 =  [\alpha] +1.
\]
The {\emph{Nikolskii-Bochner space  $N^{\alpha, p} ( I; X)$}}  is the subspace of the Lebesgue-Bochner space $L^{p} ( I; X)$ with finite $[f]_{{r_0, 1,} N^{\alpha, p} ( I; X)}$,
 \emph{i.e.}
\begin{equation*}
N^{\alpha, p} ( I; X) := \{ f \in L^{p} ( I; X) | \; \;  [f]_{{r_0, 1,} N^{\alpha, p} ( I; X)}  < \infty \}.
\end{equation*}
\end{mydef}
Observe that in the case $\alpha \in \N$ we have $r_0= \alpha + 1$, \emph{i.e.} we use higher-order differences than these giving Sobolev-Bochner space.

\section{Results}\label{sec:res}
We will use 
\[
\gamma_0 := \frac{2p}{(p-2)(d (p-2) + p)}.
\]
Observe that 
\[p \ge 2 + \frac{2}{\sqrt{d+1}} \; \; \iff \gamma_0 \le 1.\]
Let us also define
\[
\gamma_1 := \frac{2}{p} + \frac{p}{d (p-2) + 2p}.
\]
We denote by $2^*$ either the Sobolev embedding exponent $\frac{2d}{d-2}$ for $d \ge 3$ or any finite number $q$ for $d=2$.

\subsection{A lightweight version of the main results}\label{sec:jb_strong:pcase}
 For the sake of clarity we start with qualitative results, \emph{\emph{i.e.}} we restrict ourselves to indicating what regularity class solutions belong to, without presenting the inequalities.

First we cover the case $\gamma_0 \le 1$, where we get low-regularity in time.
\begin{lem}\label{for:jb_strong:low}
Let Assumption \ref{rem:uberunter} be valid.  Consider a local weak solution $u$ to \eqref{eq:jb_structure:weak} on $\O_I$. If $p \ge 2 + \frac{2}{\sqrt{d+1}}$, then for any $\alpha < \gamma_0$, locally on $\O_I$, $u$ belongs to
\[
N^{\alpha, \infty} (L^2), \; N^{1+ \alpha, p'}( W^{-1, p'}),  \; N^{\alpha, 2} ( W^{1,2} ),  \; N^{\alpha, p} (L^p),  \;  N^{\frac{2\alpha}{p}, p} (W^{1,p} ),
\]
whereas $ \V (\D u)$ is in
\[
 N^{\alpha, 2} (L^2).
\]
\end{lem}

Observe that for  $p \in [2, 2 + \frac{2}{\sqrt{d+1}})$  we have $\gamma_1 \in (1, \frac{3}{2}]$. In this range of $p's$ we get the following better regularity. 

 \begin{lem}\label{for:jb_strong:hi}
Let Assumption \ref{rem:uberunter} be valid. Consider a local weak solution $u$ to \eqref{eq:jb_structure:weak} on $\O_I$. If $p \in [2, 2 + \frac{2}{\sqrt{d+1}})$, then for any positive
$ \gamma < \gamma_1$,
 locally on $\O_I$,  $u$ belongs to
\[
N^{ \gamma, p} (L^{p}), \; { W^{2, p'} (W^{-1, p'} )}, \; { W^{1, \infty} (L^2)}, \; {W^{1, 2} (W^{1,2})}, \; L^\infty ( W^{2,2}), \; L^\infty ( W^{1,\frac{2^*p}{2} }),
\]
and for $p>2$ additionally to
\[
N^{\frac{2}{p}, p} ( W^{1,p}).
  \]
Moreover, $ \V (\D u)$ is in
\[
{W^{1, 2} (L^2)},  \; L^\infty ( W^{1,2}).
\]
\end{lem}

 For $d=2$ Lemma \ref{for:jb_strong:hi} implies that locally $\nabla u \in L^\infty ( BMO)$. In fact it can be raised to H\"older continuity of $\nabla u$. It is the subject of current research, based on quadratic approximations of $\A$.

 \begin{cor}\label{cor:uc}
 Let Assumption \ref{rem:uberunter} be valid. Consider a local weak solution $u$ to \eqref{eq:jb_structure:weak} on $\O_I$. It is H\"older continuous for $p \in [2,4)$ for $d=2$ and  $p \in [2, 3+ \frac{ \sqrt{33}-3}{4})$ for $d=3$.
 \end{cor}
 \begin{proof}
 $N^{\frac{2\alpha}{p}, p} (W^{1,p} )$ with $\alpha < \gamma_0$  of Lemma  \ref{for:jb_strong:low} gives H\"older continuity of $u$, provided $p > d$, $\alpha> \sfrac{1}{2}$, which is equivalent to $p \in (d, \frac{3+2d+\sqrt{8d+9}}{d+1}) $. Focusing on the physically relevant $d =2, 3$ we see that  the case $p \in (2, 3)$ has  already been covered by Lemma  \ref{for:jb_strong:hi} and that we can easily complete the case $p=3$. Namely for $p=3$ we use another information of Lemma  \ref{for:jb_strong:low}, \emph{i.e.} $L^p ( W^{1,\frac{2^*p}{2} })$ to raise infinitesimally in space, by interpolation, the already used $N^{\frac{2\alpha}{p}, p} (W^{1,p} )$.
\end{proof}

\subsection{Full main results}

\begin{theo}\label{theo:jb_strong:temporal}

Take an untercylinder and an  \"ubercylinder $ \underaccent{\ddot}{Q} \Subset \ddot{Q} \Subset \O_I$ fixed by Assumption \ref{rem:uberunter} and a local weak solution $u$ to \eqref{eq:jb_structure:weak} on $\O_I$. \\
(i) \emph{Case  $p \ge 2 + \sfrac{2}{\sqrt{d+1}}$}. For any $\alpha < \gamma_0$ there exist $\kappa (p, \alpha), \, C (p, \alpha, \vfi'' (0))$ that for all concentric cylinders $Q_r, Q_R$ such that $\underaccent{\ddot}{Q} \subset Q_r  \Subset Q_R \subset \ddot{Q}$ one has
  \begin{multline*}
 |u|_{ N^{\alpha, \infty} (I_{r^2};  L^2 (B_{r})) }  + |\V (\D u)|_{N^{\alpha, 2} (I_{r^2}; L^2  (B_{r}))}  +
 |u |_{N^{1+ \alpha, p'} ( I_{r^2};  W^{-1, p'}(B_{r}) )} +  |u |_{N^{\alpha, 2} ( I_{r^2}; W^{1,2} (B_{r}))} +\\
   |u |_{ N^{\alpha, p} (I_{r^2}; L^{p} (B_{r}))} +  |u |_{ N^{\frac{2\alpha}{p}, p} (I_{r^2}; W^{1,p} (B_{r}))} \le
 C (p, \alpha, \vfi'' (0))\left( \frac{ 1+  |u|_{L^\infty ( I_{R^2};  L^{2} (B_R))} +  |u |_{L^{ p} ( I_{R^2}; W^{1, p}(B_R) )} }{R-r}  \right)^{\kappa (p, \alpha)}.
 \end{multline*}
(ii) \emph{Case  $p \in \left(2, 2 + \sfrac{2}{\sqrt{d+1}}\right)$}.
For any $\gamma$ such that
\[ \gamma < \gamma_1 \, \left( = \frac{2}{p} + \frac{p}{d (p-2) + 2p} \right)\]
there exist $\kappa (p, \gamma), \; C (p, \gamma, \vfi'' (0))$ such that for all the~concentric cylinders $Q_r, Q_R$ satisfying $\underaccent{\ddot}{Q} \subset Q_r  \Subset Q_R \subset \ddot{Q}$ one has
 \begin{multline*}
 |u|_{ N^{\gamma, p} ( I_{r^2}; L^{p}(B_{r}))}  + |u|_{ W^{2, p'} \left( I_{r^2}; W^{-1, p'} (B_{r}) \right)} +  |u|_{ W^{1, \infty} (I_{r^2}; L^2  (B_{r}))}  + |\V (\D u)|_{W^{1, 2} (I_{r^2}; L^2 (B_{r}))} + \\
 |u |_{W^{1, 2} (I_{r^2}; W^{1,2}  (B_{r}))} +
  |u|_{ N^{\frac{2}{p}, p} (I_{r^2}; W^{1,p}(B_{r}))} \le
C (p, \gamma, \vfi'' (0))  \left( \frac{ 1+  |u|_{L^\infty ( I_{R^2}; L^{2} (B_R))} +  |u |_{L^{ p} (I_{R^2}; W^{1, p} (B_R)) } }{R-r}  \right)^{\kappa(p, \gamma)}.
 \end{multline*}
(iii) \emph{Case $p=2$}.  Then for any positive $ \gamma < \frac{3}{2}$ there exist $ \kappa (\gamma), \, C (\gamma, p, \vfi'' (0))$ that for all the concentric cylinders $Q_r, Q_R$ satisfying $\underaccent{\ddot}{Q} \subset Q_r  \Subset Q_R \subset \ddot{Q}$ one has
  \begin{multline*}
 |u|_{ N^{ \gamma, 2} ( I_{r^2}; L^{2}(B_{r}))}  + |u|_{ W^{2, 2} ( I_{r^2}; W^{-1, 2} (B_{r}) )} +  |u|_{ W^{1, \infty} (I_{r^2}; L^2  (B_{r}))}  + |\V (\D u)|_{W^{1, 2} (I_{r^2}; L^2 (B_{r}))} +\\
  |u |_{W^{1, 2} (I_{r^2}; W^{1,2}  (B_{r}))} \le
 C (p, \gamma,  \vfi'' (0))  \left( \frac{  1+ |u|_{L^\infty ( I_{R^2};  L^{2} (B_R))} +  |u |_{L^{ 2} ( I_{R^2}; W^{1, 2} (B_R)) } }{R-r}  \right)^{\kappa( \gamma)}.
 \end{multline*}
\end{theo}

Theorem \ref{theo:jb_strong:temporal} can be extended to systems dependent on lower order terms $z, u$. More importantly, it can be a starting point to proving temporal interior regularity results for evolutionary $p$-Laplace, $p$-Stokes and $p$-Navier-Stokes systems, that are apparently strongly needed (see for instance Acerbi, Mingione \& Seregin \cite{AceMinSer04}, Bae \& Jin \cite{BumBae_pre} and Duzaar, Mingione \& Steffen \cite{DuzMinSte11}). Of course there, the main additional difficulty would be pressure.

Now, let us present the eponymous consequences of regularity of time derivatives.
Having  results of  Theorem \ref{theo:jb_strong:temporal} (ii) and (iii), let us consider system  \eqref{eq:jb_structure:weak}
as an inhomogenous  stationary one on time levels with its r.h.s. in  $ L^2$, \emph{i.e.}
\[
-\divv \A(\D u) (t)= u,_{t} (t)
\]
and derive the regularity of its l.h.s. from the information on $u,_{t} (t)$. It is a common approach in evolutionary PDEs that proves fruitful also in the nonlinear case --- recall the paper of Ne\v cas and \v Sver\'ak \cite{NecSve91}, where they study the regularity of solutions to general nonlinear, quadratic, strongly elliptic systems. In our case, the stationary estimate reads
\begin{theo}[Regularity via the stationary estimates]\label{lem:jb_strong_stat}
Take a local weak solution $u$ to \eqref{eq:jb_structure:weak} on $\O_I$.
If $ u,_t  \in L^\infty (I; L^2(\O))$ and   $\vfi  ( |\nabla u|)  \in L^\infty (I; L^1(\O))$, then for any concentric balls $B_r \Subset B_R \Subset \O$ 
\begin{multline}\label{eq:jb_strong:stat}
\esssup_{t \in I} \int_{B_r} \left( |\nabla \V ( \D  u (t) ) |^2   + \vfi'' (0)   | \nabla \D  u (t) |^2 \right) \le\\
\Big(1 + \frac{1}{\vfi'' (0)} \Big) \frac{C ( p) }{(R-r)^2}  \esssup_{t \in I}   \int_{B_R} \left(  \vfi  ( |\nabla u  (t)|)  +    |u,_t (t)|^2   \right).
\end{multline}
\end{theo}
\begin{rem}
The statement of the main theorems, especially Theorem  \ref{theo:jb_strong:temporal},  may seem discouraging. Let us clarify that they are simply a qualitative version of Lemmas from Subsection \ref{sec:jb_strong:pcase}.
\end{rem}
\subsection{Comparison with the available techniques}\label{ssec:comp}
\subsubsection{Comparison with results by parabolic embedding}
 Let us explain here, why we believe that our approach may be useful for obtaining new results for the full-gradient case with lower-order terms in the main part $\A (z, u, \nabla u)$. The standard technique in such a case is to rely on parabolic embedding, compare for instance Acerbi, Mingione \& Seregin \cite{AceMinSer04}  and  B\"ogelein, Duzaar \& Mingione \cite{BogDuzMin13}. For some more details of it, see Subsection \ref{ssec:apriori}. This approach has certain limitations: a lot of smoothness is needed a priori and, more importantly, one can deal only with $p \le 2 + \frac{4}{d} $ (possibly $+ \delta$ thanks to a Gehring-type argument). Observe that the range of $p$'s, where we obtain full time derivatives, is larger for $d>4$ than the range given by parabolic embedding approach (compare Lemma \ref{for:jb_strong:hi}). Moreover, unlike in the parabolic embedding approach, we can obtain some gain in regularity for \emph{any} $p$, see Lemma  \ref{for:jb_strong:low}. It implies, in particular, H\"older continuity of $u$ for  $p \in [2, \frac{3+2d+\sqrt{8d+9}}{d+1}), \; $ \emph{\emph{i.e.}} $p \in [2,4)$ for $d=2$, $p \in [2, 3+ \frac{ \sqrt{33}-3}{4})$ for $d=3$, see Corollary \ref{cor:uc}. This second range is larger than $[2, \frac{10}{3}]$ of parabolic embedding. Finally, one can improve our iterative scheme and allow for larger range of $p$'s in Lemma \ref{for:jb_strong:hi}, whereas bound $p \le 2 + \frac{4}{d}  \, (+ \delta)$ in the other technique seems to be unavoidable. This improvement of our iterative scheme involves using higher space regularity, provided by \cite{BurKap14a},  at each iteration step.
\subsubsection{Comparison with Frehse \& Schwarzacher \cite{FreSchARXIV}}
The main difference between our results and these of \cite{FreSchARXIV} concerns localization. We derive space-time local estimates and actually the biggest difficulty that we face, compare Subsection \ref{ssec:apriori}, follows from space-time mismatch in lower-order terms that stem from using a cutoff function. On the contrary, \cite{FreSchARXIV} concerns global in space estimates, hence no need of a cutoff function there (localization in time there follows simply from integration over a portion of time interval).

In fact there is no intersection of our result with \cite{FreSchARXIV}. While one of our main goals is to overcome the troubles with localization and to prove that $u_t\in L^\infty(L^2)$, Frehse \& Schwarzacher consider this regularity as starting point, see \cite[Proposition 4.1]{FreSchARXIV}.

In this context it is interesting that from our results in Lemmas \ref{for:jb_strong:low} and \ref{for:jb_strong:hi} it follows that $u\in N^{\alpha,2}(L^2)$ locally on $\Omega_I$ for any $\alpha\in[0,1+1/p)$ if $p\in[2,2+2/\sqrt{d+1})$, and $\alpha\in[0,1/2+\gamma_0/2+\gamma_0/p)$ if $p\geq2+2/\sqrt{d+1}$. Which is similar to $u\in N^{3/2,2}(L^2)$ obtained in \cite{FreSchARXIV}. Since our method is completely different to the one in \cite{FreSchARXIV}, it would be interesting to find if their method could improve the result presented here.
 
 
Finally, let us observe that in \cite{FreSchARXIV} the authors do not need to set any restriction on the growth $p$, do not need to rely on structure with safety-$1$ and they provide some results also for real fluid-dynamics systems.

\section{Proofs}\label{sec:pfs}
Let us first indicate the main trouble, by trying to derive formal estimates. Compare subsection 5.1 of our paper on regularity of space derivatives \cite{BurKap14a}.
\subsection{A priori estimates}\label{ssec:apriori}
 In order to detect the difficulties that are of a structural origin, we begin with a priori estimates. Let us test formally \eqref{eq:jb_structure:weak} with $( u,_{t} \psi),_{t}$, where $\psi$ is a cutoff function for $\O'_{I'} \Subset \O_I$.We arrive at
\begin{equation}\label{eq:jb_strong:apriori_3ta}
\int_{\O} |   u,_\tau   \psi  |^2 (\tau)+  \int_{0}^\tau \int_{\O}  \vfi'' (|  \D u|) |  \D   u,_t \!|^2  \psi^2 \le
C \cdot  (|  \psi,_t |_\infty   + |\nabla \psi |^2_\infty)  \int_{\O_I}  |    u,_t  \!|^2 +  \vfi'' (|  \D u|)   |   u,_t \!|^2 
\end{equation}
for any $\tau \in I'$. Observe, that the term $ \vfi'' (|  \D u|)   | u,_t\! |^2 $ at the r.h.s. of \eqref{eq:jb_strong:apriori_3ta} is more troublesome than the respective $ \vfi'' (|  \D u|)   | \nabla u |^2 $  of (21) in  \cite{BurKap14a}. This difficulty arises   especially in the case of an unbounded $\vfi''$, \emph{i.e.}  in the superquadratic regime, which is of interest for us. Using the assumed here $p$ growth of $\A$ we have $ \vfi'' (t) \le C (1 +  t^{p-2}) $. Hence one can choose $r>2$ and split $ \vfi'' (|  \D u|)   | u,_t \!|^2 $  into
\begin{equation}\label{eq:jb_strong:apriori_3tb}
\delta | u,_t \!|^r + C( \delta) \big( 1+ |  \D u|^\frac{r(p-2)}{r-2} \big).
\end{equation}
After integration over $\O_I$, the second summand \eqref{eq:jb_strong:apriori_3tb} may be externally controlled with the use of the spatial estimate of  in  \cite{BurKap14a} and the parabolic embedding 
\[|\D u|_{L^{2 + \frac{4}{d}} (\O_I)} \le  |\D u|^\theta_{L^\infty (L^2 (\O))}  |\D u|^{1-\theta}_{L^p (W^{1,p} (\O))},\] provided $\frac{r(p-2)}{r-2} \le 2 + \frac{4}{d}$. Consequently we would end up with
\begin{equation}\label{eq:jb_strong:apriori_3t}
\sup_{t \in I'} \int_{\O'} |  u,_t \! |^2 (t)+  \int_{I'} \int_{\O'}  \vfi'' (|  \D u|) |  \D u,_t \!|^2 \le 
C (|u|_V)  (|  \psi,_t |_\infty   + |\nabla \psi |^2_\infty)  \int_{\O_I} 1 +  \delta | u,_t \!|^r,
\end{equation}
where $V = L^\infty (I_1; L ^2 (\O_1)) \cap L^p (I_1;W^{1,p}(\O_1))$.
The entire information from the l.h.s. of \eqref{eq:jb_strong:apriori_3t} controls through the parabolic embedding $ u,_t  $ in $L^{2 + \frac{4}{d}} $. Hence we can close the estimate  for certain $p$'s, dealing with the larger support of the  r.h.s. of \eqref{eq:jb_strong:apriori_3t} by means of a Giaquinta-Modica-type Lemma, where we use the smallness of $\delta$. This is essentially the parabolic embedding method, mentioned in Subsection \ref{ssec:comp}. Such method has been used in \cite{AceMinSer04} by Acerbi, Mingione and Seregin to derive the spatial regularity of a full-gradient $p(x)$-Laplacian. Unfortunately, this approach requires very smooth approximate solutions to the considered system, in order to be made rigorous. By superseding the physically justified symmetric gradient with the full one, the authors gain such a natural smooth approximate solution, \emph{i.e.} the safety-$1$ $p$-Laplacian with a~smooth dependence on lower-order terms. One can also provide enough smoothing for the symmetric-gradient case, for instance by adding the extra smoothening term  $- \eps \Delta^m u$ to the system \eqref{eq:jb_structure:weak}, with a large  $m$. It is a rather technical approach and it is still limited to $p \le 2 + \frac{4}{d} (+ \delta)$ due to the parabolic embedding.

We will evade the need of the above mentioned use of smooth approximations to  \eqref{eq:jb_structure:weak} and another technicalities with a method inspired by \cite{BEKP11} of Bul\'\i\v cek, Ettwein, Kaplick\'y and Pra\v z\'ak. 

The key observation is that we already possess a certain time differentiability (of a fractional degree) directly from a weak solution. Observe that in order to obtain the information on the time continuity of weak  $L^\infty  (I; L^2 (\O)) \cap L^p (I; W_0^{1, p} (\O))$ solutions one interpolates, roughly speaking, $u_t \in (L^p (I; W_0^{1, p} (\O)))^*$ and $u \in L^p (I; W_0^{1, p} (\O))$. However, one can interpolate between the same spaces in the scale of Besov spaces and obtain certain information on fractional time derivatives of $u$. Next we will try to  improve this information iteratively.

 \subsection{Energy estimates  in Nikolskii-Bochner spaces}\label{sec:energyinNik}
 In the next three sections we derive  regularity results in Nikolskii-Bochner spaces. We strongly rely on the assumption of the $p$-structure for the system \eqref{eq:jb_structure:weak}, in order to describe  precisely the Nikolskii-Bochner space resulting from a given estimate. As already mentioned, it is also possible to obtain analogous results for a~growth function $\vfi$ satisfying certain growth restrictions expressed via Boyd indices, but then results are less tangible, compare Remark \ref{rem:notp}.
 
 Recall that $Q_\varrho$ is the parabolic cylinder $B_\varrho \times I_{\varrho^2}$, $Q^\delta := \{ z \in Q \; | \; z+s \in Q \; \text{ for any } |s| \le \delta \}$ and for an interval $I$ we denote $I_{rh} = \{ t \in I |\; t+rh \in I \} $. Recall also that Definition \ref{def:jb_weak:st_ws} of a local weak solution to \eqref{eq:jb_structure:weak} already includes the fact that $\A$ satisfies  growth Assumption \ref{ass:struc}.
\subsubsection{High time-regularity estimates in Nikolskii-Bochner spaces}
The result of this subsection will be used in the following Section \ref{sec:jb_strong:iteration} as a~high-time and low-space regularity ingredient of interpolation Lemma \ref{lem:jb_strong:nikolski_inter}, \emph{\emph{i.e.}} a  ${N^{\alpha_2, p_2} ( I; W_0^{-1, q_2} (Q)) } $ part of the r.h.s. of \eqref{eq:jb_strong:nikolski_inter}. 
\begin{lem}\label{lem:jb_strong:iter_1} We take an untercylinder and an \"ubercylinder $\underaccent{\ddot}{Q} \Subset \ddot{Q} \Subset \O_I$ fixed by Assumption \ref{rem:uberunter}.
 Fix any  $\delta >0$ so small that $\ddot{Q}\subset \O_I^\delta$ and any smooth time-independent $\eta \in C^\infty (\O)$. Take a local weak solution $u$ to \eqref{eq:jb_structure:weak} on $\O_I$. Then for any $Q_\rho$ such that $ \underaccent{\ddot}{Q} \subset Q_\rho \subset \ddot{Q} $ and any $\alpha \in [0,1) $
\begin{multline}\label{eq:jb_strong:iter_1:superquad}
|u \eta|_{\frac{\delta}{2}, N^{1+ \alpha, p'} (I_{\rho^2}; W^{-1, p'} (B_\rho)) } \le \\
   C (p)  ( |\eta|_\infty + | \nabla \eta|_\infty ) 
[\V (\D u)]_{{1, \frac{\delta}{2},} N^{\alpha, 2} ( I_{\rho^2};  L^2( B_\rho))}    \left( 1+  \left|  \D u \right|_{ L^p (\O_I)}  \right)^{\frac{p-2}{2}} + C |u \eta|_{L^\infty (I_{\rho^2}; L^{2} (B_\rho))} 
\end{multline}
\begin{multline}\label{eq:jb_strong:iter_1:superquad:sob}
|u \eta|_{ W^{2, p'} (I_{\rho^2}; W^{-1, p'} (B_\rho)) } \le \\
   \frac{C (p)}{\delta}  ( |\eta|_\infty + | \nabla \eta|_\infty ) 
[\V (\D u)]_{{1, \frac{\delta}{2},} N^{1, 2} ( I_{\rho^2};  L^2( B_\rho))}    \left( 1+  \left|  \D u \right|_{ L^p (\O_I)}  \right)^{\frac{p-2}{2}} + C |u \eta|_{L^\infty (I_{\rho^2}; L^{2} (B_\rho))} +   \\
(  |  \eta|_\infty + | \nabla \eta|_\infty )^{p-1}    \left|  u \right|^{p-1}_{ L^p (I; W^{1,p} (\Omega))} 
\end{multline}
 provided the above r.h.s.'s are meaningful.  For the case $p=2$ we apply the convention {\rm>>}$\infty^0 = 1${\rm <<} (In the sense that the term $ ( 1+  \left|  \D u \right|_{ L^p (\O_I)} )^{\frac{p-2}{2}}$ in the  r.h.s.'s above is not present for $p=2$.).
\end{lem}
\begin{proof}
As $\ddot{Q} \subset \O_I^\delta$, we have $Q_\rho \subset  \O_I^\delta$. For an arbitrary ${w \in  L^p( (I_{\rho^2})_h; W^{1,p}_{0} (B_\rho) )}$, let us use \eqref{eq:jb_weak:lemdiff:pointwise:t} of Lemma \ref{cor:jb_strong:intro} with $\tau : = h$, time interval $(t_1, t_2 ) := (I_{\rho^2})_h$ and $w := w \eta$. We get for a fixed $\alpha \in (0,1]$ and  any $h \in (0, \sfrac{\delta}{2}]$
\begin{equation}\label{jb_strong:t:bib:1}
h^{-\alpha} \int_{ (I_{\rho^2})_h} \left\langle ( \Delta^h u )_t , w \eta \right\rangle_{(W^{1,p}_{0} (B_\rho))^*\!,\, W^{1,p}_{0} (B_\rho  )}  = - h^{-\alpha} \int_{{(I_{\rho^2})_h \times B_\rho}} \Delta^h  \left(  \A ( \D u ) \right) : \D (w \eta).
\end{equation}
Thanks to independence of $\eta$ from the time variable we can put $\eta$ inside the generalized time derivative above. Hence  \eqref{jb_strong:t:bib:1} gives
\begin{equation}\label{jb_strong:t:bib:2}
h^{-\alpha} \int_{  (I_{\rho^2})_h } \left\langle ( \Delta^h u \eta )_t , w \right\rangle_{(W^{1,p}_{0} (B_\rho))^*\!,\, W^{1,p}_{0} (B_\rho  )}  = - h^{-\alpha} \int_{{(I_{\rho^2})_h \times B_\rho}} \Delta^h  \left(  \A ( \D u ) \right) : \D (w \eta).
\end{equation}
Now we estimate the r.h.s. of \eqref{jb_strong:t:bib:2} using Assumption \ref{ass:struc} and the H\"older inequality
\begin{multline*}
\Big| h^{-\alpha} \int_{{(I_{\rho^2})_h \times B_\rho}} \Delta^h  \left(  \A ( \D u) \right) : \D (w \eta) \Big| \le  h^{-\alpha}  \int_{{(I_{\rho^2})_h \times B_\rho}}   \vfi''({|\D u| + |\D u \circ T_h| }) |\Delta^h \D u|  | \D (w \eta) | \le \\
\bigg| \frac{ \sqrt{ \vfi''({|\D u | + |\D u \circ T_h| }) }|\Delta^h \D u| }{ h^{\alpha}} \bigg|_{L^2 ({(I_{\rho^2})_h \times B_\rho} )} \left|  \sqrt{ \vfi''({|\D u | + |\D u \circ T_h| }) }|  \D (w \eta)| \right|_{L^2 ({(I_{\rho^2})_h \times B_\rho} )}. 
\end{multline*}
Further we estimate the first term of the r.h.s. above by \eqref{eq:lem:jb_weak:mon_equiv2}. This, plugged into \eqref{jb_strong:t:bib:2}, gives
\begin{multline}\label{jb_strong:t:bib:3}
h^{-\alpha} \Big|  \int_{  (I_{\rho^2})_h } \langle ( \Delta^h u \eta )_t , w  \rangle_{(W^{1,p}_{0} (B_\rho))^*\!, \, W^{1,p}_{0} (B_\rho  )}   \Big|  \le \\
 C (p)
 \Big| \frac{\Delta^h \V (\D u)}{h^\alpha} \Big|_{ L^2 ({(I_{\rho^2})_h \times B_\rho} ) }  \left|  \sqrt{ \vfi''({|\D u | + |\D u \circ T_h| }) }|  \D (w \eta)| \right|_{L^2 ({(I_{\rho^2})_h \times B_\rho} )}.
\end{multline}
So far, we have not used the $p$-structure in our estimates. We do it now, estimating the last term of \eqref{jb_strong:t:bib:3}  via $\vfi'' (t) \le  C (1 + t^{p-2}) $. We get

(i) In the case $p=2$ 
\begin{multline}\label{jb_strong:t:bib:4:2}
h^{-\alpha} \Big|  \int_{  (I_{\rho^2})_h } \langle ( \Delta^h u \eta )_t , w  \rangle_{(W^{1,2}_{0} (B_\rho))^*\!,\, W^{1,2}_{0} (B_\rho  )} \Big| \le \\
  C (p)  
\Big| \frac{\Delta^h \V (\D u)}{h^\alpha} \Big|_{ L^2 ({(I_{\rho^2})_h \times B_\rho}) }  ( |\eta|_\infty + | \nabla \eta|_\infty ) |w|_{L^2 ((I_{\rho^2})_h; W^{1,2} (B_\rho))}.
\end{multline}
(ii) For $p>2$ we have 
\begin{multline}\label{jb_strong:t:bib:4:p}
h^{-\alpha} \Big|  \int_{  (I_{\rho^2})_h } \langle ( \Delta^h u \eta )_t , w  \rangle_{(W^{1,p}_{0} (B_\rho))^*\!,\, W^{1,p}_{0} (B_\rho  )}   \Big| \le \\
 C (p)
\Big|\frac{\Delta^h \V (\D u)}{h^\alpha} \Big|_{ L^2 ({(I_{\rho^2})_h \times B_\rho}) }  \sqrt{ \int_{{(I_{\rho^2})_h \times B_\rho}}  \left( 1+   |\D u | + |\D u  \circ T_h|  \right)^{p-2} |  \D (w \eta)|^2}.
\end{multline}
Next we split the square root in \eqref{jb_strong:t:bib:4:p}  by the H\"older inequality with powers $\frac{p}{p-2}, \; \frac{p}{2}$ to obtain
\begin{equation}\label{jb_strong:t:bib:44:p}
 \sqrt{ \int_{{(I_{\rho^2})_h \times B_\rho}}  \left( 1+   |\D u | + |\D u \circ T_h |  \right)^{p-2} |  \D (w \eta)|^2} \le C \left( 1+  \left|  \D u \right|_{ L^p (\O_I)}  \right)^{\frac{p-2}{2}}  \left|  \D (w \eta) \right|_{ L^p \left({(I_{\rho^2})_h \times B_\rho} \right)  },
\end{equation}
where we have also increased the domain of integration in the terms containing $\D u$. Use \eqref{jb_strong:t:bib:44:p} in \eqref{jb_strong:t:bib:4:p} to get for any ${w \in  L^p( (I_{\rho^2})_h; W^{1,p}_{0} (B_\rho) )}$
\begin{multline}\label{jb_strong:t:bib:5:p}
h^{-\alpha} \Big|  \int_{  (I_{\rho^2})_h} \langle ( \Delta^h u \eta )_t , w  \rangle_{(W^{1,p}_{0} (B_\rho))^*, W^{1,p}_{0} (B_\rho  )} \Big| \le\\
  C (p) ( |\eta|_\infty + | \nabla \eta|_\infty ) 
\Big| \frac{\Delta^h \V (\D u)}{h^\alpha} \Big|_{ L^2 \left({(I_{\rho^2})_h \times B_\rho} \right) }  \left( 1+  \left|  \D u \right|_{ L^p (\O_I)}  \right)^{\frac{p-2}{2}} |w|_{L^p ((I_{\rho^2})_h; W_0^{1, p} (Q_\rho))}  
 \end{multline}
Both cases (i) and (ii), \emph{\emph{i.e.}} estimates \eqref{jb_strong:t:bib:4:2}, \eqref{jb_strong:t:bib:5:p} imply that for any $h \in \left(0, \sfrac{\delta}{2} \right]$
\begin{multline*}
h^{-\alpha}  \sup_{|w|_{L^p \left( (I_{\rho^2})_h ; \, W^{1, p }_{0} (B_\rho) \right) } \le 1 } \int_{  (I_{\rho^2})_h } \left| \langle \Delta^h (u \eta)_t, w  \rangle_{(W^{1,p}_{0} (B_\rho))^*\!,\, W^{1,p}_{0} (B_\rho  )}   \right|\le \\
  C (p)( |\eta|_\infty + | \nabla \eta|_\infty ) 
\Big| \frac{\Delta^h \V (\D u)}{h^\alpha} \Big|_{ L^2 \left({(I_{\rho^2})_h \times B_\rho} \right) }  \left( 1+  \left|  \D u \right|_{ L^p (\O_I)}  \right)^{\frac{p-2}{2}},
\end{multline*}
where we have used for the case $p=2$ the convention $\infty^0 = 1$. Hence, applying $\sup_{h \in \left(0, \sfrac{\delta}{2} \right]}$ to both sides of the above estimate we have
\begin{equation}\label{eq:jb_strong:iter_1:superquad_pre}
[(u \eta)_t ]_{{1, \frac{\delta}{2},} N^{\alpha, p'} ( I; W^{-1, p'} (B_\rho) )} \le 
C (p)  ( |\eta|_\infty + | \nabla \eta|_\infty ) 
[\V (\D u)]_{{1, \frac{\delta}{2},} N^{\alpha, 2} ( I_{\rho^2};  L^2( B_\rho))}  \left( 1+  \left|  \D u \right|_{ L^p (\O_I)}  \right)^{\frac{p-2}{2}}.
\end{equation}
Observe that  
\begin{equation}\label{eq:jb_strong:iter_1:superquad_pre:b_low}
| u \eta |_{L^{ p'} ( I_{\rho^2}; W^{-1, p'} (B_\rho) )} = \sup_{|w|_{L^p ( I_{\rho^2} ; W^{1, p}_{0} (B_\rho)) } \le 1} \Big| \int_{Q_\rho} (u \eta) w\Big| \le  C |u \eta|_{L^\infty (I_{\rho^2}; L^{2} (B_\rho))},
\end{equation}
where the inequality holds by a generous use of Sobolev embedding and the H\"older inequality. 

Now, in the case $\alpha \in [0,1)$, we can use the~reduction inequality \eqref{eq:redNik} of Proposition \ref{lem:jb_strong:nikolski_red} with choices $r:=1$,  $\beta:=1$, $\delta:=\frac{\delta}{2}$ and $X:=W^{-1, p'} (B_\rho) $  for the l.h.s. of \eqref{eq:jb_strong:iter_1:superquad_pre}. It gives
\[
[u \eta]_{{2, \frac{\delta}{2},} N^{1+\alpha, p'} ( I; W^{-1, p'} (B_\rho) )} \le [(u \eta)_t ]_{{1, \frac{\delta}{2},} N^{\alpha, p'} ( I; W^{-1, p'} (B_\rho) )},
 \]
so \eqref{eq:jb_strong:iter_1:superquad_pre} implies
\begin{equation}\label{eq:jb_strong:iter_1:superquad_pre:b}
[u \eta]_{{2, \frac{\delta}{2},} N^{1+\alpha, p'} ( I_{\rho^2}; W^{-1, p'} (B_\rho) )}  \le  C (p)  ( |\eta|_\infty + | \nabla \eta|_\infty ) 
[\V (\D u)]_{{1, \frac{\delta}{2},} N^{\alpha, 2} ( I_{\rho^2};  L^2( B_\rho))}    \left( 1+  \left|  \D u \right|_{ L^p (\O_I)}  \right)^{\frac{p-2}{2}}.
\end{equation}
We want to have a full norm on the l.h.s. of \eqref{eq:jb_strong:iter_1:superquad_pre:b}, therefore we add \eqref{eq:jb_strong:iter_1:superquad_pre:b_low} to \eqref{eq:jb_strong:iter_1:superquad_pre:b} and obtain \eqref{eq:jb_strong:iter_1:superquad}.

 In the case $\alpha =1$ we drop out of the Nikolskii-Bochner class into the Sobolev-Bochner class. In order to deal with first-order (lower) seminorm here,  we add \eqref{eq:jb_strong:iter_1:superquad_pre:b}  to 
\begin{equation}\label{eq:jb_weak:theo:gentd}
| (u \eta)_t |_{ (L^p (I; W_0^{1, p} (\O)))^*} \le |\D (u \eta)|^{p-1}_{L^p}
\end{equation}
that holds for any weak solution. Consequently
\begin{multline}\label{eq:jb_strong:iter_1:superquad:sobB}
 |(u \eta)_t |_{{ \frac{\delta}{2},} W^{1, p'} ( I; W^{-1, p'} (B_\rho) )} \le  \\
 C (p)  \left( 1+  \left|  \D u \right|_{ L^p (\O_I)}  \right)^{\frac{p-2}{2}}  \left( ( |\eta|_\infty + | \nabla \eta|_\infty ) 
[\V (\D u)]_{1,{ \frac{\delta}{2},} N^{1, 2} ( I_{\rho^2};  L^2( B_\rho))} 
+  1 \right) + ( |\eta|_\infty + | \nabla \eta|_\infty )^{p-1}    \left|  u \right|^{p-1}_{ L^p (I; W^{1,p} (\Omega))} ,
\end{multline}
where we can write the Sobolev norm thanks to Corollary \ref{cor65}. It says also that the l.h.s. above dominates 
 \[
 {6}{\delta^{-1}} |(u \eta),_t \!|_{ W^{1, p'} ( I; W^{-1, p'} (B_\rho) )}.
 \]
We use this in \eqref{eq:jb_strong:iter_1:superquad:sobB}, add to both sides of the resulting estimate \eqref{eq:jb_strong:iter_1:superquad_pre:b_low} and obtain \eqref{eq:jb_strong:iter_1:superquad:sob}.
\end{proof}
\subsubsection{High space-regularity estimate in Nikolskii-Bochner spaces}
We have
\begin{lem}\label{lem:jb_strong:iter_2}
We take an untercylinder and an \"ubercylinder $\underaccent{\ddot}{Q} \Subset \ddot{Q} \Subset \O_I$ fixed by Assumption \ref{rem:uberunter}.
 Fix any  $\delta >0$ so small that $\ddot{Q}\subset \O_I^\delta$.
 Take a local weak solution $u$ to \eqref{eq:jb_structure:weak} on $\O_I$. For any $Q_{\rho_1}, Q_{\rho_2} $ such that $ \underaccent{\ddot}{Q}   \subset Q_{\rho_1} \Subset Q_{\rho_2}  \subset \ddot{Q} $ one has
  \begin{multline}\label{eq:jb_strong:2+:t:frac}
 |u|_{{1, \frac{\delta}{2},} N^{\alpha, \infty} (  I_{\rho^2_1}; L^2 (B_{\rho_1}) )}  + |\V (\D u)|_{{1, \frac{\delta}{2},} N^{\alpha, 2} (  I_{\rho^2_1}; L^2 (B_{\rho_1}) )}  \le \\
  \frac{C (p) }{({\rho_2}-{\rho_1})} \big(  1+ | \D u|_{ L^p (\O_I)}^\frac{p-2}{2} \big)  [ u]_{{1, \frac{\delta}{2},} N^{\alpha, p} (  I_{\rho^2_2}; L^p (B_{\rho_2}) )} +
   |u|_{L^\infty (I; L^{2} (\O) )} +  |u|^\frac{p}{2}_{L^p (I; W^{1,p} (\O) )}
\end{multline}
for any $\alpha \in [0,1) $ and 
\begin{multline}\label{eq:jb_strong:2+:t:frac:sob}
 |u|_{{1, \frac{\delta}{2},} W^{1, \infty} (  I_{\rho^2_1}; L^2 (B_{\rho_1}) )}  + |\V (\D u)|_{{1, \frac{\delta}{2},} W^{1, 2} (  I_{\rho^2_1}; L^2 (B_{\rho_1}) )}  \le \\
  \frac{C ( p) }{({\rho_2}-{\rho_1})} \big(   1+ | \D u|_{ L^p (\O_I)}^\frac{p-2}{2} \big)  [ u]_{{1, \frac{\delta}{2},} W^{1, p} (  I_{\rho^2_2}; L^p (B_{\rho_2}) )} +
   |u|_{L^\infty (I; L^{2} (\O) )} +  |u|^\frac{p}{2}_{L^p (I; W^{1,p} (\O) )}
\end{multline}
 provided the r.h.s.'s above are meaningful. For the case $p=2$ we apply the convention {\rm >>}$\infty^0 = 1${\rm<<}.
\end{lem}
\begin{proof} Let us take smooth cutoff functions $\sigma \in \test( I^\delta)$,  $ \eta \in C^\infty_0 (\O^\delta )$. Corollary \ref{cor:jb_strong:intro} shows that  $\Delta^{(0, h)}  u(t) \eta^2 \sigma^2 (t)$ is an admissible test function in \eqref{eq:jb_weak:theo:pointwise_t:t} for a.a. $t \in I^\delta$.
A standard computation implies  for any  $h \in (0, \sfrac{\delta}{2}]$ (for more details, see Lemma 4.7 of \cite{Burphd} or compare subsection 5.2.1 of \cite{BurKap14a}).
    \begin{multline}\label{eq:jb_strong:intro:delta_id_x:t'}
\sup_{\tau \in I^\delta} \int_{\O^\delta} \left|  \left(\Delta_{h} u \right)  \eta \sigma \right|^2 (\tau)+  \int_{I^\delta} \int_{\O^{\delta}}  \left| \left( \Delta_{h} \V ( \D u ) \right) \eta \sigma \right|^2 \le \\
  C (p) \int_{I^\delta}  \int_{\O^{\delta}} ( |\Delta_{h}u |^2   |\sigma,_\tau |_\infty + \vfi'' (|\D u | + |\D u \circ T_h |)|\Delta_h u|^2  | \nabla \eta |^2_\infty.
\end{multline}
 Now we choose in \eqref{eq:jb_strong:intro:delta_id_x:t'} $\eta \sigma$ so that it equals $1$ on $(I_{\rho_1^2})_h \times B_{\rho_1}$, vanishes outside $(I_{\rho_2^2})_h \times B_{\rho_2}$ and \[|\sigma,_\tau |\le \frac{4}{\rho_2^2 - h - (\rho_1^2-h)} \vee \frac{4}{\rho_2^2 - \rho_1^2} \le \frac{4}{(\rho_2 - \rho_1)^2}.\]. It gives
 \vspace{-3pt}
   \begin{multline} \label{eq:jb_strong:iter_2:01}
 \sup_{\tau \in (I_{\rho_1^2})_h}  \int_{ B_{\rho_1}}\!\! \left| \Delta_h u  \right|^2  (\tau)+ \! \int_{(I_{\rho_1^2})_h \times B_{\rho_1}}  \! \left| \Delta_h \V(\D u) \right|^2 \le \\
 \frac{  C (p)   }{({\rho_2}-{\rho_1})^2}  \! \int_{(I_{\rho_2^2})_h \times B_{\rho_2}} \! \left( \left| \Delta_h u \right|^2 + \vfi'' (|\D u | + |\D u \circ T_h |)|\Delta_h u|^2 \right),
\end{multline}

After multiplication by $h^{- 2 \alpha}$ and use of Assumption \ref{ass:struc}, \eqref{eq:jb_strong:iter_2:01}  becomes 
   \begin{multline}  \label{eq:jb_strong:iter_2:02}
\sup_{\tau \in (I_{\rho^2_1})_h}  \int_{B_{\rho_1}} \Big| \frac{ \Delta_h u}{h^\alpha}  \Big|^2  (\tau)+  \! \int_{(I_{\rho_1^2})_h \times B_{\rho_1}}  \!  \Big| \frac{\Delta_h \V (\D u)}{h^\alpha}  \Big|^2 \le \\
  \frac{C (p) }{({\rho_2}-{\rho_1})^2}  \! \int_{(I_{\rho_2^2})_h \times B_{\rho_2}} \!  (1+|\D u | + |\D u \circ T_h |)^{p-2}  \Big|\frac{ \Delta_h u}{h^\alpha}  \Big|^2.
\end{multline}
Now we follow the proof of Lemma \ref{lem:jb_strong:iter_1}. We split the last term with the H\"older inequality with powers $\frac{p}{p-2}, \; \frac{p}{2}$ and obtain
  \begin{equation} \label{eq:jb_strong:iter_2:1}
\!\!\sup_{\tau \in (I_{\rho^2_1})_h}   \intop_{B_{\rho_1}} \!\! \Big| \frac{ \Delta_h u}{h^\alpha}  \Big|^2 \! (\tau)\,+ \!\!\!\!\!\! \intop_{(I_{\rho_1^2})_h \times B_{\rho_1}}  \!\!\!  \Big|\frac{\Delta_h \V (\D u)}{h^\alpha} \Big|^2\!\le \!
  \frac{C (p) }{({\rho_2}-{\rho_1})^2} \left(  1\!+\! | \D u|_{ L^p (\O_I)}^{p-2} \right) \!  \Big|\frac{ \Delta_h u}{h^\alpha}  \Big|_{ L^p ((I_{\rho_2^2})_h; L^p( B_{\rho_2}) )  }^2\!\!.
\end{equation}
We apply $\sup_{h \in \left(0, \sfrac{\delta}{2} \right]}$ to the both sides of the above estimate and end up with
\begin{multline}\label{eq:jb_strong:iter_2:1:1}
 [ u]_{{1, \frac{\delta}{2},} N^{\alpha, \infty} (  I_{\rho^2_1}; L^2 (B_{\rho_1}) )}  + [\V (\D u)]_{{1, \frac{\delta}{2},} N^{\alpha, 2} (  I_{\rho^2_1}; L^2 (B_{\rho_1}) )}  \le \\
  \frac{C (p) }{({\rho_2}-{\rho_1})} \big(  1+ | \D u|_{ L^p (\O_I)}^\frac{p-2}{2} \big)  [ u]_{{1, \frac{\delta}{2},} N^{\alpha, p} (  I_{\rho^2_2}; L^p (B_{\rho_2}) )}. 
\end{multline}
In order to gain control over the full norms, we add to \eqref{eq:jb_strong:iter_2:1:1} the low-order estimate
\[|u|_{L^\infty (I_{\rho^2_1}; L^{2} (B_{\rho_1}))} +  | \V (\D u)|_{L^2 (I_{\rho^2_1}; L^{2} (B_{\rho_1}))} \le    |u|_{L^\infty (I; L^{2} (\O) )} +  |u|^\frac{p}{2}_{L^p (I; W^{1,p} (\O) )}.
\]
Hence we arrive at \eqref{eq:jb_strong:2+:t:frac} (the case $\alpha \in [0,1)$). Finally we also get \eqref{eq:jb_strong:2+:t:frac:sob}  via Corollary \ref{cor65} (the case $\alpha = 1$).
\end{proof}

\subsubsection{A keystone inequality} Recall that our goal in this section is to develop estimates that shall be used in Lemma \ref{lem:jb_strong:nikolski_inter}. For the time being,  Lemma \ref{lem:jb_strong:iter_1} provides its  high-time and low-space regularity leg. To have its  low-time and high-space regularity ingredient we would like to control
\[|u|_{1, \frac{\delta}{2}, N^{\alpha_1, p_1} ( I_{\rho^2_1}; W_0^{1, q_1} (B_{\rho_1})) }\]
for certain $\alpha_1, p_1, q_1$. 
However,  Lemma \ref{lem:jb_strong:iter_2} only controls
\[
 |u|_{{1, \frac{\delta}{2},} N^{\alpha, \infty} (  I_{\rho^2_1}; L^2 (B_{\rho_1}) )}  + |\V (\D u)|_{{1, \frac{\delta}{2},} N^{\alpha, 2} (  I_{\rho^2_1}; L^2 (B_{\rho_1}) )}.
\]
We bridge this gap in the following lemma.\begin{lem}\label{cor:jb_strong:iter:p} 
Let us take an untercylinder and an \"ubercylinder $\underaccent{\ddot}{Q} \Subset \ddot{Q} \Subset \O_I$ fixed by Assumption \ref{rem:uberunter}.
 Fix any  $\delta >0$ so small that $\ddot{Q}\subset \O_I^\delta$ and a smooth time-independent $\eta \in C^\infty (\O)$ such that $|\eta| \le 1$ and $|\nabla \eta|_\infty \ge 1$. Take a local weak solution $u$ to \eqref{eq:jb_structure:weak} on $\O_I$. For any $Q_\rho$ such that $ \underaccent{\ddot}{Q} \subset Q_\rho \subset \ddot{Q} $ and any $\alpha \in [0,1) $, provided the~r.h.s.'s are meaningful, one has
  \begin{multline}\label{eq:jb_strong:iter:cor'a}
|u \eta |^2_{1, \delta, N^{\alpha, 2} (  I_{\rho^2}; W^{1,2} (B_{\rho}) )} + |u \eta |^p_{1, \delta, N^{\frac{2\alpha}{p}, p} (  I_{\rho^2}; W^{1,p} (B_{\rho}) )}   \le 
 \frac{C (p) }{\left(\vfi'' (0)\right)^2} \left(1 + | \nabla \eta|^{p}_\infty \right) \times \\
 \Big( [\V ( \D u)]^2_{1, \delta, N^{\alpha, 2} (  I_{\rho^2}; L^{2} (B_{\rho}) )} +  [u ]^2_{1, \delta, N^{\alpha, 2} (  I_{\rho^2}; L^{2} (B_{\rho}) )} + \\
  [u ]^p_{1, \delta, N^{\frac{2\alpha}{p}, p} (  I_{\rho^2}; L^{p} (B_{\rho}) )} \!+  |u |^p_{L^{ p} ( I_{\rho^2}; W^{1, p} (B_\rho)) }  \!+  |u |^2_{L^{ 2} ( I_{\rho^2}; W^{1, 2} (B_\rho)) }  \Big)
\end{multline}
\vspace{-4pt}
for  $\alpha \in [0, 1)$ and either 
\begin{multline}\label{eq:jb_strong:iter:cor'a:sob}
|u \eta |^2_{1, \delta, W^{1, 2} (  I_{\rho^2}; W^{1,2} (B_{\rho}) )} + |u \eta |^p_{1, \delta, N^{\frac{2}{p}, p} (  I_{\rho^2}; W^{1,p} (B_{\rho}) )}   \le 
 \frac{C (p) }{\left(\vfi'' (0)\right)^2} \left(1 + | \nabla \eta|^{p}_\infty \right) \times \\
 \Big( [\V ( \D u)]^2_{1, \delta, W^{1, 2} (  I_{\rho^2}; L^{2} (B_{\rho}) )} +  [u ]^2_{1, \delta, W^{1, 2} (  I_{\rho^2}; L^{2} (B_{\rho}) )} + \\
  [u ]^p_{1, \delta, N^{\frac{2}{p}, p} (  I_{\rho^2}; L^{p} (B_{\rho}) )} \!+  |u |^p_{L^{ p} ( I_{\rho^2}; W^{1, p} (B_\rho)) }  \!+  |u |^2_{L^{ 2} ( I_{\rho^2}; W^{1, 2} (B_\rho)) }  \Big)
\end{multline}
\vspace{-4pt}
in the case $p>2$ or 
\begin{multline}\label{eq:jb_strong:iter:cor'a:sob:p2}
|u \eta |^2_{1, \delta, W^{1, 2} (  I_{\rho^2}; W^{1,2} (B_{\rho}) )}    \le \\
 \frac{C (p)}{\left(\vfi'' (0)\right)^2} \left(1 + | \nabla \eta|^{2}_\infty \right) \left( [\V ( \D u)]^2_{1, \delta, W^{1, 2} (  I_{\rho^2}; L^{2} (B_{\rho}) )} +  [u ]^2_{1, \delta, W^{1, 2} (  I_{\rho^2}; L^{2} (B_{\rho}) )} +  |u |^2_{L^{ 2} ( I_{\rho^2}; W^{1, 2} (B_\rho)) }  \!\right)
\end{multline}
in the case $p=2$.
\end{lem}
\begin{proof}
Let us  consider $ \vfi (|\Delta^h \D (u \eta)|) $. It holds
\begin{equation}\label{eq:jb_strong:it_qstr2:p}
\vfi (|\Delta^h \D (u \eta)|) \le  \vfi (|\eta|_\infty |\Delta^h \D u | + | \nabla \eta|_\infty |\Delta^h u |)  \le \\
 C (\Delta_2 (\vfi)) \left(  \vfi (|\eta|_\infty |\Delta^h \D u | )+    \vfi ( | \nabla \eta|_\infty |\Delta^h u |) \right).
\end{equation}
Assumptions \ref{ass:struc} yields
\begin{equation}\label{eq:jb_strong:it_qstr2:p:0}
  \frac{ \vfi'' (0) }{p(p-1)}  (t^2+ t^{p}) \le  \vfi (t) \le  \frac{1}{2} C (t^2+ t^{p}).
\end{equation}
Hence, using $|\eta| \le 1 $  we get
\[
\vfi (|\eta|_\infty |\Delta^h \D u | ) \le C |\eta|^2_\infty \left( |\Delta^h \D u |^2 + |\Delta^h \D u |^p  \right) \le \frac{C (p)}{\vfi'' (0)} \vfi ( |\Delta^h \D u | )
\]
and similarly, via $|\nabla \eta|_\infty \ge 1$,
\vspace{-8pt}
\[
\vfi (|\nabla \eta|_\infty |\Delta^h u | ) \le \frac{C (p)}{\vfi'' (0)} |\nabla \eta|^p_\infty \vfi ( |\Delta^h u | ).
\]
Using the above two inequalities in \eqref{eq:jb_strong:it_qstr2:p}, we obtain
\begin{equation}\label{eq:jb_strong:it_qstr2:p:1}
\vfi (|\Delta^h \D (u \eta)|) \le  \frac{C (p) }{\vfi'' (0)} \left( \vfi ( |\Delta^h \D u | ) + |\nabla \eta|^p_\infty \vfi ( |\Delta^h u | ) \right).
\end{equation}
Since Assumption \ref{ass:struc} and \eqref{eq:lem:jb_weak:mon_equiv2} hold, from \eqref{eq:jb_strong:it_qstr2:p:1} we arrive at
\begin{equation}\label{eq:jb_strong:it_qstr2:p:2}
\vfi (|\Delta^h \D (u \eta)|) \le  \frac{C (p) }{\vfi'' (0)} \left( |\Delta^h \V ( \D u) |^2 + |\nabla \eta|^p_\infty \vfi ( |\Delta^h u | ) \right).
\end{equation}
We use \eqref{eq:jb_strong:it_qstr2:p:2}  to write
\begin{equation}\label{eq:jb_strong:it_qstr2:p:3}
\sup_{h \in (0, \delta)} h^{-2 \alpha} \int_{(I_{\rho^2})_h \times B_\rho} \!  \vfi (|\Delta^h \D (u\eta)|) \le  \\
\frac{C (p)  }{\vfi'' (0)}\sup_{h \in (0, \delta)} h^{- 2 \alpha} \int_{(I_{\rho^2})_h \times B_\rho} \! |\Delta^h \V ( \D u) |^2 + |\nabla \eta|^p_\infty \, \vfi ( |\Delta^h u | ).
\end{equation}
Next, we use \eqref{eq:jb_strong:it_qstr2:p:0} in the l.h.s. of \eqref{eq:jb_strong:it_qstr2:p:3} to get
\begin{multline}\label{eq:jb_strong:it_qstr2'a:p:4}
\sup_{h \in (0, \delta)} \int_{(I_{\rho^2})_h \times B_\rho} \! \Big| \frac{\Delta^h \D (u\eta)}{h^\alpha} \Big|^2 +  \Big| \frac{\Delta^h \D (u\eta)}{h^\frac{2\alpha}{p}} \Big| ^p  \le \\
 \frac{C (p) }{\left(\vfi'' (0)\right)^2} \left(1 + | \nabla \eta|^{p}_\infty \right) \sup_{h \in (0, \delta)} \int_{(I_{\rho^2})_h \times B_\rho} \! \Big|\frac{\Delta^h \V ( \D u)}{h^\alpha} \Big|^2 +\Big| \frac{\Delta^h u}{h^\alpha}\Big|^2 + \Big|\frac{\Delta^h u}{h^\frac{2\alpha}{p}} \Big|^p.
\end{multline}
We can add to both sides harmlessly 
\[\sup_{h \in (0, \delta)} \int_{(I_{\rho^2})_h \times B_\rho} \!   \Big| \frac{\Delta^h (u \eta)}{h^\alpha}\Big|^2  +\Big|\frac{\Delta^h (u \eta)}{h^\frac{2\alpha}{p}}\Big|^p, 
\]
 because $\eta$ is time-independent. Hence 
 \begin{multline}\label{eq:jb_strong:it_qstr2'a:p:5}
\sup_{h \in (0, \delta)} \int_{(I_{\rho^2})_h } \int_{ B_\rho} \! \Big| \frac{\Delta^h \D (u\eta)}{h^\alpha}\Big|^2 +  \Big|\frac{\Delta^h (u \eta)}{h^\alpha}\Big|^2 +  \Big|\frac{\Delta^h \D (u\eta)}{h^\frac{2\alpha}{p}} \Big|^p  +\Big|\frac{\Delta^h (u \eta)}{h^\frac{2\alpha}{p}} \Big|^p  \le \\
 \frac{C (p) }{\left(\vfi'' (0)\right)^2} \left(1 + | \nabla \eta|^{p}_\infty \right) \sup_{h \in (0, \delta)} \int_{(I_{\rho^2})_h} \int_{ B_\rho} \! \Big|\frac{\Delta^h \V ( \D u)}{h^\alpha} \Big|^2 +\Big|\frac{\Delta^h u}{h^\alpha}\Big|^2 + \Big|\frac{\Delta^h u}{h^\frac{2\alpha}{p}}\Big|^p.
\end{multline}
The commutativity of $\Delta^h $ and $\D$ and the Korn inequality allow us to obtain from \eqref{eq:jb_strong:it_qstr2'a:p:5}
  \begin{multline}\label{eq:jb_strong:it_qstr2'a:p:6}
[u \eta ]^2_{1, \delta, N^{\alpha, 2} (  I_{\rho^2}; W^{1,2} (B_{\rho}) )} + [u \eta ]^p_{1, \delta, N^{\frac{2\alpha}{p}, p} (  I_{\rho^2}; W^{1,p} (B_{\rho}) )}   \le \\
 \frac{C (p)  }{\left(\vfi'' (0)\right)^2} \left(1 + | \nabla \eta|^{p}_\infty \right) \Big( [\V ( \D u)]^2_{1, \delta, N^{\alpha, 2} (  I_{\rho^2}; L^{2} (B_{\rho}) )} +  [u ]^2_{1, \delta, N^{\alpha, 2} (  I_{\rho^2}; L^{2} (B_{\rho}) )}+  [u ]^p_{1, \delta, N^{\frac{2\alpha}{p}, p} (  I_{\rho^2}; L^{p} (B_{\rho}) )} \Big)
\end{multline}
 Adding to both sides of  \eqref{eq:jb_strong:it_qstr2'a:p:6}
\[
|u \eta |^2_{L^{ 2} (  I_{\rho^2}; W^{1, 2} (B_\rho)) } + |u \eta |^p_{L^{ p} (  I_{\rho^2}; W^{1, p} (B_\rho)) }   \le C \left( 1+ |\nabla \eta|^p_\infty \right) \left(|u |^2_{L^{ 2} (  I_{\rho^2}; W^{1, 2} (B_\rho)) }  + |u |^p_{L^{ p} (  I_{\rho^2}; W^{1, p} (B_\rho)) } \right)
 \]
 we form the Nikolskii-Bochner norms in the l.h.s. of  \eqref{eq:jb_strong:it_qstr2'a:p:6}. It gives in the case $\alpha \in [0,1)$ the wanted estimate \eqref{eq:jb_strong:iter:cor'a} and in the case $\alpha =1$, via Corollary \ref{cor65},  estimates \eqref{eq:jb_strong:iter:cor'a:sob}, \eqref{eq:jb_strong:iter:cor'a:sob:p2}.
 \end{proof}
 \begin{rem}\label{rem:notp}
One can  restate Lemma \ref{cor:jb_strong:iter:p} without Assumption \ref{ass:struc} using instead of the $p$-growth -- the Boyd indices $q_1, q_2$ in an Orlicz-growth setting.
\end{rem}

\subsection{Iteration in  Nikolskii-Bochner spaces }\label{sec:jb_strong:iteration} 
Endowed with the estimates of Section \ref{sec:energyinNik}, we are now ready to use the interpolation Lemma  \ref{lem:jb_strong:nikolski_inter} recursively,  raising the time-regularity of a local weak solution to \eqref{eq:jb_structure:weak}.
For clarity of exposition, let us begin the presentation of our iteration procedure without the quantitative details. 
\subsubsection{A qualitative iteration}\label{ssec:idea:ind}  First we present how one raises the regularity in a single step of our iteration.
\paragraph{An iteration step} We take the inductive assumption that 
\begin{equation}\label{idea:indAss}
u \in N^{\alpha_i, p} (  I_{\rho^2_i}; L^p (B_{\rho_i}) ) \tag{$\alpha_i$}
\end{equation}
for a certain $\alpha \in [0,1)$ and we want to raise the fractional time differentiability of $u$ up to ($\alpha_{i+1}$) with the aid of energy estimates of Section \ref{sec:energyinNik}. We do the following
 \begin{itemize}
\item[(i)] Assumption \eqref{idea:indAss} implies via Lemma \ref{lem:jb_strong:iter_2}  that on a smaller cylinder $I_{\tilde \rho^2_i} \times B_{\tilde \rho_i}$ we have
\begin{equation}\label{idea:indAss:i} 
\V (\D u) \in N^{\alpha_i, 2} (I_{\tilde \rho^2_i} ; L^2 (B_{\tilde \rho_i})).
\end{equation}
\item[(ii)] The estimate \eqref{idea:indAss:i} via Lemma \ref{lem:jb_strong:iter_1}  gives the high-time and low-space information 
\begin{equation}\label{eq:jb_strong:htls}
u \eta \in {N^{1+ \alpha_i, p'} \left(I_{\tilde \rho^2_i}; W^{-1, p'} (B_{\tilde \rho_i})\right) }.
\end{equation}
Moreover, \eqref{idea:indAss:i} via Lemma \ref{cor:jb_strong:iter:p} and the assumption  \eqref{idea:indAss} gives also the low-time and high-space information  
 \begin{equation}\label{eq:jb_strong:hslt}
u \eta \in {N^{\frac{2\alpha_i}{p}, p}(I_{\tilde \rho_i^2}; W^{1, p} (B_{\tilde \rho_i})) }.
 \end{equation}
\item[(iii)] We intend to interpolate the high-time and low-space \eqref{eq:jb_strong:htls} with the low-time and high-space \eqref {eq:jb_strong:hslt}  pieces of information by means of Lemma \ref{lem:jb_strong:nikolski_inter}. In order to use it we need zero space-trace functions. Therefore we choose an appropriate spatial cutoff function $\eta$, which equals $1$ on $B_{ \rho_{i+1}}$ slightly smaller than $B_{\tilde \rho_i}$. Now we take in Lemma  \ref{lem:jb_strong:nikolski_inter} the following parameters
\begin{equation}\label{eq:jb_strong:inter:choice}
\alpha_1 :=  \frac{2 \alpha_i}{p},  \quad \alpha_2 := 1+ \alpha_i, \quad p_1 = q_1  := p, \; \quad p_2 = q_2 :=  p'.
\end{equation}
The choice \eqref{eq:jb_strong:inter:choice} is admissible in  Lemma \ref{lem:jb_strong:nikolski_inter}.
As a result we obtain  
\begin{equation}\label{eq:jb_strong:inter:result:a}
u \eta \in {N^{\alpha', p_0} ( I_{\tilde \rho_i^2}; L^{q_0} (B_{\tilde \rho_i}))},
  \end{equation}
 where
 \begin{itemize}
  \item[(iii.i)] For the case $p=2$ we simply choose $\theta=1$ in  Lemma \ref{lem:jb_strong:nikolski_inter} and get
\begin{equation}\label{eq:jb_strong:inter:result:p2}
\alpha' = \alpha_i+ \frac{1}{2}, \qquad p_0 =  2, \qquad q_0 = 2.
 \end{equation}
Hence we can proceed immediately to the next step ($\alpha_{i+1}$) with $\alpha_{i+1} = \alpha_i + \sfrac{1}{2}$.
 \item[(iii.ii)] In the case $p>2$ we have the following parameters in Lemma \ref{lem:jb_strong:nikolski_inter}
\begin{equation}\label{eq:jb_strong:inter:result}
\begin{aligned}
&\alpha' = \frac{\alpha_i}{p} \left(2 + \frac{p-2}{2}  \theta   \right) + \frac{\theta}{2}\\
&p_0 =  \frac{2p }{\theta (p-2) + 2}, \qquad 
q_0 = \begin{cases} \tilde q_0:=  \frac{2p }{\theta \left( p-2  + \frac{2p}{d} \right)+ 2 (1-  \frac{p}{d}) } &\text{ when } \tilde q_0 \in (1, \infty ), \\   \text{ any finite number } &\text{ otherwise. } \end{cases}
\end{aligned}
  \end{equation}
Observe that for $\theta \le 1$ one has $p_0 \le q_0$.
\end{itemize}
  \item[(iv)] We look for such $\theta$ in  \eqref{eq:jb_strong:inter:result}  that  the space ${N^{\alpha', p_0}  (L^{q_0})}$, resulting from the interpolation \eqref{eq:jb_strong:inter:result:a}, imbeds continuously into ${N^{\alpha_{i+1}, p}  ( L^{p} )}$ for a certain $\alpha_{i+1}$. We can take now two approaches.
\begin{itemize}
\item[(iv.i)] The first one is to choose $\theta$ that gives $p_0 =  p$. Knowing $ q_0 \ge p_0$, we then get \[{N^{\alpha', p_0} ( L^{q_0 } )} \csubset{N^{\alpha',  p } ( L^p )}.\] Unfortunately, in order to have $p_0 =  p$ we need $\theta =0$, see \eqref{eq:jb_strong:inter:result}. Consequently $\alpha' < \alpha_i$, so there is no increase of fractional time differentiability. 
\item[(iv.ii)] The second one is choose $\theta$ that gives $q_0 =  p$. This condition used in   \eqref{eq:jb_strong:inter:result} implies
\begin{equation}\label{eq:jb_strong:ivi:pq1}
\theta = \frac{2p}{d(p-2) + 2p}
\end{equation}
and consequently
\begin{equation}\label{eq:jb_strong:ivi:pq2}
\alpha' = \alpha_i \left( \frac{2}{p} + \frac{p-2}{d (p-2) + 2p}  \right) + \frac{p}{d (p-2) + 2p}, \qquad p_0 =  \frac{2p^2 + pd (p-2)}{p^2 + d(p-2)}.
\end{equation}
We have thus obtained that 
${N^{\alpha', p_0}  (L^{q_0})}$ of \eqref{eq:jb_strong:inter:result:a} is ${N^{\alpha', p_0}  (L^{p})}$. Recall that $p_0 \le q_0 \;(=p)$ now. Therefore to end up with $(\alpha_{i+1})$ we need to decrease $\alpha' $ to an $\alpha_{i+1}$ such that
\begin{equation}\label{eq:jb_strong:ivi}
{N^{\alpha', p_0} ( L^p)} \csubset N^{ \alpha_{i+1}, p} \left( L^p \right) 
\end{equation}
Lemma \ref{lem:nik:emb} provides for \eqref{eq:jb_strong:ivi} the following condition
\begin{equation}\label{eq:jb_strong:ivi2}
 \alpha_{i+1} =\alpha'  - \frac{1}{p_0} +  \frac{1}{p} =  \alpha_i \left( \frac{2}{p} + \frac{p-2}{d (p-2) + 2p}  \right) + \frac{2}{d (p-2) + 2p}.
 \end{equation}
\end{itemize}
Hence we have $u \eta \in {N^{\alpha_{i+1}, p} ( I_{\tilde \rho_i^2}; L^{p} (B_{\tilde \rho_i}))}$. 
 \end{itemize}
Particularly, as $\eta \equiv 1$ on $B_{ \rho_{i+1}}$ we have reached the next step
\begin{equation}\label{idea:indAss:next}
u \in N^{\alpha_{i+1}, p} (  I_{\rho^2_{i+1}}; L^p (B_{\rho_{i+1}}) ) \tag{$\alpha_{i+1}$}.
\end{equation}

Briefly, the iteration step $(\alpha_{i}) \implies (\alpha_{i+1})$ reads
 \begin{equation}\label{eq:jb_strong:iter:cor_scheme:bf}
u \!\in\! N^{ \alpha_i, p} \left( L^p \right)  \!\!  \stackrel{(i)}{\implies} \!\! \V (\D u) \!\in\! N^{\alpha_i, 2} (L^2) \!\! \stackrel{(ii)}{\implies}  \!\!\!
\begin{cases}
u \eta \!\in\! {N^{1+ \alpha_i, p'} \left( W^{-1, p'} \right) }\\
u \eta \!\in\! {N^{\frac{2\alpha_i}{p}, p}(W^{1, p} ) } 
\end{cases}
\!\!\!\!\!\! \stackrel{(iii) + (iv)}{\implies} \! u \!\in\! N^{ \alpha_{i+1}, p} \left( L^p \right),
\end{equation}
where the obtained information holds on a slightly smaller cylinder.
\paragraph{The upper bound for the possible regularity gain}\label{ssec:crone} We can carry on with our iterative scheme unless $\alpha_i  \ge 1$, because the energy estimates for the Nikolskii-Bochner spaces, derived in Section \ref{sec:energyinNik}, hold within the range $\alpha \in [0,1)$. Nevertheless, if  the iteration provides us with strictly increasing fractional differentiabilities $\alpha_i <1$, $i=0,1,\dots, i_0 -1$,  $\alpha_{i_0} \le 1$, $\alpha_{i_0+1} > 1$, we can cross the~full differentiability and reach any $\alpha < \alpha_{i_0+1}$ by interpolation. More precisely, the energy estimates hold only for $\alpha <1$, but there is no such bound in the interpolation Lemma \ref{lem:jb_strong:nikolski_inter}. Hence an allowed iteration step is $a_{i_0} - \delta (\epsilon) \to a_{i_0+1} - \epsilon$, where for any small  $\epsilon > 0$ we have $ \delta (\epsilon)>0$, due to the strict monotonicity and continuity of the iteration step with respect to the parameter of the fractional differentiability.

In other words, if the iteration process formally exceeds $1$ and the first obtained there fractional differentiability is $1 + \eps$, we can obtain for any $\eps^- < \eps$ that 
\begin{equation}\label{cross1:3}
u \in N^{1 + \eps^- , p} \left( L^{ p} \right)
\end{equation}

\paragraph{The iterative regularity gain}\label{ssec:grain} 
The $L^p (W^{1,p})$ regularity of a local weak solution to \eqref{eq:jb_structure:weak} allows us to start the iteration process \eqref{eq:jb_strong:iter:cor_scheme:bf} with $\alpha_0 = 0$.

In the case $p=2$ the fractional differentiability grows along the iteration \eqref{eq:jb_strong:iter:cor_scheme:bf} arithmetically 
\begin{equation}\label{idea:regGain:p2}
\alpha_{i+1} = \alpha_i + \frac{1}{2},
\end{equation}
compare \eqref{eq:jb_strong:inter:result:p2} of the substep (iii.i) in Subsection \ref{ssec:idea:ind}. Here we perform two steps of iteration that begins with $\alpha_0 = 0$ to get ${N^{1 , p} \left( L^{ p} \right)} $ and decrease  a little this regularity, as explained before, to cross the full differentiability. Hence we get $u \in {N^{\frac{3}{2}-\delta,\, p} \left( L^{ p} \right)} $ for an arbitrarily small $\delta$.

In the case $p>2$ the fractional differentiability grows along the iteration \eqref{eq:jb_strong:iter:cor_scheme:bf} slowly geometrically
\begin{equation}\label{idea:regGain}
 \alpha_{i+1} = \alpha_i A + B, \qquad A :=\frac{2}{p} + \frac{p-2}{d (p-2) + 2p}, \quad  B:= \frac{2}{d (p-2) + 2p},
\end{equation}
compare substep \eqref{eq:jb_strong:ivi2}  in Subsection \ref{ssec:idea:ind}. We have $A \in (0,1)$, $B> 0$ hence 
 \[
 \alpha_n = \alpha_0 A^n + B \frac{1- A^n}{1-A}.
 \]
For $\alpha_0 = 0$ we get
\begin{equation}\label{idea:regGain:p}
 \alpha_n =  B \frac{1- A^n}{1-A} \; \stackrel{n \to \infty}{\nearrow} \; \frac{B}{1-A} = \frac{2p}{(p-2)(d (p-2) + p)}.
 \end{equation}
Therefore, in the case $ \frac{B}{1-A} \le 1$, for any $\alpha < \frac{B}{1-A}$, after a finite number of iterations we get
\begin{equation}\label{idea:regGain:plow}
u \in {N^{\alpha , p} \left( L^{ p} \right)}.
 \end{equation}
Otherwise, \emph{\emph{i.e.}}  when $ \frac{B}{1-A} > 1$, we can cross the full differentiability in a finite number of iterations, as outlined at the end of the previous section, to get 
\[u \in {N^{1 + \eps^- , p} \left( L^{ p} \right)}. \]
As we have the explicit iteration formula \eqref{idea:regGain}, we realize that 
\begin{equation}\label{A+B}
1 + \eps = A+B = \frac{2}{p} + \frac{p}{d (p-2) + 2p}.
\end{equation}
Observe that it also agrees with the result of the case $p=2$, as there $ A+B = \frac{3}{2}$.
\begin{rem}\label{rem:iter:finite}
In both cases we terminate the iterations after a finite number of steps. This will be important in the following subsection on the qualitative iterations (\emph{i.e.} iterations, where we also keep track of inequalities), because it will allow to deal easily with the constants of estimates.
\end{rem}
\paragraph{A remark on a trouble}\label{ssec:trouble} 
Observe that for $p > 2$ we are in a worse situation than in \cite{BEKP11} in two aspects. 
\begin{itemize}
\item[(i)]  Firstly, because the localization produces in estimates a lower order term in $N^{ \alpha, p} \left(L^p \right)$ that we need to control, contrary to  \cite{BEKP11}, where the bad term belongs to $N^{ \alpha, 2} \left(L^2 \right)$. 
\item[(ii)] Secondly, because  our low-time and high-space information has $\frac{2\alpha}{p} < \alpha$ order of fractional time-differentiability. 
\end{itemize}
Consequently
\begin{itemize}
\item[(i)] In order to control space regularity of the  lower order term in $N^{ \alpha, p} \left(L^p \right)$ we need to interpolate known $ {N^{\frac{2\alpha}{p}, p}(W^{1, p}) } $ and $N^{1+ \alpha, p'} ( W^{-1, p'})$ not with $\theta = \frac{1}{2}$, but closer to the high-space regularity leg. This diminishes the gained fractional-time differentiability\footnote{One can choose $ {N^{\alpha, 2}(W^{1, 2}) } $ instead of $ {N^{\frac{2\alpha}{p}, p}(W^{1, p}) } $ as the high-space-regularity endpoint for interpolation. Then a formal computation indicates that a possible gain of fractional time differentiability, \emph{\emph{i.e.}} $\alpha_{i+1} \ge \alpha_{i}$ occurs only for $p \le 2 + \frac{4}{d+1}$. This is always majorized by $2 + \frac{4}{d}$ of parabolic embedding and worse than our high-regularity bound $2 + \frac{2}{\sqrt{d+1}}$ for $d \ge 3$, so of no interest to us.}.
\item[(ii)] Moreover, the gained fractional-time differentiability is smaller because the high-space regularity leg enjoys lower fractional-time differentiability $\frac{2\alpha}{p} < \alpha$.
\end{itemize}

\subsubsection{Quantitative iteration}
 For brevity, we write
\[
Y (X) (Q_{\rho}) \quad \text{for} \quad Y (  I_{\rho^2}; X (B_{\rho}) ).
\] 
Let us also recall that
\begin{equation}\label{eq:nik:g0}
\gamma_0 := \frac{2p}{(p-2)(d (p-2) + p)}
\end{equation}
which is $\frac{B}{1-A}$ from \eqref{idea:regGain:p}. Observe that
\begin{equation}\label{eq:nik:gammaP}
\gamma_0 \le 1 \quad \iff \quad p \ge 2 + \frac{2}{\sqrt{d+1}}
\end{equation}
Having presented the qualitative iteration, we are ready to provide a quantitative regularity result based on our iterative scheme. It is divided into two lemmas. The first one considers the situation when we cannot reach  the full differentiability by our iteration (case $\gamma_0 \le 1$). The second one covers the full differentiability case (case $\gamma_0 > 1$).

\begin{lem}[Fractional differentiability in time]\label{lem:jb_strong:iter_end:low}
Let us take an \"ubercylinder and an untercylinder $ \underaccent{\ddot}{Q} \Subset \ddot{Q} \Subset \O_I$ fixed by the Assumption \ref{rem:uberunter} and a local weak solution $u$ to \eqref{eq:jb_structure:weak} on $\O_I$. Moreover, assume that  $p \ge 2 + \frac{2}{\sqrt{d+1}}$.  Then for any $\alpha < \gamma_0$ of \eqref{eq:nik:g0} there exist $\kappa (p, \alpha), \; C (p, \alpha,  \vfi'' (0))$ such that for all the concentric cylinders $Q_r, Q_R$ such that $\underaccent{\ddot}{Q} \subset Q_r  \Subset Q_R \subset \ddot{Q}$ one has
  \begin{multline}\label{eq:jb_strong:iter:lem_gen_low}
 |u|_{ N^{\alpha, \infty} (L^2)(Q_{r}) }  + |\V (\D u)|_{N^{\alpha, 2} ( L^2 )(Q_{r})}  +
 |u |_{N^{1+ \alpha, p'} \left(W^{-1, p'} \right) (Q_{r})} + |u |_{N^{\alpha, 2} ( W^{1,2} )(Q_{r})} + \\
 |u |_{ N^{\alpha, p} (L^{p} )(Q_{r})} +  |u |_{ N^{\frac{2\alpha}{p}, p} (W^{1,p} )(Q_{r})}
 \le C (p, \alpha, \vfi'' (0))\left( \frac{ 1+  |u|_{L^\infty ( L^{2} )(Q_R)} +  |u |_{L^{ p} (W^{1, p} )(Q_R) } }{R-r}  \right)^{\kappa (p, \alpha)}.
 \end{multline}
\end{lem}
\begin{proof}
We choose $\delta'$ so small that $ \ddot{Q} \Subset \O^{\delta'}_I$. We take $Q_r \Subset Q_R$  and  $\alpha < \gamma_0$ as in the assumptions. From Subsection \ref{ssec:idea:ind} on qualitative iterations we already know that $\alpha$ is accessible in a finite number of steps $N (\alpha)$, compare \eqref{idea:regGain:plow}, with specified intermediate $\alpha_i, \; i= 1, \dots, N (\alpha)$. This allows us to deal with the constants that follow from the embedding and interpolation used in the iterations, \emph{i.e.} from Lemmas \ref{lem:jb_strong:nikolski_inter}, \ref{lem:nik:emb} (compare Remark \ref{rem:iter:finite}). We simply take the most restrictive constant from any of the iteration steps of \eqref{eq:jb_strong:iter:cor_scheme:bf} and denote it by $C_1$. Similarly, we fix the smallest $\delta_i \le \delta'$ needed there and denote it by $\delta$.

Recall that at every iteration step \eqref{eq:jb_strong:iter:cor_scheme:bf} we decrease a little the underlying cylinder. It needs to be taken care of, because we want to obtain our regularity result for $Q_r$ starting from $Q_R$, where $r, R$ may be arbitrarily close to each other. 
Therefore we define the intermediate cylinders $Q_{\rho^i}$ in-between $Q_r, Q_R$ as follows. The cylinder $Q_{\rho^i}$ is concentric with $Q_r$ and $\rho^i := r + (R-r) 2^{-i}$. These intermediate cylinders decrease from $Q_{\rho^0} = Q_R$  to  $Q_{\rho^\infty} = Q_r$. We will use also a cylinder $Q_{\rho^{(i, i+1)}}$ which is in the~middle between $Q_{\rho^i}$ and $Q_{\rho^{i+1}}$, \emph{i.e.} $\rho^{(i, i+1)} = \frac{\rho^i + \rho^{i+1}}{2}$. Now let us rewrite the  iteration step from Subsection \ref{ssec:idea:ind} qualitatively. For brevity, we write 
\[
|u|_E : = |u|_{L^\infty (I; L^{2} (\O) )} , \qquad |u|_M  := |\D u|_{ L^p (\O_I)},  \qquad |u|_V  := |u|_E + |u|_M.
\]
We take the inductive assumption
\begin{equation}
u \in N^{\alpha_i, p} (L^p )  (Q_{\rho_i}) \tag{$\alpha_i$}
\end{equation}
\begin{itemize}
\item[(i)] The assumption \eqref{idea:indAss} with \eqref{eq:jb_strong:2+:t:frac} of Lemma \ref{lem:jb_strong:iter_2}  gives on $Q_{\rho^{(i, i+1)}}$
\begin{equation}\label{idea:ind:det:0}
 |\V (\D u)|_{{1, \frac{\delta}{2},} N^{\alpha_i, 2} (L^2 ) (Q_{\rho^{(i, i+1)}})} \le 
2^i  \frac{C (p) }{R-r} \left(  1+  |u|_M\right)^\frac{p-2}{2}   [ u]_{{1, \frac{\delta}{2},} N^{\alpha_i, p} (L^p)  (Q_{\rho_i}) } +  C |u|_V.
\end{equation}
\item[(ii)] Let us choose smooth $\eta \in [0,1]$ such  that $\supp(\eta) \Subset B_{\rho^{(i, i+1)}}$ and $\eta \equiv 1$ on $B_{\rho_{i+1}}$, where $ |\nabla \eta| \le C 2^i  \frac{1 }{R-r}$. 

Firstly,  \eqref{eq:jb_strong:iter_1:superquad} of Lemma \ref{lem:jb_strong:iter_1} gives via (i)
\begin{multline}\label{idea:ind:det:1}
|u \eta|_{\frac{\delta}{2}, N^{1+ \alpha_i, p'} \left( W^{-1, p'} \right)  (Q_{\rho^{(i, i+1)}}) } \le 2^i  \frac{C (p) }{R-r}  [\V (\D u)]_{{1, \frac{\delta}{2},} N^{\alpha_i, 2} (L^2 ) (Q_{\rho^{(i, i+1)}})}    \left(  1+  |u|_M\right)^\frac{p-2}{2} +  C |u|_V \\
\le 4^i  \frac{C (p) }{(R-r)^2} \left( \left(  1+  |u|_M \right)^{p-2} [ u]_{{1, \frac{\delta}{2},} N^{\alpha_i, p} (L^p)  (Q_{\rho_i}) } +   |u|_V  \left(  1+  |u|_M \right)^\frac{p-2}{2}   \right),
\end{multline}
where we use also the fact that that we work with $R \le 1$, as Assumption \ref{rem:uberunter} includes $\ddot{Q}  \subset Q_1$.

Secondly, \eqref{eq:jb_strong:iter:cor'a} of Lemma \ref{cor:jb_strong:iter:p} gives 
  \begin{multline}\label{idea:ind:det:2}
 |u \eta |_{\frac{\delta}{2}, N^{\frac{2\alpha_i}{p}, p} ( W^{1,p})  (Q_{\rho^{(i, i+1)}}) }   \le 
2^{i} \frac{C (p, \vfi'' (0))}{R-r} \times \\
 \Big( [\V ( \D u)]_{1, \frac{\delta}{2}, N^{\alpha_i, 2} ( L^{2}  ) (Q_{\rho^{(i, i+1)}}) } +  [u ]_{1, \frac{\delta}{2}, N^{\alpha_i, 2} (  I_{\rho^2}; L^{2} (B_{\rho}) )}\!+  [u ]_{1, \frac{\delta}{2}, N^{\frac{2\alpha_i}{p}, p} (  I_{\rho^2}; L^{p} (B_{\rho}) )} \!+  |u |_V + 1\! \Big),
\end{multline}
where we have increased the powers using $x^\frac{2}{p} \le 1 + x$ for $p \ge 2$. We estimate the first summand at the r.h.s. of \eqref{idea:ind:det:2} with \eqref{idea:ind:det:0} and the next two ones with  \eqref{idea:indAss} as well as the~embedding Lemma \ref{lem:nik:emb}. Hence we obtain the low-time and high-space information 
\begin{equation}\label{idea:ind:det:3}
 |u \eta |_{\frac{\delta}{2}, N^{\frac{2\alpha_i}{p}, p} ( W^{1,p})  (Q_{\rho^{(i, i+1)}}) }   \le 
4^{i} \frac{C (p, \vfi'' (0))}{ (R-r)^2} \left( \left(  1+  |u|_M \right)^\frac{p-2}{2}  |u|_{{\frac{\delta}{2},} N^{\alpha_i, p} (L^p)  (Q_{\rho_i}) } +  |u |_V + 1\!\right).
  \end{equation}
\item[(iii)+(iv)]We use the interpolation Lemma  \ref{lem:jb_strong:nikolski_inter} with the parameters given by \eqref{eq:jb_strong:inter:choice}, $q_0 = p$, \eqref{eq:jb_strong:ivi:pq1},  \eqref{eq:jb_strong:ivi:pq2} for $u \eta$ and obtain 
\begin{equation}\label{idea:ind:det:4}
 |u \eta|_{\frac{\delta}{2}, N^{\alpha' (i), p_0} ( L^{p}) (Q_{\rho^{(i, i+1)}})}  \le  C_1 |u \eta |^{1-\frac{\theta}{2}}_{\frac{\delta}{2}, N^{\frac{2\alpha_i}{p}, p} ( W^{1,p})  (Q_{\rho^{(i, i+1)}}) }   |u \eta|^\frac{\theta}{2}_{\frac{\delta}{2}, N^{1+ \alpha_i, p'} \left( W^{-1, p'} \right)  (Q_{\rho^{(i, i+1)}}) }.
  \end{equation}
  \end{itemize}
We estimate the r.h.s. of \eqref{idea:ind:det:4} with \eqref{idea:ind:det:1}, \eqref{idea:ind:det:3}. Next we use the embedding Lemma \ref{lem:nik:emb} and choice of $\eta$ to obtain the qualitative ($\alpha_{i+1}$)
\begin{multline}\label{idea:ind:det:5}
 |u|_{\frac{\delta}{2}, N^{\alpha_{i+1}, p} ( L^{p}) (Q_{\rho^{ i+1}})}  \le \\
 \frac{  C_1 4^{i} C (p)  }{(\vfi'' (0))^{\beta(p)} (R-r)^2} 
   \left( \left(  1+  |u|_M\right)^{p-2}   |u|_{{\frac{\delta}{2},} N^{\alpha_i, p} (L^p)  (Q_{\rho_i}) } +   |u|_V  \left(  1+  |u|_M\right)^\frac{p-2}{2}   +1 \right)
  \end{multline}
with $ \alpha_{i+1}$ given by \eqref{eq:jb_strong:ivi2} and $\beta (p) =  \frac{2}{p} (1 - \frac{\theta}{2}) =  \frac{2}{p} - \frac{2}{d (p-2) + 2p}$. Hence for any $\alpha < \gamma_0$ after $N (\alpha)$ iterations we obtain 
\begin{equation}\label{idea:ind:det:6}
 |u|_{\frac{\delta}{2}, N^{\alpha, p} ( L^{p}) (Q_{\rho^{ N (\alpha) +1}})}  \le \frac{ C (C_1, p, \alpha)  }{\left((\vfi'' (0))^{\beta( p)} (R-r)^2 \right)^{N(\alpha)}} 
  \left(  1+  |u|_V \right)^{(p-2) N(\alpha)}  \left( |u|_{{ \frac{\delta}{2},} N^{\alpha_0, p} (L^p)  (Q_{R}) } +  1\right).
  \end{equation}
On $Q_{\rho^{ N (\alpha) + 2}}$, the l.h.s. of \eqref{idea:ind:det:6}  controls 
\begin{equation}\label{idea:ind:det:6N}
 |u|_{{\frac{\delta}{2},} N^{\alpha, \infty} (L^2)}  + |\V (\D u)|_{{\frac{\delta}{2},} N^{\alpha, 2} ( L^2 )}  +
 |u|_{\frac{\delta}{2}, N^{1+ \alpha, p'} \left(W^{-1, p'} \right) } +
 |u |^\frac{2}{p}_{\frac{\delta}{2}, N^{\alpha, 2} ( W^{1,2} )} +  |u |_{\frac{\delta}{2}, N^{\alpha, p} (L^{p} )}  + |u |_{\frac{\delta}{2}, N^{\frac{2\alpha}{p}, p} (W^{1,p} )}  
  \end{equation}
through \eqref{eq:jb_strong:2+:t:frac}, \eqref{eq:jb_strong:iter_1:superquad}, \eqref{eq:jb_strong:iter:cor'a}. Recall that $Q_r \subset Q_{\rho^{ N (\alpha) + 2}}$, $\alpha_0 =0 $ and that  $C_1$ depends on $\alpha, p$. Hence \eqref{idea:ind:det:6} via \eqref{idea:ind:det:6N} gives for a certain ${\kappa (p, \alpha)}$
  \begin{multline}\label{idea:ind:det:6N2}
 |u|_{{\frac{\delta}{2},} N^{\alpha, \infty} (L^2)(Q_{r}) }  + |\V (\D u)|_{{\frac{\delta}{2},} N^{\alpha, 2} ( L^2 )(Q_{r})}  +
 |u|_{\frac{\delta}{2}, N^{1+ \alpha, p'} \left(W^{-1, p'} \right) (Q_{r})} +
 |u |_{\frac{\delta}{2}, N^{\alpha, 2} ( W^{1,2} )(Q_{r})} +\\
  |u |_{\frac{\delta}{2}, N^{\alpha, p} (L^{p} )(Q_{r})}  + |u |_{\frac{\delta}{2}, N^{\frac{2\alpha}{p}, p} (W^{1,p} )(Q_{r})}
 \le C (p, \alpha, \vfi'' (0))  \Big( \frac{ 1+ |u|_{L^\infty ( L^{2} )(Q_R)} +  |u |_{L^{ p} (W^{1, p} )(Q_R) } }{R-r}  \Big)^{\kappa (p,\, \alpha)}.
  \end{multline}
Finally, we use \eqref{eq:jb_strong:nik_bas:n:diff:both:delta} that allows to control the standard norm  $|\cdot|_{ N^{\alpha, q} (X)}$ with  $|\cdot|_{{\frac{\delta}{2},} N^{\alpha, q} (X)}$ at the cost of the term $\delta^{-\alpha}$. The choice of $\delta$ depends on the iteration parameters and on the (fixed) distance of the \"ubercylinder $ \ddot{Q}$ from $\Subset \O_I$. Hence the desired \eqref{eq:jb_strong:iter:lem_gen_low} has the form \eqref{idea:ind:det:6N2} (recall that constants that depend on parameters may vary within this dependence).
\end{proof}
Now let us consider the situation when we can gain the full differentiability. This is the case \[\gamma_0 > 1 \quad \iff \quad p < 2 + \frac{2}{\sqrt{d+1}}.\]
\begin{lem}[Full differentiability in time]\label{lem:jb_strong:iter_end:hi}
Take an \"ubercylinder, an untercylinder $ \underaccent{\ddot}{Q} \Subset \ddot{Q} \Subset \O_I$ fixed by the Assumption \ref{rem:uberunter} and a local weak solution $u$ to \eqref{eq:jb_structure:weak} on $\O_I$. Moreover, assume that  $p \in \left[2, 2 + \frac{2}{\sqrt{d+1}}\right)$.  Then for any $\gamma$ such that
\[ \gamma < \gamma_1 \, \left(= \frac{2}{p} + \frac{p}{d (p-2) + 2p} \right)\]
there exist $\kappa (p, \gamma), \; C (p, \gamma, \vfi'' (0))$ such that for all the concentric cylinders $Q_r, Q_R$ satisfying $\underaccent{\ddot}{Q} \subset Q_r  \Subset Q_R \subset \ddot{Q}$ one has
 \begin{multline}\label{eq:jb_strong:iter:lem_gen_hip}
 |u|_{ N^{ \gamma, p} ( L^{p})(Q_{r})}  + |u|_{ W^{2, p'} \left( W^{-1, p'}\right) (Q_{r})} +  |u|_{ W^{1, \infty} ( L^2 )(Q_{r})}  + |\V (\D u)|_{W^{1, 2} (L^2 )(Q_{r})} +\\
  |u |_{W^{1, 2} (W^{1,2} )(Q_{r})} +  |u|_{ N^{\frac{2}{p}, p} (W^{1,p})(Q_{r})} 
 \le
C (p, \gamma, \vfi'' (0))  \left( \frac{  |u|_{L^\infty ( L^{2}) (Q_R)} +  |u |_{L^{ p} (W^{1, p} )(Q_R) } +1}{R-r}  \right)^{\kappa (p, \gamma)}
 \end{multline}
in the case $p>2$ or 
  \begin{multline}\label{eq:jb_strong:iter:lem_gen_hi2}
   |u|_{ N^{ \gamma, 2} ( L^{2})(Q_{r})}  + |u|_{ W^{2, 2} \left( W^{-1, 2}\right) (Q_{r})} +  |u|_{ W^{1, \infty} ( L^2 )(Q_{r})}  + |\V (\D u)|_{W^{1, 2} (L^2 )(Q_{r})} +\\
    |u |_{W^{1, 2} (W^{1,2} )(Q_{r})} \le
 C (\gamma, \vfi'' (0))  \left( \frac{  |u|_{L^\infty ( L^{2} (Q_R))} +  |u |_{L^{ 2} (W^{1, 2} (Q_R)) } +1}{R-r}  \right)^{\kappa (\gamma)}
 \end{multline}
in the case $p=2$.
\end{lem}
\begin{proof} 
Along the proof of Lemma \ref{lem:jb_strong:iter_end:low}, for any $\alpha_{i_0} <1$ we needs a finite number of steps of the iteration \eqref{eq:jb_strong:iter:cor_scheme:bf} to reach
\begin{equation}\label{idea:ind:det:hi:1}
 |u|_{\frac{\delta}{2}, N^{\alpha_{i_0}, p} ( L^{p}) (Q_{\rho^{ N (\alpha_{i_0}) }})}  \le
 C (p, \alpha_{i_0}, \vfi'' (0)) \Big( \frac{ 1+  |u|_{L^\infty ( L^{2} )(Q_R)} +  |u |_{L^{ p} (W^{1, p} )(Q_R) } }{R-r} \Big)^{\kappa (p, \alpha_{i_0})}.
  \end{equation}
 It has been already explained in Subsection \ref{ssec:crone}, that performing one iteration step beyond  the full differentiability is allowed. Let us recall briefly the argument. At the level of any $\alpha_{i_0} <1$ the fractional energy estimates of Section \ref{sec:energyinNik} are valid. Therefore we take in \eqref{eq:jb_strong:iter:cor_scheme:bf} the  iterative assumption ($\alpha_{i_0}$) and proceed with the step  ${i_0} \to {i_0+1} $, because both interpolation and embeddings work without the upper bound $\alpha <1$ for the fractional differentiability. For the quantitative result, we additionally take into account the~possible increase of the constant, caused by this last step ${i_0} \to {i_0+1} $. There is also no problem with second-order differences, as Lemma \ref{lem:jb_strong:nik_bas} allows to increase the order of differences in a Nikolskii-Bochner space. Hence one obtains for any $\gamma < \gamma_1$,  compare \eqref{A+B}  with $1+\eps =\gamma_1$,
\begin{equation}\label{idea:ind:det:hi:2}
\!\! |u|_{\frac{\delta}{2}, N^{\gamma, p} ( L^{p}) (Q_{\rho^{ N (\alpha_{i_0+1})}})}    \le 
 C (p, \gamma, \vfi'' (0))  \Big( \frac{ 1+  |u|_{L^\infty ( L^{2} )(Q_R)} +  |u |_{L^{ p} (W^{1, p} )(Q_R) } }{R-r}  \Big)^{\kappa (p, \gamma)}.
  \end{equation}
  The r.h.s. of \eqref{idea:ind:det:hi:2} controls 
  \[
   |u|_{\frac{\delta}{2}, W^{1, p} ( L^{p}) (Q_{\rho^{ N (\alpha_{i_0+1})}})} 
  \]
  thanks to \eqref{eq:nik:embW} of the embedding Lemma \ref{lem:nik:emb}. It gives thesis via \eqref{eq:jb_strong:2+:t:frac:sob}, \eqref{eq:jb_strong:iter_1:superquad:sob}, \eqref{eq:jb_strong:iter:cor'a:sob}, \eqref{eq:jb_strong:iter:cor'a:sob:p2}. We use \eqref{eq:jb_strong:nik_bas:n:diff:both:delta:sob}, \eqref{eq:jb_strong:nik_bas:n:diff:both:delta} to pass from the norms containing term $ \frac{\delta}{2}$ to the standard norms, where applicable\footnote{Observe that neither \eqref{eq:jb_strong:nik_bas:n:diff:both:delta:sob} nor \eqref{eq:jb_strong:nik_bas:n:diff:both:delta} covers the case $ |u|_{ W^{2, p'} \left( W^{-1, p'}\right) (Q_{r})}$, but we do not need to get rid of $ \frac{\delta}{2}$-term there, because \eqref{eq:jb_strong:iter_1:superquad:sob} concerns already the standard  $\delta$-free norm of $ |u|_{ W^{2, p'} \left( W^{-1, p'}\right) }$.}, and get thesis.
\end{proof}
\subsection{Proof of Theorem \ref{theo:jb_strong:temporal} }
\markboth{\MakeUppercase{Chapter 4. \, Temporal Regularity}}{\MakeUppercase{4.4. \: Proof of Theorem 4.1}}
It is a combination of Lemma \ref{lem:jb_strong:iter_end:hi} and Lemma \ref{lem:jb_strong:iter_end:low}.

\subsection{Proof of Theorem \ref{lem:jb_strong_stat}}
\begin{proof}
We follow closely the proof of Theorem 1 of \cite{BurKap14a}. The main difference lies in the fact that now we do not obtain any regularity from the evolutionary part, but we assume it. This in turn allows us to use an arbitrary $\sigma \in \test(I^\eps) $ in place of the cutoff function $\sigma$ from the proof of Theorem 1 of \cite{BurKap14a}.

Now we proceed with some more details. However, since the approach is very close to \cite{BurKap14a}, we refer to \cite{BurKap14a} or to Chapter  5 of \cite{Burphd} for full explanations. Let us define the space difference
\[
\Delta^{(s, 0)} f = T_{(s,0)} f - f,
\]
where 
\[
( T_{(s,0)} f )(x,t) = f (x+s, t) 
\]
is the space shift.

\underline{Step 1. (Estimate for differences with growing spatial support.)} 
Let us fix $\eps > 0$ so small that the~\"ubercylinder $Q_{R+ 2 \eps} \Subset \O_I$. We take any concentric $Q_{\rho_1} \Subset Q_{\rho_2}$, such that $Q_{\rho_2} \subset Q_R$. A smooth $\eta$ cuts off between $B_{\rho_1}$ and $B_{\rho_2}$ with  $| \nabla \eta |_\infty \le  \frac{2}{{\rho_2}-{\rho_1}}$ and $\sigma$ is an arbitrary function from $ \test( I^\eps)$.
Similarily as in Lemma \ref{cor:jb_strong:intro}, a local weak solution $u$ of \eqref{eq:jb_structure:weak} on $\O_I$ satisfies for real $|l| \le \frac{\eps}{2}$, where $e_i$ denotes $i$-th canonical vector in $\er^d$

\begin{multline}\label{stat:0}
 \int_{I^\eps} \int_{\Omega^\eps } \Delta^{(le_i, 0)}  \! \left(  \A ( \D u  ) \right):   \left( \Delta^{(le_i, 0)} \D u \right)   \eta^2 \sigma^2    =\\
- \int_{I^\eps} \! \langle  (\Delta^{(le_i, 0)}  u)_t  ,  \left( \Delta^{(le_i, 0)}  u \right)    \eta^2 \sigma^2   \rangle_{(W^{1,p}_{0} (\Omega^\eps  ))^*\!, \; W^{1,p}_{0} (\Omega^\eps )} +  2 \int_{I^\eps} \!  \int_{\Omega^\eps }  \Delta^{(le_i, 0)} \! \left(  \A ( \D u  ) \right):  \left(  \Delta^{(le_i, 0)}  u \; \hat \otimes \; \nabla \eta \right) \eta \sigma^2\\
 =: I + II.
\end{multline}
Let us focus on the term $I$ at the r.h.s. of \eqref{stat:0}. Commutativity  of time derivatives with space differences yields
\[
I = \int_{I^\eps} \left(\int_{\Omega^\eps }  \left( \Delta^{(le_i, 0)}   u,_t \right) \cdot   \left( \Delta^{(le_i, 0)}  u \right)    \eta^2 \right)\sigma^2.
\]
Our choice of the cutoff function gives $\supp \eta \subset B_R$. Hence our choice of $Q_{R + 2 \eps}$ gives $\supp \eta \subset \Omega^{2\eps}$. This along with the assumed $ |s| \le \sfrac{\delta }{2}$, after the discrete integration by parts, yields
\[
I = \int_{I^\eps} \left(\int_{\Omega^\eps }  u,_t  \cdot  \Delta^{(-le_i, 0)} \left(\left( \Delta^{(le_i, 0)}   u \right)    \eta^2 \right) \right)\sigma^2 \le 
\int_{I^\eps} \int_{\Omega^\eps }  |u,_t| \left|   \Delta^{(-l e_i, 0)} \left(\left( \Delta^{(le_i, 0)}  u \right)    \eta^2 \right) \right| \sigma^2.
\]
Next we estimate the difference $\Delta^{(-l e_i, 0)} $ with the derivative and obtain
\begin{multline*}
I \le \int_{I^\eps} \int_{\Omega^\eps }  |u,_t| |l| \left( \dashint_0^l | \left(\left( \Delta^{(le_i, 0)}  u \right)    \eta^2\right)\!,_{x_i}  \! \circ \, T_{(-\lambda e_i, 0)}|d \lambda \right)  \sigma^2 \le \\
 \delta  \int_{I^\eps}  \dashint_0^l   \int_{\Omega^\eps} | \left(\left( \Delta^{(le_i, 0)}  u \right)    \eta^2\right)\!,_{x_i}  \! \circ \, T_{(-\lambda e_i, 0)}|^2   d \lambda \;  \sigma^2 + C (\delta) |l|^2 \int_{I^\eps} \int_{B_{R +  \eps}}  |u,_t|^2  \sigma^2.
\end{multline*}
In the above estimate  we have also used  that  the spatial support of the integrand of the last term is contained in $ B_{R+ \eps} \subset \Omega^{\eps}$. It follows from $\supp \eta \subset B_{\rho_2} \subset B_R$. This last information and the~Tonelli Theorem allow us to write\vspace{-1mm}
\begin{multline}\label{stat:05}
I \le \delta  \int_{I^\eps}  \int_{\Omega^\eps }  |  \nabla \Delta^{(le_i, 0)}   u|^2      |\eta|^2 \sigma^2 + \delta  \int_{I^\eps}  \int_{\Omega^\eps }  | \Delta^{(le_i, 0)}  u |^2   |\nabla \eta|^2 \sigma^2 + C (\delta)  |l|^2 \int_{I^\eps} \int_{B_{R +  \eps}}  |u,_t \!|^2  \sigma^2  \le \\
\delta  \int_{I^\eps}  \int_{B_{\rho_2} }  | \nabla \Delta^{(le_i, 0)}    u|^2  \sigma^2  + \frac{ \delta }{({\rho_2}-{\rho_1})^2} \int_{I^\eps}  \int_{B_{\rho_2} }  | \Delta^{(le_i, 0)}  u |^2   \sigma^2  + C (\delta)   |h|^2  \int_{I^\eps} \int_{B_{R +  \eps}}  |u,_t \!|^2  \sigma^2.
\end{multline}
For the second inequality in \eqref{stat:05} we have used characteristics of $\eta$ and we have increased $l $ to $h$ such that $\frac{\eps}{2} \ge h \ge l$. This $h$ will be used later as the proper denominator of difference quotients. 
The~Korn inequality gives us subsequently 
\begin{multline}\label{stat:1}
I \le \delta  \int_{I^\eps}  \int_{B_{\rho_2} }  | \D \Delta^{(le_i, 0)}    u|^2  \sigma^2  + \frac{ C \delta }{({\rho_2}-{\rho_1})^2} \int_{I^\eps}  \int_{B_{\rho_2} }  | \Delta^{(le_i, 0)}  u |^2   \sigma^2  + C (\delta)   |h|^2  \int_{I^\eps} \int_{B_{R +  \eps}}  |u,_t \!|^2  \sigma^2  \\
\le \delta  \int_{I^\eps}  \int_{B_{\rho_2} }  | \Delta^{(le_i, 0)}   \D  u|^2  \sigma^2  +   \frac{C (\delta)   |h|^2  }{({\rho_2}-{\rho_1})^2}  \int_{I^\eps} \int_{B_{R +  \eps}}  |u,_t \!|^2  \sigma^2 +   | \nabla  u |^2 \sigma^2.
\end{multline}
For the second inequality in \eqref{stat:1} we have estimated the difference $| \Delta^{(le_i, 0)}  u |^2$ with the space gradient.

Let us now focus on the term $II$ on the r.h.s. of  \eqref{stat:0}. Using Assumption \ref{ass:struc} on growth of $\A$, we have   \begin{multline}\label{stat:2}
II \le 4 \delta  \int_{I^\eps}   \int_{B_{\rho_2}}  |\Delta^{(l e_i, 0)}  \V ( \D u) |^2   \sigma^2  + 6 \delta  \int_{I^\eps}  \frac{|l|}{|h|}  \dashint_0^l     \intop_{B_{\rho_2}}   |\Delta^{(\lambda e_i, 0)}  \V ( \D u) |^2 d \lambda  \; \sigma^2  \\
   +  
C (\delta, p)\Big(1 + \frac{1}{\vfi'' (0)}\Big) \frac{|h|^2}{({\rho_2}-{\rho_1})^2}    \int_{I^\eps} \int_{B_{R+ \eps}}    \vfi \left(  | \nabla u  | \right)    \sigma^2.
\end{multline}
For the time being, \eqref{stat:2} holds  for $|h| \le h_0 = {\rho_2}-{\rho_1}$.
Simultaneously, in view of the pointwise estimate 
\[
\Delta^{(s,0)} \left(  \A ( \D u  ) \right)\!:\! ( \Delta^{(s,0)} \D u )    \ge \frac{1}{C(G(\vfi'))}  | \Delta^{(s, 0)} \V ( \D u ) |^2  +  \frac{\vfi'' (0)}{C}  | \Delta^{(le_i, 0)} \D u |^2
\]

one has
\begin{multline}\label{stat:3}
 \int_{I^\eps}   \int_{B_{\rho_1}} | \Delta^{(le_i, 0)}  \V ( \D u ) |^2  \sigma^2 +  \vfi'' (0)  \int_{I^\eps}   \int_{B_{\rho_1}}  \left| \Delta^{(le_i, 0)} \D u \right|^2 \sigma^2 \le \\
  C(p) \int_{I^\eps}   \int_{B_{\rho_1}} \Delta^{(le_i, 0)}  \left(  \A ( \D u  ) \right)\!:\!  \left( \Delta^{(le_i, 0)}  \D u \right)  \sigma^2.
\end{multline}
 We put together  \eqref{stat:1},  \eqref{stat:2} and  \eqref{stat:3} with the use of \eqref{stat:0} to get
\begin{multline}\label{stat:4}
\int_{I^\eps}   \int_{B_{\rho_1}} | \Delta^{(le_i, 0)}  \V ( \D u ) |^2  \sigma^2 +  \vfi'' (0)   \int_{I^\eps}   \int_{B_{\rho_1}}  \left| \Delta^{(le_i, 0)} \D u \right|^2 \sigma^2 \le \\
 \delta_0  \int_{I^\eps}  \int_{B_{\rho_2} }  |  \Delta^{(le_i, 0)}   \D u|^2  \sigma^2  + 4 \delta_0  \int_{I^\eps}   \int_{B_{\rho_2}}  |\Delta^{(l e_i, 0)}  \V ( \D u) |^2   \sigma^2  + 6 \delta_0  \int_{I^\eps}  \frac{|l|}{|h|}  \dashint_0^l     \int_{B_{\rho_2}}   |\Delta^{(\lambda e_i, 0)}  \V ( \D u) |^2 d \lambda   \sigma^2  \\
   +  
C (\delta_0, p) \left(1 + \frac{1}{\vfi'' (0)} \right) \frac{|h|^2}{({\rho_2}-{\rho_1})^2}    \int_{I^\eps} \int_{B_{R+ \eps}}  \left(  \vfi \left(  | \nabla u  | \right)    \sigma^2 +   |u,_t\!|^2  \sigma^2 \right),
\end{multline}
where we have chosen $\delta_0 : = \frac{ \delta}{C(p)}$. 

\underline{Step 2. (Homogenization and uniformization.)} We apply $\dashint_0^h d l$ to \eqref{stat:4} and obtain
\begin{multline}\label{stat:5}
\int_{I^\eps}  \dashint_0^h \int_{B_{\rho_1}} | \Delta^{(le_i, 0)}  \V ( \D u ) |^2  dl \sigma^2 +  \vfi'' (0)   \int_{I^\eps}  \dashint_0^h  \int_{B_{\rho_1}}  \left| \Delta^{(le_i, 0)} \D u \right|^2 dl \sigma^2 \le \\
 \delta_0  \int_{I^\eps} \dashint_0^h  \int_{B_{\rho_2} }  |  \Delta^{(le_i, 0)}   \D u|^2 dl \sigma^2  + 10 \delta_0  \int_{I^\eps} \dashint_0^h   \int_{B_{\rho_2}}  |\Delta^{(l e_i, 0)}  \V ( \D u) |^2 dl  \sigma^2   \\
   +  
C (\delta_0, p) \Big(1 + \frac{1}{\vfi'' (0)}\Big) \frac{|h|^2}{({\rho_2}-{\rho_1})^2}    \int_{I^\eps} \int_{B_{R+ \eps}} \left(   \vfi \left(  | \nabla u  | \right)    \sigma^2 +   |u,_t \!|^2  \sigma^2 \right). 
\end{multline}
Next, we choose $l:= h$ in \eqref{stat:4}  and add the result to \eqref{stat:5}. The obtained estimate works only for $|h| \le {\rho_2}-{\rho_1}$. But the case $|h| >{\rho_2}-{\rho_1}$, due to large $h$'s, is in fact harmless (for details, subsection 3.2.3 of \cite{Burphd}). Consequently it holds for any $|h| \in (0, \frac{\eps}{2} ]$ and $\rho_1 < \rho_2 \le R$. Finally, the choice $\delta_0 = \frac{1}{22}$ in it gives for 
\begin{multline*}
f^{h,\sigma} (\rho) : =\\
  \int_{I^\eps} \left[  \int_{B_{\rho}} \! \left( \left| \Delta^{(he_i, 0)} \V(\D u) \right|^2 \!+\!  \vfi'' (0) \left| \Delta^{(he_i, 0)} \D u \right|^2\!+  \dashint_0^h \left( \left| \Delta^{(\lambda e_i, 0)} \V(\D u) \right|^2 \!+\!  \vfi'' (0)  \left| \Delta^{( \lambda e_i, 0)} \D u \right|^2    \right) d \lambda \!  \right) \! \right] \!\sigma^2
\end{multline*}\vspace{-3mm}
that
\[
f^{h,\sigma} (\rho_1) \le  \frac{1}{2} f^{h,\sigma} (\rho_2) +  C ( p) \Big(1 + \frac{1}{\vfi'' (0)} \Big) \frac{|h|^2}{({\rho_2}-{\rho_1})^2}    \int_{I^\eps} \left[ \int_{B_{R+ \eps}}    \vfi \left(  | \nabla u  | \right)+   |u,_t\!|^2 \right] \sigma^2
\]
for any $h$ such that $|h| \in (0, \frac{\eps}{2} ]$, any $\sigma \in \test (I^\eps)$ and any $0 < \rho_1 < \rho_2 \le R$.

\underline{Step 3. (Conclusion via a DuBois-Reymond-type argument.)} 
For a fixed $h$, $\sigma$, $f^{h,\sigma} (\rho)\!: [0, R] \to \er_+$  is non-decreasing and bounded. Hence vie the Giaquinta-Modica argument one has
\[
f^{h,\sigma} (r) \le \Big( 1 + \frac{1}{\vfi'' (0)}  \Big) \frac{C ( p)  |h|^2}{(R-r)^2}    \int_{I^\eps} \left[ \int_{B_{R+ \eps}}  (  \vfi \left(  | \nabla u  | \right)+   |u,_t \!|^2 ) \right] \sigma^2
\]
for any $r \le R$. After dropping the homogenization terms in $f^{h,\sigma} (r) $ we arrive at
\begin{multline}\label{stat:6}
 \int_{I^\eps}  \int_{B_{r}} \! \left( \left| \Delta^{(he_i, 0)} \V(\D u) \right|^2 +  \vfi'' (0) \left| \Delta^{(he_i, 0)} \D u \right|^2 \right) \sigma^2 \\
 \le  \int_{I^\eps}  \Big( 1 + \frac{1}{\vfi'' (0)}  \Big) \frac{C ( p) |h|^2}{(R-r)^2}   \int_{B_{R+ \eps}}  \left(    \vfi \left(  | \nabla u  | \right)+   |u,_t\!|^2 \right) \sigma^2
\end{multline}
for any $h$ such that $|h| \in (0, \frac{\eps}{2} ]$, any $\sigma \in \test (I^\eps)$. For a fixed, positive  $h$, the estimate \eqref{stat:6} can be written as
\[
 \int_{I^\eps}  (A_h- B_h) (t) \psi (t) \le 0,
\]
holding for any smooth nonnegative $\psi \in \test (I^\eps)$. The term $ B_h$ belongs to  $L^\infty (I^\eps) $ thanks to our assumptions. Similarly $A_h \in L^\infty (I^\eps)$, where we use Assumption \ref{ass:struc} and $p \ge 2$. The arbitrariness of a nonnegative $\psi \in \test (I^\eps)$ implies therefore that for a.e. $t \in I$
\[
-\infty <-B_h (t) \le A_h(t)-B_h(t)  \le 0.
\]
Hence for a.e. $t \in I$ and $h$ such that $|h| \in (0, \frac{\eps}{2} ]$ we get
\[
\int_{B_{r}} \! \left( \left| \frac{\Delta^{(he_i, 0)} \V(\D u)}{h} \right|^2 \!+\!  \vfi'' (0) \left| \frac{\Delta^{(he_i, 0)} \D u}{h} \right|^2 \right) (t)  \! \le  \Big( 1 + \frac{1}{\vfi'' (0)} \Big) \frac{C ( p) }{(R-r)^2}   \int_{B_{R+ \eps}}    \vfi \left(  | \nabla u (t) | \right)+   |u (t),_t \!|^2,
\]
which gives our thesis via the difference quotient argument.
\end{proof}

\subsection{Proofs of Lemmas \ref{for:jb_strong:low}, \ref{for:jb_strong:hi}}
Lemma \ref{for:jb_strong:low} is a qualitative version of  Theorem \ref{theo:jb_strong:temporal}.

\begin{proof}[Proof of Lemma \ref{for:jb_strong:hi}]
 Theorem \ref{theo:jb_strong:temporal} (ii), (iii) implies that $u$ is in $N^{\gamma, p} (L^{p}), { W^{2, p'} (W^{-1, p'} )}, { W^{1, \infty} (L^2)}, {W^{1, 2} (W^{1,2})}$ and, for  $p>2$, $N^{\frac{2}{p}, p} ( W^{1,p})$. The information that $u$ comes from ${ W^{1, \infty} (L^2)}$ and 
 \begin{itemize}
 \item $N^{\frac{2}{p}, p} ( W^{1,p}) \csubset L^\infty ( W^{1,p})$ in the case   $p>2$,
 \item ${W^{1, 2} (W^{1,2})} \csubset   L^\infty ( W^{1,2}) $ in the case   $p=2$
 \end{itemize}
 allows us to use Theorem \ref{lem:jb_strong_stat} and obtain membership of $u$ in $L^\infty ( W^{2,2})$ and of $ \V (\D u)$ in $L^\infty ( W^{1,2})$. This, via Sobolev embedding, gives $u \in L^\infty ( W^{1,\frac{2^*p}{2} })$.
 \end{proof}

\section{Appendix}\label{sec:app}
Here we provide needed results on Nikolskii-Bochner spaces and some comments

In order to put the Nikolskii-Bochner spaces in a broader analytic context, let us briefly mention Besov-Bochner spaces. We introduce them by the real interpolation, along  \cite{Ama00} of Amann. By $(\cdot, \cdot)_{\theta, q}$ we understand the real interpolation functor. For $k >0$ we fix natural $k_1, k_2$ such that $k_1 < k < k_2$ and put $\theta = \sfrac{k-k_1}{k_2-k_1}$. Any intermediate space \[(W^{k_1,p} ( I; X); W^{k_2,p}  ( I; X))_{\theta, q}\] is the same one; we call it the Besov-Bochner space  $B_q^{k, p} ( I; X)$ (see \cite{Ama00}, especially Section 4 there). The Nikolskii-Bochner space with a non-integer $\alpha$ can now be identified with a Besov-Bochner space with $q=\infty$. Indeed, formula (6.1) of \cite{Ama00}, Section 2 gives \[N^{\alpha, p} ( I; X) \simeq B_\infty^{\alpha, p} ( I; X)\] for $\alpha \notin \N$. Observe that in fact Nikolskii-Bochner spaces are defined in \cite{Ama00} in a slightly different way then in our case. Instead of higher-order differences of \eqref{eq:jb_strong:nik_semi:gen}, Amann uses in his definition differences of derivatives. As Besov spaces admit numerous equivalent descriptions, these definitions are equivalent\footnote{For references, see the remainder of this subsection. For rigorous results in the case interesting for us, see subsection \ref{subs:diff_der_eq}}.

\subsection{References}\label{bes:ref} The literature on standard Besov spaces (as opposed to the Besov-Bochner case) is broad. In this case one can trace back the theory to the original works of Nikolskii \cite{Nik51} and Besov \cite{Bes59}\footnote{The Nikolskii's paper is older. It reflects the fact, that Nikolskii spaces are the most straightforward generalization of  Sobolev spaces over the fractional differentiabilities, within the Besov scale.}, with important earlier contributions by Marchaud \cite{Mar27} and Zygmund \cite{Zyg45}. The~classical references for function spaces are the books by Triebel, for instance \cite{Tri83}. For anisotropic Besov spaces, useful for the evolutionary problems, one refers to Besov, Illin, Nikolskii \cite{BesIliNik75}.

 One can define a Besov space, using equivalently:
\begin{itemize}
\item interpolation, see \cite{AdaFou03}, Section 7.30 and the following ones,
 \item differences (that corresponds to our Definition \ref{def:nik}), see \cite{Tri83}, Subsection 2.2.2,
 \item a mixture of differences and weak derivatives, see \cite{Gol62},
 \item  moduli of continuity (see \cite{BenSha88}, Chapter 5, Section 4),
 \item the Fourier transform (see \cite{Tri83}, section 2.3). 
 \end{itemize}
 This versatility carries over the Besov-Bochner spaces (which are commonly referred to in literature as the vector-valued Besov spaces), but the results are less clearly stated and scattered, compare \cite{Ama97},  \cite{Mur74}, \cite{SchSchSic12}. This forces us to prove the needed results by ourselves.
 \subsection{Results on  Nikolskii-Bochner  spaces}
 Here we gather needed results. In order to allow this section to serve as a quick reference, the proofs of these results are gathered in subsection \ref{sec:nikpfs}.
 \subsubsection{Basic properties of Nikolskii-Bochner  spaces}
The symbol $[ r ]$ denotes here entire value of a real number $r$ and should not be confused with the notation for a seminorm. We have
\begin{lem}\label{lem:jb_strong:nik_bas}
Take $\alpha > 0$, $p \in [1, \infty]$ and  $r_0 = \min \{n \in \N | \; n > [\alpha] \}$. $N^{\alpha, p} ( I; X) $ endowed with the norm 
\begin{equation}\label{eq:jb_strong:nik_bas:n:d}
|f|_{N^{\alpha, p} ( I; X)} :=  [f]_{{r_0, 1,} N^{\alpha, p} ( I; X)}  + |f|_{L^{ p} ( I; X)}
\end{equation}
becomes a Banach space. Moreover, for any $r \ge r_0$, $\delta \in (0, 1]$ the equivalent norm is
\begin{equation}\label{eq:jb_strong:nik_bas:n:diff}
|f|_{r, \delta, N^{\alpha, p} ( I; X)} :=  [f]_{{r, \delta,} N^{\alpha, p} ( I; X)}  + |f|_{L^{ p} ( I; X)}.
 \end{equation}
More precisely, for a fixed $r \ge r_0$ and any  $\delta_1 \le \delta_2 $ from $(0, 1]$
 \begin{equation}\label{eq:jb_strong:nik_bas:n:diff:both:delta}
 |f|_{r,\delta_1, N^{\alpha, p} ( I; X)}   \le    |f|_{r,\delta_2,  N^{\alpha, p} ( I; X)}   \le  \frac{3^{r }}{\delta_1^\alpha}   |f|_{r,\delta_1, N^{\alpha, p} ( I; X)}  
 \end{equation}
and for any $r \ge r_0$ and a fixed  $\delta \in (0, 1]$
 \begin{equation}\label{eq:jb_strong:nik_bas:n:diff:both:r}
2^{ r_0 - r} |f|_{r,\delta, N^{\alpha, p} ( I; X)}   \le    |f|_{ r_0 ,\delta, N^{\alpha, p} ( I; X)}   \le \Big(  \frac{4 r^2}{r_0 - \alpha} \frac{1}{\delta^\alpha} \Big)^{r- r_0} |f|_{r,\delta, N^{\alpha, p} ( I; X)},
 \end{equation}
 where the right inequality holds for  
  \begin{equation}\label{eq:jb_strong:smallness}
       \delta \le \frac{|I|}{r 2^{ \left[ \frac{1}{2 (r_0 - \alpha)} \right] + 2}}.
  \end{equation}
\end{lem}

\subsubsection{Equivalences between differences and derivatives}\label{subs:diff_der_eq}
Now let us compare our Definition \ref{def:nik} of a Nikolskii-Bochner space (the definition by differences) with a definition that involves a~mixture of weak derivatives and differences. This is a result that corresponds to one of the equivalences mentioned in the subsection \ref{bes:ref}. 
Recall that we work all the time within Remark \ref{rem:uberunter}, so the length of $I$, present implicitly in the estimates, is bounded from below and from above.  We denote the $\beta$-th order weak time derivative by $f,_t^{(\beta)}$.

The fact that derivatives control differences reads
\begin{prop}[Reduction]\label{lem:jb_strong:nikolski_red}
For $r, \beta ,j \in \N$,  $\alpha \in \er$ such that $r> \alpha \ge 0$
\begin{equation}\label{eq:redNik}
[f]_{{r+\beta, \delta,} N^{\beta+\alpha, p} ( I; X)} \le \left[f,_t^{(\beta)}\! \right]_{{r, \delta,} N^{\alpha, p} ( I; X)},
\end{equation}
provided the r.h.s. is meaningful. Hence in particular
\[
f,_t \, \in N^{\alpha, p} ( I; X) \; \wedge \; f \in L^{ p} ( I; X) \implies f \in N^{1+\alpha, p} ( I; X).
\]
\end{prop}
We call Proposition \ref{lem:jb_strong:nikolski_red} the reduction proposition, because it says that  derivatives  can be reduced to differences.

In order to provide the result converse to the reduction\footnote{Hence \emph{accession}.} Propositon \ref{lem:jb_strong:nikolski_red}, we resort to an extension of a function from $N^{\alpha, p} ( I; X)$ to an element of $ N^{\alpha, p} ( \er; X) $ and the interchangeability of derivatives with differences in general Besov-type spaces. After having invoked these general results we can also immediately reprove and generalize Propositon \ref{lem:jb_strong:nikolski_red}. We keep it however for the~shortness of its proof. The reverse of Proposition \ref{lem:jb_strong:nikolski_red} reads
\begin{prop}[Accession]\label{lem:jb_strong:nikolski_ared}
Take $r, \beta \in \N$, $\alpha \in \er$. For any $r > \alpha> \beta \ge 0$ one can access differences to derivatives as follows
\begin{equation}\label{eq:accNik}
\left| f,_t^{(\beta)} \right|_{{r - \beta,} \delta, N^{\alpha-\beta, p} ( I; X)}  \le  \frac{C (\alpha - \beta)  6^{r }}{\delta^\alpha}  | f|_{{r, \delta,} N^{\alpha, p} ( I; X)},
\end{equation}
provided the r.h.s. is meaningful. Hence, in particular
\[
 f \in N^{1+\alpha, p} ( I; X)  \implies f,_t \in N^{\alpha, p} ( I; X).
\]
\end{prop}
\subsubsection{Interpolation}
Besov spaces are stable under the real interpolation. In case of Nikolskii spaces, their interpolation actually reduces to the H\"older inequality. More precisely,
we have
\begin{prop}\label{prop:abstractNikInter}
Assume that  norms of Banach spaces $Z, X, Y$ satify
\begin{equation}
|g|_Z \le C_{Z,X,Y} |g|^{1-b}_X |g|^b_Y, 
\end{equation}
with $b \in [0,1]$. Then we have for any $\alpha_1, \alpha_2 \ge 0$ and any $r >  \alpha_1 \vee \alpha_2 $,   $\delta > 0$
\begin{equation}
[f]_{r, \delta, N^{\alpha_b, p_b} ( I; Z)} \le C_{Z,X,Y} [f]^{1-b}_{r, \delta, N^{\alpha_1, p_1} ( I; X)} [f]^b_{r, \delta, N^{\alpha_2, p_2} ( I; Y)},
\end{equation}
where 
\vspace{-2mm}
\[
\alpha_b:= (1-b) \alpha_1 + b \alpha_2, \qquad \frac{1}{p_b} = \frac{1-b}{p_1} + \frac{b}{p_2}.
\]
\end{prop}
Let us now derive from Proposition \ref{prop:abstractNikInter} an interpolation result, that is tuned to our further needs (and  not the most general possible). From now on we use the notation
\[
 [f]_{\delta, N^{\alpha, p} ( I; X)}  : =  [f]_{r_0, \delta, N^{\alpha, p} ( I; X)}, \qquad  |f|_{\delta, N^{\alpha, p} ( I; X)}  : =  |f|_{r_0, \delta, N^{\alpha, p} ( I; X)},
\]
where, as before,  $ r_0 = [\alpha] +1 $ is the >>natural<< Nikolskii-Bochner difference step.
\begin{lem} \label{lem:jb_strong:nikolski_inter}
Assume that we have the parameters 
\[p_1, \,q_1, \,p_2, \,q_2 \in (1, \infty) \quad \text{ such that } \quad q_1 \ge q'_2, \quad p_1 \ge p'_2 \quad \text{ and } \quad \alpha_1  \in [ 0, 1 ), \quad \alpha_2  \in [1,2).\]
Take $\theta \in [0,1]$ and \[f \in N^{\alpha_1, p_1} (I; W_0^{1, q_1} (Q)) \cap N^{\alpha_2, p_2} (I; W^{-1, q_2}(Q)).\] Then for 
\vspace{-2mm}
  \begin{equation}\label{eq:jb_strong:smallness'}
       \delta \le \frac{|I|}{2^{ \left[ \frac{1}{2 ( [\alpha] - \alpha +1)} \right] + 3}} 
  \end{equation}
  one has
\begin{equation}\label{eq:jb_strong:nikolski_inter}
 |f|_{\delta, N^{\alpha, p_0} ( I; L^{q_0} (Q)) }  \le  C \Big( \frac{1}{\left( [\alpha] - \alpha +1 \right) \delta^\alpha} \Big)^{1 - [\alpha]}  |f|^{1-\frac{\theta}{2}}_{\delta, N^{ \alpha_1, p_1} ( I; W_0^{1, q_1} (Q)) }  |f|^\frac{\theta}{2}_{\delta, N^{ \alpha_2, p_2} ( I; W^{-1, q_2} (Q)) }
\end{equation}
with
\begin{subequations}\label{eq:jb_strong:nikolski_inter:para}
\begin{align}
\alpha &= \frac{\theta}{2} \alpha_2 + \Big(1 - \frac{\theta}{2} \Big) \alpha_1, \\
 p_0 &=  \frac{2 p_1}{\theta p_1  + 2 (1- \theta)} \label{eq:jb_strong:nikolski_inter:parab}, \\
  q_0 &= \begin{cases} \tilde q_0:= \frac{2 d q_1}{ q_1 \theta d +  (1- \theta) (2 d - 2 q_1)} &\text{ when } \tilde q_0 \in (1, \infty ), \\   \text{ any finite number } &\text{ otherwise. } \end{cases} \label{eq:jb_strong:nikolski_inter:parac}
  \end{align}
\end{subequations}
\end{lem}

The formulation of Lemma \ref{lem:jb_strong:nikolski_inter} may appear slightly awkward. For instance, one would expect interpolation parameter from $[0, 1]$ and we have it effectively in $[0, \frac{1}{2}]$  as well as some presence of $p_2$ in (\ref{eq:jb_strong:nikolski_inter:parab}). This, together with an unnatural (from the viewpoint of interpolation theory) requirement  $q_1 \ge q'_2$, $p_1 \ge p'_2$, is caused by the fact that we deal with the negative differentiability $W^{-1, q_2}$  by a duality formula and next we use a standard interpolation. None of these suboptimalities affect our further results. To prove them  we use the case   $q_1 = q'_2$, $p_1 = p'_2$, where $p_0, q_0$ coincide with the expected interpolation values. Of course our energy estimates depend on peculiarities of Lemma \ref{lem:jb_strong:nikolski_inter}, but regularity classes remain intact.

  \subsubsection{Embeddings}
As before, in the following embedding result we do not write explicitly the~dependence of constants on size of the underlying domain $ I$. In view of Assumption \ref{rem:uberunter} this dependence is irrelevant. We have
 \begin{lem}\label{lem:nik:emb} Fix a non-empty unterinterval $\underaccent{\ddot}{I}$ and an  \"uberinterval $\ddot{I}$.  Take any interval $I$ such that $\underaccent{\ddot}{I} \subset I \subset \ddot{I}$, a Banach space $X$ and $p \in [1, \infty]$,   $\alpha \in \er_+$. \\
 
(i) \; For any $\gamma \in (0,1)$ there is an embedding $N^{[\alpha]+\gamma, p} ( I; X) \csubset W^{[\alpha], p} ( I; X)$ with 
\begin{equation}\label{eq:nik:embW}
\left| f \right|_{ W^{[\alpha], p} ( I; X)} \le  \frac{C\,  [\alpha]  \, 6^{[\alpha] }}{ \gamma \, \delta^{[\alpha] + \gamma}}   | f|_{{ \delta,} N^{[\alpha] + \gamma, p} ( I; X)}.
\end{equation}

(ii) \; Consider the Nikolskii-Bochner spaces $N^{\alpha, p} ( I; X)$, $N^{\alpha', q} ( I; X) $ such that
\[
p, q \in [1, \infty], \quad  \alpha \ge \alpha' \ge 0 , \quad \beta :=  \alpha - \frac{1}{p} - \Big(\alpha' - \frac{1}{q} \Big) >0.
\]
Then there is the continuous embedding into $N^{\alpha, p} ( I; X)  \csubset N^{\alpha', q} ( I; X) $. More precisely,
  \begin{equation}\label{lem:nik:emb:th:low}
 |f|_{\delta, N^{\alpha', q} ( I; X) } \le C_\eqref{lem:nik:emb:th:low} (\alpha, \alpha', \beta, \delta)  |f|_{ \delta, N^{\alpha, p} ( I; X) },
 \end{equation}
provided
\[
\delta \le    \frac{|I|}{ 2^{ \left[ \frac{\tilde C_\eqref{lem:nik:emb:th:low} (\alpha, \alpha')}{ \beta} \right] + 3}} \wedge 1
\]
where
 \[
C_\eqref{lem:nik:emb:th:low} (\alpha, \alpha', \beta, \delta) :=            \begin{cases}
      \left(    \left(  \frac{(\alpha - \alpha'+2)^2}{ ([ \alpha'] +1 -  \alpha')^3} \vee \frac{1}{\alpha -  \alpha'}  \right)  \frac{C 6^{\alpha}}{((\beta \wedge (\alpha - \alpha'))^2 \wedge 1) \delta^{1+\alpha}} \right)^{\alpha - \alpha' + 2}   \quad\quad\quad\quad \text{for} \;\;\, \alpha > \alpha', \; \alpha = [\alpha],\\
 \left(     \left(  \frac{(\alpha - \alpha'+2)^2}{ \left( ([ \alpha'] +1 -  \alpha')  (\alpha - [\alpha] )\right)^3} \vee \frac{1}{\alpha -  \alpha'} \right)   \frac{C 6^{\alpha}}{((\beta \wedge (\alpha - \alpha'))^2 \wedge 1) \delta^{1+\alpha}} \right)^{\alpha - \alpha' + 2}  \quad \text{for} \; \;\, \alpha > \alpha', \alpha > [\alpha],\\
 C   \qquad  \qquad \qquad \qquad \qquad \qquad \qquad \qquad \qquad  \quad  \qquad \; \quad \text{for} \quad \alpha = \alpha'
 \end{cases}
\]
   \[
 \tilde  C_\eqref{lem:nik:emb:th:low}(\alpha, \alpha') =
           \begin{cases}
            \frac{\alpha - \alpha'+2}{[ \alpha'] +1 -  \alpha'} \quad\quad\quad\quad \text{for} \quad \alpha = [\alpha],\\
        \frac{\alpha - \alpha'+2}{ ([ \alpha'] +1 -  \alpha')  (\alpha - [\alpha] )} \quad \text{for} \quad  \alpha > [\alpha].
 \end{cases}
 \]
 
(iii) \; Take the Nikolskii-Bochner space $N^{\alpha, p} ( I; X)$ with $\alpha \in ( 0, 1)$. If $\alpha - \frac{1}{p}> 0 $, then every element of $N^{\alpha, p} ( I; X)$ has  a pointwise representative that belongs to the space of H\"older continuous functions $C^{ 0, \; \alpha - \frac{1}{p}} ( I; X) $ and the estimate
  \begin{equation}\label{lem:nik:embC}
 |f|_{C^{\alpha - \frac{1}{p}} ( I; X) } \le  \frac{3}{\delta^\alpha}  |f|_{\delta, N^{\alpha, p} ( I; X) } 
 \end{equation}
 holds.
\end{lem}

 \subsubsection{Sobolev-Bochner space $W^{1,p}(I; X)$} The space $N^{1,p}(I; X)$ is defined with the second order differences. As we use only the first order estimates in the energy estimates of Section \ref{sec:energyinNik}, let us mention here the standard Sobolev-Bochner space $W^{1,p}(I; X)$.
We have
 \begin{cor}\label{cor65} Take a Banach space $X$ and $p \in [1, \infty]$. 
 The equivalent norm in $W^{1,p}(I; X)$ is
  \begin{equation}\label{eq:jb_strong:nik_bas:n:diff:sob}
|f|_{1, W^{1, p} ( I; X)} :=  [f]_{{1, 1,} N^{1, p} ( I; X)}  + |f|_{L^{ p} ( I; X)}.
 \end{equation}
We have for any  $\delta_1 \le \delta_2 $ from $(0, 1]$
 \begin{equation}\label{eq:jb_strong:nik_bas:n:diff:both:delta:sob}
 |f|_{\delta_1, W^{1, p} ( I; X)}   \le    |f|_{\delta_2,  W^{1, p} ( I; X)}   \le  \frac{3}{\delta_1}   |f|_{\delta_1, W^{1, p} ( I; X)}.  
 \end{equation}
 In the case $p=\infty$ any element of $W^{1, p} ( I; X) $ has a  $C^{0, 1} ( I; X)$ (Lipschitz) representative and there is equivalence of norms. 
 \end{cor}

\subsection{Proofs} \label{sec:nikpfs}
\begin{proof}[Proof of Lemma \ref{lem:jb_strong:nik_bas}]

\underline{Step 1.} We show that ${N^{\alpha, p} ( I; X)} $ is a Banach space by using the fact that ${L^{ p} ( I; X)} $ is a Banach space, compare Proposition 23.2 in \cite{Zei90}. More precisely, we take a Cauchy sequence $\{ f_n \}_{n \in \N}$ in ${N^{\alpha, p} ( I; X)} $. The norm of ${N^{\alpha, p} ( I; X)} $ contains the ${L^{ p} ( I; X)} $-norm, so we have $f$ such that $\lim_{n \to \infty} f_n = f \; \text{ in } \;{L^{ p} ( I; X)}$.  As  $\{ f_n \}_{n \in \N}$ is the  Cauchy sequence in  ${N^{\alpha, p} ( I; X)} $, we have
\[
 h^{- \alpha} | \Delta_h^r f_n  |_{L^{p} ( I_{rh}; X)} \le C
\]
for any ${n \in \N}$. For any  $h \in ( 0, 1) $ one finds $n_h$ such that
\[
| \Delta_h^r (f_{n_h} - f) |_{L^{p} ( I_{rh}; X)} \le 2^r | f_{n_h} - f |_{L^{p} ( I; X)} \le  h^\alpha
\]
the above two inequalities give
\[
 h^{- \alpha} | \Delta_{h}^r f |_{L^{p} ( I_{rh}; X)}  \le  h^{- \alpha} | \Delta_h^r (f_{n_h} - f) |_{L^{p} ( I_{rh}; X)}  +  h^{- \alpha} | \Delta_h^r f_{n_h}  |_{L^{p} ( I_{rh}; X)} \le 1 + C,
 \]
 so $f$ is indeed in ${N^{\alpha, p} ( I; X)} $.
 
\underline{Step 2.} Now we show \eqref{eq:jb_strong:nik_bas:n:diff:both:delta}. The important part of the semi norm $  [f]_{{r, \delta,} N^{\alpha, p} ( I; X)} $ is the one for small $\delta$'s. More precisely, for $\delta_1 \le \delta_2$
 \begin{equation}\label{eq:jb_strong:nik_bas:n:diff:both:p1}
 [f]_{{r, \delta_1,} N^{\alpha, p} ( I; X)} \le  [f]_{{r, \delta_2,} N^{\alpha, p} ( I; X)} \le  [f]_{{r, \delta_1,} N^{\alpha, p} ( I; X)} + \delta_1^{- \alpha}  \sup_{h \in [ \delta_1, \delta_2) } | \Delta_h^r f |_{L^{p} ( I_{rh}; X)}
 \end{equation}
in view of the formula \eqref{eq:jb_strong:nik_semi:gen}. For the last summand of the r.h.s. of \eqref{eq:jb_strong:nik_bas:n:diff:both:p1} we have
\[
\sup_{h \in [ \delta_1, \delta_2) } | \Delta_h^r f |_{L^{p} ( I_{rh}; X)} \le 2^r | f |_{L^{p} ( I; X)}.
\]
Using this in  \eqref{eq:jb_strong:nik_bas:n:diff:both:p1} and adding to the resulting formula  $ |f|_{L^{ p} ( I; X)}$,  we have for $|f|_{r, \delta, N^{\alpha, p} ( I; X)} $ of \eqref{eq:jb_strong:nik_bas:n:diff}
 \begin{equation*}
 |f|_{{r, \delta_1,} N^{\alpha, p} ( I; X)} \le  |f|_{{r, \delta_2,} N^{\alpha, p} ( I; X)} \le ( 2^r+1)  \delta_1^{- \alpha}   |f|_{{r, \delta_1,} N^{\alpha, p} ( I; X)}, 
 \end{equation*}
 which gives \eqref{eq:jb_strong:nik_bas:n:diff:both:delta}.
 
\underline{Step 3.} Finally, we want to show  \eqref{eq:jb_strong:nik_bas:n:diff:both:r}. Decreasing the number of differences is easy because
 \[
| \Delta_h^{k+l} u |_{L^{p} ( I_{(k+l)h}; X)} =  | \Delta_h^{k} \Delta_h^{l} u |_{L^{p} ( I_{(k+l)h}; X)}   \le 2^k | \Delta_h^{l} u |_{L^{p} ( I_{rh}; X)}.
 \]
 Hence 
 \[
  [f]_{{r, \delta,} N^{\alpha, p} ( I; X)} \le 2^{r-r_0}   [f]_{{r_0, \delta,} N^{\alpha, p} ( I; X)}.
  \]
 In order to add a number of differences, we derive a Marchaud-type inequality. A computation gives
 \[
 \Delta^r_{2h} - 2^r \Delta^r_h = \sum_{k=1}^r {{r}\choose{k}} \sum_{m=0}^{k-1} \Delta^{r+1}_h  T_{mh}.
 \]
 The above formula can be found in the monograph \cite{BenSha88} by Bennett and Sharpley, p. 333. It implies
  \begin{equation}\label{eq:nik:m1}
\left| \Delta^r_{2h} f- 2^r \Delta^r_h f \right|_{L^p (I_{2rh}, X)}\le  \sum_{k=1}^r {{r}\choose{k}} \sum_{m=0}^{k-1} |\Delta^{r+1}_h  T_{mh} f|_{L^p (I_{2rh}, X)} \le  \sum_{k=1}^r {{r}\choose{k}} k |\Delta^{r+1}_h  f|_{L^p (I_{(r+1)h}, X)}.
 \end{equation}
 Observe that the used interval $I_{2hr}$ is the largest admissible, if we want to stay in \eqref{eq:nik:m1} within $I$, \emph{i.e.} the domain of $f$. Formula \eqref{eq:nik:m1} gives
   \begin{equation}\label{eq:nik:m2}
|  \Delta^r_h f -  2^{-r} \Delta^r_{2h} f |_{L^p (I_{2rh}, X)} \le   \frac{r}{2} |\Delta^{r+1}_h  f|_{L^p (I_{(r+1)h}, X)}.
 \end{equation}
 Let us take any natural $j$. One has
   \begin{equation}\label{eq:nik:m3}
|  \Delta^r_h f -  2^{-jr} \Delta^r_{2^jh} f |_{L^p (I_{r2^jh}, X)} \le   \sum_{m=0}^{j-1} 2^{-mr} |   \Delta^r_{2^mh} f - 2^{-r} \Delta^r_{2 2^{m}h} f |_{L^p (I_{r2^mh}, X)}. 
 \end{equation}
 Using for the r.h.s. of \eqref{eq:nik:m3} the inequality  \eqref{eq:nik:m2} with $h: = 2^{m}h$ we get
    \begin{equation}\label{eq:nik:m4}
 \left|  \Delta^r_h f \right|_{L^p (I_{r2^jh}, X)} \le  2^{-jr}\left| \Delta^r_{2^jh} f \right|_{L^p (I_{2^jhr}, X)} +  \frac{r}{2}  \sum_{m=0}^{j-1} 2^{-mr} |\Delta^{r+1}_{2^m h}  f|_{L^p (I_{(r+1)2^m h}, X)}.
 \end{equation}
 We want to increase in the l.h.s. of \eqref{eq:nik:m4} the domain of integration to $I_{hr}$. For the missing part we have
     \begin{equation}\label{eq:nik:m5}
 \left|  \Delta^r_h f \right|_{L^p (I_{rh} \setminus I_{r2^jh}, X)} \le  |  \Delta^r_{-h} f |_{L^p ( I_{-r2^jh}, X)} , \end{equation}
where $I_{-a} := (I_L + a, I_R)$, provided 
      \begin{equation}\label{eq:nik:m55}
     r 2^j |h| \le \frac{|I|}{2}
 \end{equation}
Next, we consider \eqref{eq:nik:m4} with $h$ and with $-h$, add and via \eqref{eq:nik:m5} arrive at
    \begin{multline}\label{eq:nik:m6}
 \left|  \Delta^r_h f \right|_{L^p (I_{hr}, X)} \le  \\
 2^{-jr} \left( \left| \Delta^r_{2^jh} f \right|_{L^p (I_{r2^jh}, X)} +  | \Delta^r_{-2^jh} f |_{L^p (I_{-r2^jh}, X)} \right) + \\ \frac{r}{2}  \sum_{m=0}^{j-1} 2^{-mr} \left(  |\Delta^{r+1}_{2^m h}  f|_{L^p (I_{(r+1)2^m h}, X)}  +  |\Delta^{r+1}_{-2^m h}  f|_{L^p (I_{-(r+1)2^m h}, X)}  \right).
 \end{multline}
 After translations, the "$-h$" terms give a copy of the "$+h$" terms, so we have
\[
 \left|  \Delta^r_h f \right|_{L^p (I_{rh}, X)} \le  \\
 2 2^{-jr}  \left| \Delta^r_{2^jh} f \right|_{L^p (I_{r2^jh}, X)} + r  \sum_{m=0}^{j-1} 2^{-mr}  |\Delta^{r+1}_{2^m h}  f|_{L^p (I_{(r+1)2^m h}, X)}.
\]
Consequently, for any $h>0$
\[
 h^{-\alpha} \left|  \Delta^r_h f \right|_{L^p (I_{rh}, X)} \le  \\
 2 2^{j(\alpha -r)}  (2^j h)^{-\alpha}  \left| \Delta^r_{2^jh} f \right|_{L^p (I_{r2^jh}, X)} + r  \sum_{m=0}^{j-1} 2^{m(\alpha -r)}  (2^m h)^{-\alpha} |\Delta^{r+1}_{2^m h}  f|_{L^p (I_{(r+1)2^m h}, X)}.
\]
 which gives
     \begin{multline}\label{eq:nik:m7}
\sup_{h \le \delta} h^{-\alpha} \left|  \Delta^r_h f \right|_{L^p (I_{rh}, X)} \le  \\
 2 2^{j(\alpha -r)}  \sup_{h_0 \le 2^j \delta} h_0^{-\alpha} | \Delta^r_{h_0} f |_{L^p (I_{rh_0 }, X)} + r  \sum_{m=0}^{j-1} 2^{m(\alpha -r)}  \sup_{h_0 \le 2^m \delta} h_0^{-\alpha} |\Delta^{r+1}_{h_0}  f|_{L^p (I_{(r+1) h_0}, X)}.
 \end{multline}
 We split and estimate the first term of the r.h.s. of \eqref{eq:nik:m7} as follows 
    \begin{multline}\label{eq:nik:m8}
  2 2^{j(\alpha -r)}  \sup_{h_0 \le \delta} h_0^{-\alpha}  | \Delta^r_{h_0} f |_{L^p (I_{rh_0 }, X)}  +   2 2^{j(\alpha -r)} \! \sup_{\delta \le h_0 \le 2^j \delta} h_0^{-\alpha}   | \Delta^r_{h_0 } f |_{L^p (I_{rh_0 }, X)} \le \\
    \frac{1}{2}  \sup_{h_0 \le \delta} h_0^{-\alpha}  | \Delta^r_{h_0} f |_{L^p (I_{r h_0 }, X)}  +      \frac{1}{2}  \delta^{-\alpha}  2^r  | f |_{L^p (I, X)},
    \end{multline}
    where to obtain the factor $\frac{1}{2}$ we need $j \ge \frac{1}{2 (r - \alpha)}$; we choose 
         \begin{equation}\label{eq:nik:m75}
    j =  \left[ \frac{1}{2 (r - \alpha)} \right] + 1.
    \end{equation}
    In \eqref{eq:nik:m75} we can see clearly why one needs for the definition of Nikolskii spaces the difference step $r$ to be sharply larger than the differentiability parameter $\alpha$. The estimate  \eqref{eq:nik:m8} gives via  \eqref{eq:nik:m7}
       \begin{multline}\label{eq:nik:m9}
    \frac{1}{2} \sup_{h \le \delta} h^{-\alpha} \left|  \Delta^r_h f \right|_{L^p (I_{rh}, X)} \le
\frac{1}{2}  \delta^{-\alpha}  2^r  \left| f \right|_{L^p (I, X)}+ rj  \sup_{h_0 \le 2^m \delta} h_0^{-\alpha} |\Delta^{r+1}_{h_0}  f|_{L^p (I_{(r+1) h_0}, X)} \le \\ 
\frac{1}{2}  \delta^{-\alpha}  2^r  \left| f \right|_{L^p (I, X)}+  rj \Big(   \sup_{h_0 \le \delta} h_0^{-\alpha} |\Delta^{r+1}_{h_0}  f|_{L^p (I_{(r+1) h_0}, X)} + \delta^{-\alpha} 2^{r+1} | f|_{L^p (I, X)} \Big).
 \end{multline}
 This estimate yields
        \begin{equation}\label{eq:nik:m77}
  [f]_{{r, \delta,} N^{\alpha, p} ( I; X)} \le  \frac{2r^2}{r-\alpha}\left(     [f]_{{r+1, \delta,} N^{\alpha, p} ( I; X)} +  \delta^{-\alpha} |f|_{L^{ p} ( I; X)} \right),
    \end{equation}
  where we use \eqref{eq:nik:m75} and $r \ge 1$ to write $rj  \le  \sfrac{2r^2}{r-\alpha}$.  Hence we have the missing right inequality in  \eqref{eq:jb_strong:nik_bas:n:diff:both:r} in the case $r -r_0 =1$. The general case, \emph{i.e.}
  \[
   |f|_{ r_0 ,\delta, N^{\alpha, p} ( I; X)}   \le \Big(  \frac{4 r^2}{r_0 - \alpha} \frac{1}{\delta^\alpha} \Big)^{r- r_0} |f|_{r,\delta, N^{\alpha, p} ( I; X)},
  \]
  with any $r \ge r_0$
  follows from an iteration of  \eqref{eq:nik:m77}. We need there \eqref{eq:nik:m55} to hold at every iteration step. We choose generously at every step of iteration \eqref{eq:nik:m75} with $r:= r_0$. To keep \eqref{eq:nik:m55}, it is hence sufficient to require   
  \[
      \delta \le \frac{|I|}{r 2^{j+1}} = \frac{|I|}{r 2^{ \left[ \frac{1}{2 (r_0 - \alpha)} \right] + 2}} ,
  \]
 which is the assumed \eqref{eq:jb_strong:smallness}.
\end{proof}

\begin{proof}[Proof of Proposition \ref{lem:jb_strong:nikolski_red}] 
The fact that $u$ has a weak time derivative $u,_t$ is equivalent to the fact that $u$ has a~representative $u \in AC(I; B)$ for which holds
\begin{equation}\label{lem:jb_weak:wtd:form}
u(t) = \int_{t_0}^t \!u,_t \!(s) \, ds+ u(t_0)
\end{equation}
for $t_0, t \in I$, where the equality is in $B$.

First consider the case $\beta =1$. The assumed finiteness of the r.h.s. of \eqref{eq:redNik} and the definition of a Nikolskii-Bochner allow us to use the representation formula \eqref{lem:jb_weak:wtd:form}. This and the Tonelli Theorem give
\[
\int_{I_{(r+1)h}} \! | \Delta_h \left( \Delta^{r}_h f (\tau)\right) \!|_X^p d \tau = \int_{I_{(r+1)h}} \! h^p \Big|   \dashint_\tau^{\tau +h } \! \Delta_h^{r} f,_s (s) ds  \Big|_X^p d \tau \le \int_{I_{rh}} h^p \left| \Delta_h^{r} f,_s (s)  \right|_X^p \left(  \dashint_{s-h}^{s} \! d \tau \right)  ds.
\]
Hence 
\[
 h^{- (1+\alpha)} | \Delta_h^{r+1} f |_{L^{p} ( I_{(r+1)h}; X)} \le h^{- \alpha} | \Delta_h^r f,_t \!|_{L^{p} ( I_{rh}; X)}.
\]
which via the definition \eqref{eq:jb_strong:nik_semi:gen} of the Nikolskii-Bochner seminorm gives
\[
[f]_{{r+1, \delta,} N^{1+\alpha, p} ( I; X)} \le [f,_t\!]_{{r, \delta,} N^{\alpha, p} ( I; X)}.
\]
Iterating this formula we obtain \eqref{eq:redNik}.
\end{proof}

\begin{proof}[Proof of Proposition \ref{lem:jb_strong:nikolski_ared}]
In contrary to the previous straightforward proofs, here we blend results from \cite{Joh72} by Johnen and \cite{Gol62} by Golovkin. The former defines Besov-Bochner spaces via a modulus of continuity. Therefore, to use its results, we need a connection of our differences-based definition of a Nikolskii-Bochner space with that using a modulus of continuity. For a function $f \in L^{ p} ( I; X) $ and a positive number $h$, the $r$-th order modulus of continuity $\omega_r (f, h)$ reads
\[
\omega_r (f, h) := \sup_{t \in (0, h]} | \Delta^r_t f|_{ L^{ p} ( I_{rt}; X) },
\]
in accordance\footnote{In the $I = \er$, $X= \er$ case it is consistent with Definition 4.2 in \cite{BenSha88}, page 332 or \cite{AdaFou03}, Section 7.46. In the latter the differences are denoted with $\omega$, whereas the modulus of continuity with $\omega^*$.} with formula (1.3) of \cite{Joh72}.
To get the definition of a Besov-Bochner space based on the modulus of continuity, we simply use $ \omega_r (f, h)$ in place of  $ | \Delta_h^r f |_{L^{p} ( I_{rh}; X)} $ in \eqref{eq:jb_strong:nik_semi:gen}. Hence
\begin{equation}\label{eq:jb_strong:nik_semi:gen:nik}
[f]^\omega_{{r, \delta,} N^{\alpha, p} ( I; X)} := \sup_{h \in ( 0, \delta) } h^{- \alpha} \omega_r (f, h) 
\end{equation}
as well as
\[
|f|^\omega_{{r, \delta,} N^{\alpha, p} ( I; X)} := [f]^\omega_{{r, \delta,} N^{\alpha, p} ( I; X)}  +  |f|_{L^{ p} ( I; X)}.
\]
As usual, we will constraint ourselves to $\delta \le 1$, without loss of generality.

\underline{Step 1.} First we observe that our Definition \ref{def:nik} of a Nikolskii-Bochner space is equivalent with the one employing the modulus of continuity. The modulus of continuity majorizes the difference by definition, so we need only the converse.  In the case of Nikolskii spaces, it is also immediate\footnote{In the general case of Besov spaces, one needs to resort to a Steklov-type argument, compare Appendix of \cite{Gol62}.}. We fix $h$ from $(0, \delta)$ and write\vspace{-0pt}
\begin{multline*}
h^{-\alpha} \omega_r (f, h) =  h^{-\alpha}   \sup_{t \in (0, h]}  | \Delta^r_t f|_{ L^{ p} ( I_{rt}; X) }  \le       \sup_{t \in (0, h]} t^{-\alpha}  | \Delta^r_t f|_{ L^{ p} ( I_{rt}; X) } \le   \sup_{t \in (0, \delta)} t^{-\alpha}  | \Delta^r_t f|_{ L^{ p} ( I_{rt}; X) } =\\ [f]_{{r, \delta,} N^{\alpha, p} ( I; X)},
\end{multline*}
hence
\begin{equation}\label{eq:nik:mod1}
 [f]^\omega_{{r, \delta,} N^{\alpha, p} ( I; X)} =  [f]_{{r, \delta,} N^{\alpha, p} ( I; X)}.
\end{equation}

\underline{Step 2. (Extension)} Now we can use the results of \cite{Joh72}. Proposition 5.2.(ii), p. 298 there allows us to extend $f \in N^{\alpha, p} ( I; X) $ to $\bar f \in N^{\alpha, p} ( \er; X) $; its formula (5.14) reads for $h \le 1$
\[
\omega_r (\bar f, h) \le C \left( \omega_r ( f, h) + h^r |f|_{L^{ p} ( I; X) } \right).
\]
Consequently we have
\begin{equation}\label{eq:nik:johnen:1}
\left[\bar f \right]^\omega_{{r, \delta,} N^{\alpha, p} ( \er; X)} \le C | f|^\omega_{{r, \delta,} N^{\alpha, p} ( I; X)}
\end{equation}
and by continuity of this extension from ${L^{ p} ( I; X) }$ to ${L^{ p} ( \er; X) }$ also
\begin{equation}\label{eq:nik:johnen:2}
|\bar f|^\omega_{{r, \delta,} N^{\alpha, p} ( \er; X)} \le C | f|^\omega_{{r, \delta,} N^{\alpha, p} ( I; X)}.
\end{equation}

\underline{Step 3. (Interchange of differences and derivatives.)} Now, we have $\bar f$ defined on the  full real-line, as needed in  \cite{Gol62}. To use results of \cite{Gol62} we go back from the definition by the modulus of continuity to the definition by differences, employing equivalence from the step 1. Next, we use Theorem 1 of \cite{Gol62}, p. 372. It concerns a family of Besov-type spaces. The particular space is realized by the choice of the functional $I$. The case of a Nikolskii space is given by $I:=I_\infty$ from Example 1, \cite{Gol62}, p. 267. We obtain via inequality (11) there the following estimate
\begin{equation}\label{eq:nik:gol:1}
\left[\bar f,_t^{(\beta)} \!\right]_{{r - \beta, \infty,} N^{\alpha-\beta, p} ( \er; X)}  \le C\left[\bar f\, \right]_{{r, \infty,} N^{\alpha, p} ( \er; X)} 
\end{equation}
with a natural $\alpha \ge \beta$ and $r \in \er$, $r> \alpha$. We combine \eqref{eq:nik:mod1}, an analogue of right inequality of   \eqref{eq:jb_strong:nik_bas:n:diff:both:delta} for $ | \bar f|^\omega_{{r, \delta,} N^{\alpha, p} ( \er; X)}  $ and \eqref{eq:nik:johnen:2} to get
\begin{equation}\label{eq:nik:gol:2}
 |\bar f|_{{r, \infty,} N^{\alpha, p} ( \er; X)} =  | \bar f|^\omega_{{r, \infty,} N^{\alpha, p} ( \er; X)}  \le   \frac{3^{r }}{\delta^\alpha}  | \bar f|^\omega_{{r, \delta,} N^{\alpha, p} ( \er; X)} \le C \frac{3^{r }}{\delta^\alpha}  | f|^\omega_{{r, \delta,} N^{\alpha, p} ( I; X)} = C \frac{3^{r }}{\delta^\alpha}  | f|_{{r, \delta,} N^{\alpha, p} ( I; X)}
\end{equation}
 Next, we use \eqref{eq:nik:gol:2} to estimate the r.h.s. of  \eqref{eq:nik:gol:1} and get
 \begin{equation}\label{eq:nik:gol:3}
\left[\bar f,_t^{(\beta)}\right]_{{r - \beta, \infty,} N^{\alpha-\beta, p} ( \er; X)}  \le  C \frac{3^{r }}{\delta^\alpha}  | f|_{{r, \delta,} N^{\alpha, p} ( I; X)},
\end{equation}
which after restriction of $\er$ to $I$, where $\bar f = f$, implies 
\begin{equation}\label{eq:accNik:semi}
\left[ f,_t^{(\beta)}\right]_{{r - \beta,} \delta, N^{\alpha-\beta, p} ( I; X)}  \le  C \frac{3^{r }}{\delta^\alpha}  | f|_{{r, \delta,} N^{\alpha, p} ( I; X)}.
\end{equation}
\underline{Step 4. (Control of the full norm.)} Observe that the l.h.s. of \eqref{eq:accNik:semi} contains only a seminorm. Now we gain control over the full norm via formula (7.2) of Proposition 7.1 from \cite{Joh72}. This Proposition is in fact stated for the case of $I$ being the real line or a half-line, hence we need to execute the formula (7.2) on the extended function $\bar f$. It reads
\[
|\bar f,_t^{(\beta)}|_{L^{p} ( \er; X)} \le C \int_\er h^{-(\beta+1)} \omega_r (\bar f, h) dh \le C \int_0^1 h^{\alpha-(\beta+1)} h^{-\alpha} \omega_r (\bar f, h) dh  + C \int_1^\infty \omega_r (\bar f, h) dh. 
\]
Next, we use the defining formula \eqref{eq:jb_strong:nik_semi:gen:nik} for the first integral on the r.h.s. above. The remaining term there is integrable for $\alpha > \beta$. For the second integral  of the r.h.s. above we use the triangle inequality. We get
\vspace{-2mm}
\begin{multline}\label{eq:nik:gol:5}
|\bar f_t^{(\beta)}|_{L^{p} ( \er; X)} \le  \frac{C}{\alpha - \beta} [\bar f]^\omega_{{r, 1,} N^{\alpha, p} ( \er; X)}   + C 2^r |\bar f|_{L^{p} ( \er; X)} \le  C (\alpha - \beta) 2^r | \bar f|^\omega_{{r, \infty,} N^{\alpha, p} ( \er; X)} \le \\ \frac{C (\alpha - \beta)  6^{r }}{\delta^\alpha}  | f|_{{r, \delta,} N^{\alpha, p} ( I; X)}.
\end{multline}
The last inequality follows from \eqref{eq:nik:gol:2}. Putting together \eqref{eq:accNik:semi} and \eqref{eq:nik:gol:5} we arrive at \eqref{eq:accNik}.
\end{proof}

\noindent
\begin{proof}[Proof of Proposition \ref{prop:abstractNikInter}]
We write
\[
\frac{1}{h^{\alpha_b}} \left( \int_{I_{rh}} \left| \Delta^r_h f (t) \right|^{p_b}_{Z} dt \right)^\frac{1}{p_b} \le C_{Z,X,Y} \frac{1}{h^{ (1-b) \alpha_1 }} \left( \int_{I_{rh}} \left| \Delta^r_h f (t) \right|^{p_1}_{X} dt \right)^\frac{1-b}{p_1}  \frac{1}{h^{ b \alpha_2} }\left( \int_{I_{rh}} \left| \Delta^r_h f (t) \right|^{p_2}_{Y} dt \right)^\frac{b}{p_2},
\]
where we have used the H\"older inequality.
\end{proof}

\begin{proof}[Proof of Lemma \ref{lem:jb_strong:nikolski_inter}]
 We get our result in two steps: a duality and a standard interpolation. 
More precisely, the idea of the proof is as follows. First we obtain an interpolation of the type
\begin{equation}\label{eq:jb_strong:nikolski_inter:p1}
 |f|_{\delta, N^{\alpha (\tilde \theta), p({\tilde \theta})} ( I; W^{k (\tilde \theta),  q (\tilde \theta)} (Q)) }  \le C  |f|^{1- \tilde \theta}_{\delta, N^{ \alpha_1, p_1} ( I; W_0^{1, q_1} (Q)) }  |f|^{\tilde \theta}_{\delta, N^{ \alpha_2, p_2} ( I; W^{-1, q_2} (Q)) }
\end{equation}
with $k (\tilde \theta)$ positive. Next we check what is the optimal $q^* (\tilde \theta)$ of embedding $W^{k (\tilde \theta),  q (\tilde \theta)} (Q) \csubset L^{q^* (\tilde \theta)}$. Unfortunately, on domains,  the interpolation theory of Sobolev spaces with a real differentiability parameter seems incomplete. There is a problem with the local description for spaces with negative differentiability, see the Triebel's classic \cite{Tri83}, p. 209. This has not changed much in modern times, compare Section 4.1.4. of \cite{Tri06}. Therefore, we carry out  our interpolation into two steps. First we use simply a duality formula which deals with the negative Sobolev space $W^{-1, q_2}$. Next we use the~complex interpolation.

\underline{Step 1. (Use of duality.)} Recall that we have restricted ourselves to $\alpha_1 \le \alpha_2 < 2$, so it is sufficient to use  the difference step $r=2$ in Nikolskii seminorms. In view of the assumed $q_1 \ge q'_2$, $p_1 \ge p'_2$ and the~H\"older inequality  we have
\begin{equation}\label{eq:jb_strong:nikolski_inter:1}
 [f]_{2, \delta, N^{ \alpha_1, p'_2} ( I; W_0^{1, q'_2} (Q)) } \le C  [f]_{2, \delta, N^{\alpha_1, p'_2} ( I; W_0^{1, q_1} (Q)) } \le C  [f]_{2, \delta, N^{ \alpha_1, p_1} ( I; W_0^{1, q_1} (Q)) }.
\end{equation}
For the  domain-independence of the constant above, recall Remark \ref{rem:uberunter} (this applies to the further embedding constants as well). 
The duality between  $W_0^{1, q'_2} (Q)$ and $W^{-1, q_2} (Q)$ gives for  any $ g \in~W_0^{1, q'_2} (Q)$ via the identification
\[\langle g , g  \rangle_{(W_0^{1, q'_2} (Q))^*, \,W_0^{1, q'_2} (Q)} =  \langle g , g  \rangle_{L^2 (Q), \, L^2 (Q)} \]
the estimate
\[
|g|_{L^2 (Q)} \le |g|^\frac{1}{2}_{W_0^{1, q'_2} (Q)}  |g|^\frac{1}{2}_{(W_0^{1, q'_2} (Q))^*}.
\]
It allows us to write  via Proposition \ref{prop:abstractNikInter} 
\begin{equation}\label{eq:jb_strong:nikolski_inter:12}
 [f]_{2, \delta, N^{\frac{\alpha_1 + \alpha_2}{2}, 2} ( I; L^2 (Q)) } \le    [f]^\frac{1}{2}_{2, \delta, N^{\alpha_1, p'_2} ( I; W_0^{1, q'_2} (Q)) }  [f]^\frac{1}{2}_{2, \delta, N^{ \alpha_2, p_2} ( I; W^{-1, q_2} (Q)) }.
 \end{equation}
Inequality \eqref{eq:jb_strong:nikolski_inter:12} with \eqref{eq:jb_strong:nikolski_inter:1} gives
\begin{equation}\label{eq:jb_strong:nikolski_inter:2}
 [f]_{2, \delta, N^{\frac{\alpha_1 + \alpha_2}{2}, 2} ( I; L^2 (Q)) } \le   C  [f]^\frac{1}{2}_{2, \delta, N^{ \alpha_1, p_1} ( I; W_0^{1, q_1} (Q)) }   [f]^\frac{1}{2}_{2, \delta, N^{ \alpha_2, p_2} ( I; W^{-1, q_2} (Q)) }.
 \end{equation}
 As intended, now we don't have any dependence on negative differentiability in the l.h.s. of \eqref{eq:jb_strong:nikolski_inter:2}.
 
 \underline{Step 2. (The complex interpolation.)} 
Now we interpolate between the  l.h.s. of \eqref{eq:jb_strong:nikolski_inter:2} and $N^{\alpha_1, p_1} (I; W_0^{1, q_1} (Q)) $. We want to use again Proposition \ref{prop:abstractNikInter}, this time using the complex interpolation for the spaces on $Q$. Hence first we need to identify the space
 \[Y = [L^2 (Q), W^{1, q_1} (Q)]_\theta\]
 where $ [\cdot, \cdot]_\theta$ denotes the complex interpolation functor. Here we see the benefit of splitting our proof into two steps. Namely, for a natural, nonnegative $k$ and $l \in (1, \infty) $ we have $W^{k,l} = F^{k,l}_2$ (see formula (1) of Section 3.4.2 of \cite{Tri83}), where $ F^{k,l}_r$ stands for the Triebel-Lizorkin space (for its definition on a smooth domain $Q$, see Section 3.4.2 of \cite{Tri83}). Triebel-Lizorkin spaces are well-behaved under the~complex interpolation, see  formula (8) of the theorem in Section 3.3.6 of \cite{Tri83}. Therefore 
 \begin{equation}
 Y = [F^{0, 2}_2 (Q), F^{1, q_1}_2 (Q)]_\theta =F^{1-\theta, \tilde q}_2 (Q) \quad \text{ with  } \quad \frac{1}{\tilde q} = \frac{\theta}{2} + \frac{1-\theta}{q_1}.
 \end{equation} 
 In view of the theorem in Section 3.3.1 of \cite{Tri83}, for  positive $p_0, p_1$ and real $s_0 > s_1$ we have the~continuous embedding $ F^{s_0, p_0}_2 (Q) \csubset F^{s_1, p_1}_2 (Q) $, provided $s_0 - \frac{d}{p_0} \ge s_1 - \frac{d}{p_1} $. This gives in our case
  \begin{equation}\label{eq:jb_strong:nikolski_inter:2''}
  Y \csubset  F^{0,\tilde q_0}_2 (Q) = L^{\tilde q_0} (Q)
 \end{equation}
for any finite $\tilde q_0$, provided
  \[
   1- \theta - \frac{d}{\tilde q} \ge 0. \] 
   Otherwise  the largest allowed  $\tilde q_0 $ (finite and $>1$) is
    \[
 \tilde q_0 = \frac{\tilde q d}{d - \tilde q (1- \theta)} = \frac{2 d q_1}{ q_1 \theta d +  (1- \theta) (2 d - 2 q_1)}.
 \]

 Embedding \eqref{eq:jb_strong:nikolski_inter:2''} in these cases can be written together as
 \begin{equation}\label{eq:jb_strong:nikolski_inter:23}
  Y \csubset L^{q_0} (Q)
 \end{equation}
with $q_0= \tilde q_0$ if $\tilde q_0$ is positive and any finite $q_0$ otherwise. This is precisely our condition   \eqref{eq:jb_strong:nikolski_inter:parac}. Quantitatively, \eqref{eq:jb_strong:nikolski_inter:23} reads
 \[
 |g|_{ L^{q_0} (Q) }  \le  C\, |g|^\theta_{ L^{2} (Q) } |g|^{1 - \theta}_{ W^{1, q_1} (Q)}.
\]
We use it in Proposition \ref{prop:abstractNikInter} to get
\begin{equation}\label{eq:jb_strong:nikolski_inter:02}
 [f]_{2, \delta, N^{\alpha, p_0} ( I; L^{q_0} (Q)) }  \le  C \,[f]^\theta_{2, \delta, N^{\frac{\alpha_1+ \alpha_2}{2}, 2} ( I; L^2 (Q)) }   [f]^{1 - \theta}_{2, \delta, N^{\alpha_1, p_1} (I; W^{1, q_1} (Q))},
\end{equation}
where
\begin{equation}\label{eq:jb_strong:nikolski_inter:2'}
\alpha = \theta \frac{\alpha_1 + \alpha_2}{2} + (1- \theta) \alpha_1, \; \; \frac{1}{p_0} = \frac{\theta}{2} + \frac{1-\theta}{p_1}. 
\end{equation}
Estimating the first term of the r.h.s. of \eqref{eq:jb_strong:nikolski_inter:02} with \eqref{eq:jb_strong:nikolski_inter:2}, we arrive at
\begin{equation}\label{eq:jb_strong:nikolski_inter:3}
 [f]_{2, \delta, N^{\alpha, p_0} ( I; L^{q_0} (Q)) }  \le C \,  [f]^\frac{\theta}{2}_{2, \delta, N^{ \alpha_2, p_2} ( I; W^{-1, q_2} (Q)) }  [f]^{1-\frac{\theta}{2}}_{2, \delta, N^{ \alpha_1, p_1} ( I; W^{1, q_1} (Q)) }
 \end{equation}
 Analogously to \eqref{eq:jb_strong:nikolski_inter:3} we obtain the estimate for the low-order part of the Nikolskii-Bochner norm, \emph{i.e.}
 \begin{equation}\label{eq:jb_strong:nikolski_inter:3:low}
 |f|_{L^{ p_0} ( I; L^{q_0} (Q)) }  \le C \,  |f|^\frac{\theta}{2}_{L^{  p_2} ( I; W^{-1, q_2} (Q)) }  |f|^{1-\frac{\theta}{2}}_{L^{ p_1} ( I; W^{1, q_1} (Q)) }.
 \end{equation}
Formulas \eqref{eq:jb_strong:nikolski_inter:3} and \eqref{eq:jb_strong:nikolski_inter:3:low} give together
 \begin{equation}\label{eq:jb_strong:nikolski_inter:3a}
 |f|_{2, \delta, N^{\alpha, p_0} ( I; L^{q_0} (Q)) }  \le C \,  |f|^\frac{\theta}{2}_{2, \delta, N^{ \alpha_2, p_2} ( I; W^{-1, q_2} (Q)) }  |f|^{1-\frac{\theta}{2}}_{2, \delta, N^{ \alpha_1, p_1} ( I; W^{1, q_1} (Q)) }.
 \end{equation}
Where applicable, we use \eqref{eq:jb_strong:nik_bas:n:diff:both:r} of Lemma \ref{lem:jb_strong:nik_bas} to decrease in \eqref{eq:jb_strong:nikolski_inter:3a} the differentiability step from $2$ to the >>natural one<<, which is either $1$ or $2$. It forces us to assume \eqref{eq:jb_strong:smallness}, which is reflected here by \eqref{eq:jb_strong:smallness'}, as we have restricted ourselves to $r \le 2$. Hence we get from \eqref{eq:jb_strong:nikolski_inter:3a}
 \begin{equation}\label{eq:jb_strong:nikolski_inter:3b}
 |f|_{\delta, N^{\alpha, p_0} ( I; L^{q_0} (Q)) }  \le C \Big( \frac{1}{\left( [\alpha] - \alpha +1 \right) \delta^\alpha} \Big)^{1 - [\alpha]} |f|^\frac{\theta}{2}_{\delta, N^{ \alpha_2, p_2} ( I; W^{-1, q_2} (Q)) }  |f|^{1-\frac{\theta}{2}}_{\delta, N^{ \alpha_1, p_1} ( I; W^{1, q_1} (Q)) }.
 \end{equation}
The fact that $f$ has zero space-trace gives via  \eqref{eq:jb_strong:nikolski_inter:3b} the desired estimate \eqref{eq:jb_strong:nikolski_inter}. Recall that  $p_0, \alpha$ come from \eqref{eq:jb_strong:nikolski_inter:2'}.
 \end{proof}

\begin{proof}[Proof of Lemma \ref{lem:nik:emb}] Ad (i). The embedding  $N^{1+\gamma, p} ( I; X) \csubset W^{1, p} ( I; X)$ with its estimate  \eqref{eq:nik:embW} follows from the accession Proposition \ref{lem:jb_strong:nikolski_ared}, because for any natural $\beta \le [\alpha]$  it gives
 \begin{equation}\label{lem:nik:emb:05}
\left| f,_t^{(\beta)}\! \right|_{ L^{p} ( I; X)} \le  \left| f,_t^{(\beta)} \! \right|_{[\alpha]+1 - \beta, \delta, N^{[\alpha] + \gamma - \beta, p} ( I; X)}  \le   \frac{C 6^{[\alpha]+1 }}{([\alpha] + \gamma - \beta) \delta^{[\alpha] + \gamma}} | f|_{{ \delta,} N^{[\alpha] + \gamma, p} ( I; X)}.
 \end{equation}
We sum \eqref{lem:nik:emb:05} over $\beta$ varying from $0$ to $[\alpha]$ to get \eqref{eq:nik:embW}.

Ad (ii). Let us first observe that in the case $ \alpha =\alpha' $ the embedding  follows from the H\"older inequality. Moreover, it suffices to show \eqref{lem:nik:emb:th:low} in the case $p \le q$, because having this case we show the complementary one $p > q$ with the H\"older inequality. For the case $p \le q$ and $ \alpha >\alpha' $, we proceed in a few steps.

\underline{Step 1. (The low differentiability case)}
For the case
 \begin{equation}\label{lem:nik:emb:215}
 \alpha - \frac{1}{p} \ge \alpha' - \frac{1}{q}, \qquad 1 > \alpha >  \alpha' >0
  \end{equation}
we quote \cite{Sim90}. More precisely, Corollary 22 there gives
 \begin{equation}\label{lem:nik:emb:2}
 [f]_{1,|I|, N^{\alpha', q} ( I; X) } \le \frac{18}{\alpha'} |I|^{\alpha - \frac{1}{p} -( \alpha' - \frac{1}{q})}  |f|_{1, |I|, N^{\alpha, p} ( I; X) }.
 \end{equation}
This and \eqref{eq:jb_strong:nik_bas:n:diff:both:delta} give
 \begin{equation}\label{lem:nik:emb:21}
 |f|_{1, \delta, N^{\alpha', q} ( I; X) } \le \frac{C}{\alpha' \delta^\alpha} |f|_{1, \delta, N^{\alpha, p} ( I; X) } \le \frac{C}{\alpha' \delta} |f|_{1, \delta, N^{\alpha, p} ( I; X) }.
 \end{equation}
Next, we extend \eqref{lem:nik:emb:21} over $\alpha' =0 $ and  $\alpha =1 $ and modify it, so that its constant does not blow up as $\alpha' \to 0$ (with all the other parameters fixed). This extension will be performed at the cost of assuming the sharp inequality in the first condition of \eqref{lem:nik:emb:215}. Namely, we substitute \eqref{lem:nik:emb:215} with
 \begin{equation}\label{lem:nik:emb:216}
\alpha - \frac{1}{p} > \alpha' - \frac{1}{q}, \qquad 1 \ge \alpha > \alpha' \ge 0.
  \end{equation}
\underline{Substep 1.1} \, Let us consider  \eqref{lem:nik:emb:216} in the case $ p=q $. First, for  $\alpha < 1$, in both norms of \eqref{lem:nik:emb:th:low} we have the difference step $r=1$, compare Definition \ref{def:nik}. Therefore the definition of a Nikolskii-Bochner seminorm, where we use the restriction $\delta \le 1$, gives
 \begin{equation}\label{lem:nik:emb:216pp}
|f|_{{1, \delta,} N^{\alpha', p} ( I; X)} \le |f|_{{1, \delta,} N^{\alpha, p} ( I; X)}.
  \end{equation}
Next, for  $\alpha = 1$ and $\alpha' < 1$, we need additionally to change the difference step by
\eqref{eq:jb_strong:nik_bas:n:diff:both:r}. It gives
 \begin{equation}\label{lem:nik:emb:2165}
   |f|_{1 ,\delta, N^{\alpha', p} ( I; X)}   \le \Big(  \frac{16}{(1 - \alpha')  \delta^{\alpha'} }\Big) |f|_{2, \delta, N^{\alpha', p} ( I; X)}   \le \Big(  \frac{16}{(1 - \alpha')  \delta }\Big) |f|_{2, \delta, N^{\alpha', p} ( I; X)},  
  \end{equation}
where we need \eqref{eq:jb_strong:smallness}, \emph{i.e.}
\[
       \delta \le \frac{|I|}{ 2^{ \left[ \frac{1}{2 (1 - \alpha')} \right] + 3}}.
\]
We estimate the r.h.s. of \eqref{lem:nik:emb:2165} via the definition of a Nikolskii-Bochner seminorm, analogously to \eqref{lem:nik:emb:216pp} and get
\[
   |f|_{1, \delta, N^{\alpha', p} ( I; X)}   \le  \Big(  \frac{16}{(1 - \alpha')  \delta }\Big) |f|_{2, \delta, N^{1, p} ( I; X)}.
 \]
\underline{Substep 1.2} \, We are left with extending \eqref{lem:nik:emb:21} over $1 \ge \alpha > \alpha'  \ge 0 $ in the case $p < q$. Let us choose $\tilde \alpha, \tilde \alpha'$ such that $\tilde \alpha - \sfrac{1}{p}, \, \tilde \alpha' - \sfrac{1}{q}$ are equal to the midpoint of the latter condition in \eqref{lem:nik:emb:216}, \emph{i.e.}
 \begin{equation}\label{lem:nik:emb:2169}
\tilde \alpha - \frac{1}{p} = \tilde \alpha' - \frac{1}{q}  = \frac{1}{2} \Big( \Big(\alpha - \frac{1}{p} \Big) + \Big(\alpha' - \frac{1}{q} \Big)  \Big).
  \end{equation}
This choice of $\tilde \alpha, \tilde \alpha'$, the assumption $(\alpha - \sfrac{1}{p} ) - (\alpha' - \sfrac{1}{q} )  > 0$, compare
\eqref{lem:nik:emb:216}, and $p < q$ give 
 \begin{equation}\label{lem:nik:emb:217}
\tilde \alpha - \frac{1}{p} = \tilde   \alpha' - \frac{1}{q}, \qquad 1 \ge \alpha > \tilde \alpha > \tilde \alpha' > \alpha' \ge 0.
  \end{equation}
 The $(\tilde \alpha, \tilde \alpha')$-part of \eqref{lem:nik:emb:217} complies with \eqref{lem:nik:emb:215}. Hence \eqref{lem:nik:emb:21}  gives
 \begin{equation}\label{lem:nik:emb:2171p}
 |f|_{1, \delta, N^{\tilde \alpha', q} ( I; X) } \le \frac{C}{\tilde \alpha' \delta} |f|_{1, \delta, N^{\tilde \alpha, p} ( I; X) } 
  \end{equation}
  The l.h.s. of \eqref{lem:nik:emb:2171p} controls the wanted $N^{ \alpha', q}$ via the definition of a Nikolskii-Bochner space, because $\alpha' < \tilde \alpha'$. Similarly, the r.h.s. of \eqref{lem:nik:emb:2171p} is controlled with $N^{ \alpha, p}$ norm; here for $\alpha=1$ we need also to change the order of differences from $1$ to $2$ via   \eqref{eq:jb_strong:nik_bas:n:diff:both:r}. More precisely, we get
   \begin{equation}\label{lem:nik:emb:2171}
  |f|_{1, \delta, N^{\alpha', q} ( I; X) } \le \frac{C}{\tilde \alpha' \delta} |f|_{1, \delta, N^{ \alpha, p} ( I; X) } 
    \end{equation}
  in the case $\alpha <1$ and 
     \begin{equation}\label{lem:nik:emb:2172}
 |f|_{1, \delta, N^{\tilde \alpha', q} ( I; X) } \le \frac{C}{ \tilde \alpha' \delta} |f|_{1, \delta, N^{\tilde \alpha, p} ( I; X) } \le\frac{C}{({1 - \tilde \alpha})\tilde \alpha' \delta^{2}} |f|_{2, \delta, N^{\tilde \alpha, p} ( I; X) } \le \frac{C}{({1 - \tilde \alpha})\tilde \alpha' \delta^{2}} |f|_{2, \delta, N^{1, p} ( I; X) } 
  \end{equation}
  in the case $\alpha =1$,  where we need \eqref{eq:jb_strong:smallness}, \emph{i.e.}
\[
       \delta \le \frac{|I|}{ 2^{ \left[ \frac{1}{2 (1 -\tilde \alpha)} \right] + 3}} =  \frac{|I|}{ 2^{ \left[ \left(\left(\alpha - \frac{1}{p} \right) - \left(\alpha' - \frac{1}{q} \right)\right)^{-1} \right] + 3}}.
\]
Substeps 1.1 --- 1.2 yield
    \begin{equation}\label{lem:nik:emb:235}
 |f|_{\delta, N^{\alpha', q} ( I; X) } \le C_\eqref{lem:nik:emb:235} (\alpha, \alpha', p, q, \delta)  |f|_{ \delta, N^{\alpha, p} ( I; X) },
  \end{equation}
provided \eqref{lem:nik:emb:216} holds,  where
    \begin{equation}\label{lem:nik:emb:235a}
    C_\eqref{lem:nik:emb:235}  (\alpha, \alpha', p, q, \delta) := \!
  \begin{cases}
C \quad \text{for} \quad   1> \alpha > \alpha', \;   p = q,\\
    \frac{C}{\left(1 - \alpha' \right)\delta}   \quad \text{for} \quad   1= \alpha > \alpha', \;    p = q,\\
  \frac{C}{\tilde \alpha' \delta}  =   \frac{C}{\left(\alpha - \frac{1}{p} + \alpha'+ \frac{1}{q}  \right)\delta}  \quad \text{for} \quad   1> \alpha > \alpha', \;    p < q,\\
    \frac{C}{\tilde \alpha'  \left(1 - \tilde \alpha \right) \delta^{2}} =  \frac{C}{\left(\alpha - \frac{1}{p} + \alpha'+ \frac{1}{q}  \right) \left(\left(\alpha - \frac{1}{p} \right) - \left(\alpha' - \frac{1}{q} \right) \right) \delta^{2}}    \quad \text{for} \quad   1= \alpha > \alpha', \;    p < q,\\
  \end{cases}
  \end{equation}
 where $\tilde \alpha, \tilde \alpha'$ are given by  \eqref{lem:nik:emb:2169}. In  the case $1= \alpha > \alpha'$ we need also
    \begin{equation}\label{lem:nik:emb:236}
  \delta \le
       \begin{cases}
        \frac{|I|}{ 2^{ \left[ (2 (1 - \alpha' ))^{-1} \right] + 3}} \quad \text{for} \quad   1= \alpha > \alpha', \;   p = q, \\
        \frac{|I|}{ 2^{ \left[ (2 (1 - \tilde \alpha ))^{-1} \right] + 3}} \quad\; \,\text{for}  \quad   1= \alpha > \alpha', \;  p < q.
 \end{cases}
  \end{equation}
We can generously simplify \eqref{lem:nik:emb:235},  \eqref{lem:nik:emb:236}, observing that  for $\beta := \left(\alpha - \sfrac{1}{p} \right) - \left(\alpha' - \sfrac{1}{q} \right)$
    \begin{equation}\label{lem:nik:emb:235'}
    C_\eqref{lem:nik:emb:235}  (\alpha, \alpha', p, q, \delta) \le     C_\eqref{lem:nik:emb:235'}  (\beta, \delta) := \frac{C}{\beta^2 \delta}   
  \end{equation}
  and \eqref{lem:nik:emb:236} follows from
    \begin{equation}\label{lem:nik:emb:236'}
  \delta \le
        \frac{|I|}{ 2^{ \left[ \frac{1}{2 (1 - \tilde \alpha )} \right] + 3}} =   \frac{|I|}{ 2^{ \left[ \frac{1}{\beta } \right] + 3}} \quad \text{for} \quad   1= \alpha > \alpha'.
  \end{equation}
  
\underline{Conclusion of step 1.} After generous estimates on constants the outcome of this step reads
      \begin{equation}\label{lem:nik:emb:235gen}
 |f|_{\delta, N^{\alpha', q} ( I; X) } \le C_\eqref{lem:nik:emb:235'} (\beta, \delta)  |f|_{ \delta, N^{\alpha, p} ( I; X) } 
  \end{equation}
provided \eqref{lem:nik:emb:216} holds. For the case  $1= \alpha> \alpha'$ we need additionally \eqref{lem:nik:emb:236'}.
  
\underline{Step 2. (The general case)}
Now we generalize \eqref{lem:nik:emb:235gen} over any $\alpha > \alpha' \ge 0$ such that $\alpha - \frac{1}{p} > \alpha' - \frac{1}{q}$. 

\underline{Substep 2.1}  \, First, consider the case when $\alpha> \alpha'$ both belong to the interval $[k, k+1]$, with an arbitrary $k \in \N$. We use the result of step 1 and the reduction/accession Propositions \ref{lem:jb_strong:nikolski_red}, \ref{lem:jb_strong:nikolski_ared}. More precisely, we write
 \begin{multline}\label{lem:nik:emb:3}
 [f]_{{\delta,}  N^{\alpha', q} ( I; X) } \le  \left[ f,_t^{(k)} \right]_{{\delta,} N^{\alpha'- k, q} ( I; X)} \le C_\eqref{lem:nik:emb:235}  (\alpha-k, \alpha'-k, p, q, \delta)  |f,_t^{(k)}|_{ \delta, N^{\alpha-k, p} ( I; X) } \le \\
C_\eqref{lem:nik:emb:235}  (\alpha-k, \alpha'-k, p, q, \delta)   \frac{C 6^{k }}{(\alpha - k) \delta^\alpha}  | f|_{{ \delta,} N^{\alpha, p} ( I; X)}.
  \end{multline}
 In order to obtain the first inequality of \eqref{lem:nik:emb:3} we have used the reduction Proposition \ref{lem:jb_strong:nikolski_red}. For the~second one we have employed \eqref{lem:nik:emb:235}. The last one follows from Proposition \ref{lem:jb_strong:nikolski_ared}. As the r.h.s. of \eqref{lem:nik:emb:3} contains already the lower-order part of the full Nikolskii-Bochner norm,   \eqref{lem:nik:emb:3}  gives with $k = [ \alpha']$
 \vspace{-4pt}
\begin{multline}\label{lem:nik:emb:3a}
 |f|_{{\delta,}  N^{\alpha', q} ( I; X) } \le C_\eqref{lem:nik:emb:235}  (\alpha-[ \alpha'], \alpha'-[ \alpha'], p, q, \delta)   \frac{C 6^{[ \alpha'] }}{(\alpha - [ \alpha']) \delta^\alpha}  | f|_{{ \delta,} N^{\alpha, p} ( I; X)} \le\\
   \frac{C}{\beta^2 \delta}     \frac{6^{\alpha }}{(\alpha - \alpha') \delta^\alpha}  | f|_{{ \delta,} N^{\alpha, p} ( I; X)}, 
\end{multline}
where the second inequality follows from \eqref{lem:nik:emb:235'}. In view of step 1 we see that   \eqref{lem:nik:emb:3a} 
 holds, provided the shifted  \eqref{lem:nik:emb:216} is valid, \emph{i.e.}
  \begin{equation}\label{lem:nik:emb:216k}
\alpha - \frac{1}{p} > \alpha' - \frac{1}{q}, \qquad k+1 \ge \alpha > \alpha' \ge k.
  \end{equation}
 In the $k+1= \alpha > \alpha'$ case we also need the shifted condition \eqref{lem:nik:emb:236'}, \emph{i.e.}
    \begin{equation}\label{lem:nik:emb:236s}
  \delta \le  \frac{|I|}{ 2^{ \left[ \frac{1}{\beta} \right] + 3}}.
  \end{equation}
   \underline{Substep 2.2}  In the situation complementary to the one considered in substep 2.1, there exists a~natural $K \ge 1$ that 
\[
[\alpha'] \le \alpha' < [\alpha'] + 1 < \dots [\alpha'] + K < \alpha \le [\alpha'] + K + 1.
\]
Let us realize that
\begin{equation}\label{eq:Knat}
 [\alpha'] + K =
        \begin{cases}
         [\alpha] \quad \quad\;\;\, \text{for} \quad  \alpha > [\alpha],\\
        [\alpha]-1 \quad \text{for} \quad  \alpha = [\alpha].
 \end{cases}
\end{equation}
Now we intend to use the previous substep $K+1$ times. To this end we define
\[
\alpha_0 := \alpha', \qquad  \alpha_i :=  [\alpha'] + i, \quad i = 1, \dots K, \qquad \alpha_{K+1} := \alpha
\]
and, recalling that $\beta$ defined as $\alpha - \sfrac{1}{p} - \left(\alpha' - \sfrac{1}{q} \right) $ is strictly positive by assumption, we define intermediary 
\[
\beta_i := \beta \frac{\alpha_{{i+1}} - \alpha_i}{\alpha_{K+1} - \alpha_0} = (\alpha_{{i+1}} - \alpha_i) \left(1+\Big( \frac{1}{q} - \frac{1}{p} \Big) \frac{1}{\alpha_{K+1} - \alpha_0} \right) , \qquad \text{for} \quad i =0, \dots, K.
\]
Observe that $\sum_{i=0}^K \beta_i = \beta$. Next, we use $\beta_i$'s to define
  \begin{equation}\label{lem:nik:emb:216pk}
\frac{1}{q_0} := \frac{1}{q}, \quad  \frac{1}{q_{i+1}}  :=  \alpha_{i+1}  - \beta_i - \Big( \alpha_{i}  - \frac{1}{q_i} \Big) =  \frac{1}{q_i} +  \Big(\frac{1}{p} - \frac{1}{q} \Big) \frac{ \alpha_{{i+1}} - \alpha_i}{\alpha_{K+1} - \alpha_0}  , \quad \text{for} \quad i =0, \dots, K.
\end{equation}
Observe that $\sfrac{1}{q_{i}} \le \sfrac{1}{q_{i+1}}  $ because $p \le q$. Moreover, the second inequality in \eqref{lem:nik:emb:216pk} yields
\[
 \frac{1}{q_{i+1}}  -  \frac{1}{q_i} =  \Big( \frac{1}{p} - \frac{1}{q} \Big) \frac{ \alpha_{{i+1}} - \alpha_i}{\alpha_{K+1} - \alpha_0},  
\]
which after summing up over $ i =0, \dots, K$ gives 
\[
 \frac{1}{q_{K+1}} = \frac{1}{p}.
 \]
Hence $q_i$'s decrease from $q$ to $p$ as $i$ increases from $0$ to $K+1$. Therefore $q_i$'s are admissible as (intermediary) integrability parameters. Now we use \eqref{lem:nik:emb:3a} iteratively  for the intermediate pairs $N^{\alpha_i, q_i} ( I; X)$, $N^{\alpha_{i+1}, q_{i+1}} ( I; X) $ to get
  \begin{equation}\label{lem:nik:emb:3al}
 |f|_{{\delta,}  N^{\alpha', q} ( I; X) } \le \Big(  \frac{C_{\eqref{lem:nik:emb:3al}}(\alpha, \alpha') }{\beta^2 \delta^{1+\alpha}}\Big)^{\alpha - \alpha' + 2}   | f|_{{ \delta,} N^{\alpha, p} ( I; X)},
  \end{equation}
   where larger than $1$
   \[
   C_{\eqref{lem:nik:emb:3al}}(\alpha, \alpha') =
           \begin{cases}
     C 6^{\alpha}       \frac{(\alpha - \alpha')^2}{ ([ \alpha'] +1 -  \alpha')^3} \quad\quad\quad\quad \text{for} \quad \alpha = [\alpha],\\
   C 6^{\alpha}     \frac{(\alpha - \alpha')^2}{ \left( ([ \alpha'] +1 -  \alpha')  (\alpha - [\alpha] )\right)^3} \quad \text{for} \quad  \alpha > [\alpha].
 \end{cases}
 \]
 
 Let us justify the constant in \eqref{lem:nik:emb:3al}. For the $i$-th intermediate pair the constant of  \eqref{lem:nik:emb:3a} reads
   \begin{equation}\label{lem:nik:emb:3alC}
\frac{C}{\beta_i^2 \delta}     \frac{6^{\alpha_{i+1} }}{(\alpha_{i+1} - [ \alpha_i]) \delta^{\alpha_{i+1}} }  \le  \frac{C}{\beta^2 \delta}    \frac{(\alpha_{K+1} - \alpha_0)^2}{(\alpha_{{i+1}} - \alpha_i)^2}  \frac{6^{\alpha}}{(\alpha_{i+1} - [ \alpha_i]) \delta^{\alpha} }  \le   \frac{C 6^{\alpha}}{\beta^2 \delta^{1+\alpha}}    \frac{(\alpha - \alpha')^2}{(\alpha_{{i+1}} - \alpha_i)^3}, 
  \end{equation}
where for the above inequality we have used definitions of $\beta_i, \,\alpha_0, \,\alpha_{K+1}$ and inequalities $\alpha_i \le \alpha$, $\alpha_{{i+1}} - \alpha_i \le \alpha_{{i+1}} - [\alpha_i]$. We get rid of the intermediary $\alpha_{i+1} $ in the r.h.s. of \eqref{lem:nik:emb:3alC}, observing that
  \[
\alpha_{{i+1}} - \alpha_i \ge  ([ \alpha'] +1 -  \alpha' ) \wedge ( \alpha - ([ \alpha'] +K))
  \]
Using \eqref{eq:Knat} in the above estimate and next employing the result in \eqref{lem:nik:emb:3alC}, we arrive at
   \begin{equation}\label{lem:nik:emb:3alC'}
\frac{C}{\beta_i^2 \delta}     \frac{6^{\alpha_{i+1} }}{(\alpha_{i+1} - [ \alpha_i]) \delta^{\alpha_{i+1}} }   \le  \frac{ C_{\eqref{lem:nik:emb:3al}} (\alpha, \alpha') }{\beta^2 \delta^{1+\alpha}}. 
  \end{equation}
We use $K+1$ times the constant of the r.h.s. of \eqref{lem:nik:emb:3alC'}. As $K+1 \le \alpha - \alpha' + 2$, we have justified the~form of the constant in \eqref{lem:nik:emb:3al}.

Using the estimate \eqref{lem:nik:emb:3al}, we must also keep in mind the condition \eqref{lem:nik:emb:236s}. For its validity 
   \begin{equation}\label{lem:nik:emb:3alC'2}
  \delta \le  \frac{|I|}{ 2^{ \left[ (\beta_i)^{-1} \right] + 3}} \le  \frac{|I|}{ 2^{ \left[ \frac{  C_{\eqref{lem:nik:emb:3alC'2}}(\alpha, \alpha')}{ \beta} \right] + 3}},
  \end{equation}
suffices,
where
\vspace{-2mm}
   \[
C_{\eqref{lem:nik:emb:3alC'2}} (\alpha, \alpha') =
           \begin{cases}
            \frac{\alpha - \alpha'}{[ \alpha'] +1 -  \alpha'} \quad\quad\quad\quad \text{for} \quad \alpha = [\alpha],\\
        \frac{\alpha - \alpha'}{ ([ \alpha'] +1 -  \alpha')  (\alpha - [\alpha] )} \quad \text{for} \quad  \alpha > [\alpha].
 \end{cases}
 \]
The second inequality in \eqref{lem:nik:emb:3alC'2} is valid in view of
\[
\frac{1}{ \beta_i} = \frac{1}{\beta}    \frac{(\alpha_{K+1} - \alpha_0)}{(\alpha_{{i+1}} - \alpha_i)} \le \tilde C(\alpha, \alpha') \frac{1}{\beta},
\]
which we justify analogously to \eqref{lem:nik:emb:3alC}.

\underline{Conclusion of Step 2.} We put together \eqref{lem:nik:emb:3a} and the condition \eqref{lem:nik:emb:236s} with \eqref{lem:nik:emb:3al} augmented with the condition \eqref{lem:nik:emb:3alC'2}. Hence, after a generous estimate of constants, we obtain part (ii) of the thesis in the case  $ \alpha >\alpha' $ and  $p \le q$. The rest follows from the H\"older inequality  (compare two first sentences of this proof of part (ii) of the thesis). Specifically, for the case $p > q$ we use first the $ \alpha >\alpha' $ and  $p = q$ estimate, which has $\beta =  \alpha - \alpha' $ and next the H\"older inequality. This is reflected in presence of the term $\beta \wedge (\alpha - \alpha')$ in the definition of the constant $C_\eqref{lem:nik:emb:th:low}$.

For the part (iii) of the thesis we quote Corollary 26 of \cite{Sim90}. The estimate there reads in our language
\[
 |f|_{C^{\alpha - \frac{1}{p}} ( I; X) } \le  |f|_{|I|, N^{\alpha, p} ( I; X) } \le  {3}{\delta^{-\alpha}}  |f|_{\delta, N^{\alpha, p} ( I; X) }. 
\]
One can also provide a higher-order version of this estimate, analogous to the one from the part (ii), but we do not need it here.
 \end{proof}

 \begin{proof}[Proof of Corollary \ref{cor65}]
 The equivalence of norms is standard. Therefore let us only observe that one direction is a verbatim of our reduction result. The other follows from the definition of a weak derivative, the~formula for the discrete integration by parts and a limit passage. For \eqref{eq:jb_strong:nik_bas:n:diff:both:delta:sob}, we repeat proof of \eqref{eq:jb_strong:nik_bas:n:diff:both:delta} of Lemma \ref{lem:jb_strong:nik_bas}. The Rademacher part of this corollary for the one-dimensional case follows from  from the absolute continuity of $W^{{1, 1}} ( I; X)$ functions and \eqref{lem:jb_weak:wtd:form}.
 \end{proof}



\end{document}